**UNIVERSIDAD NACIONAL AUTÓNOMA DE MÉXICO**
PROGRAMA DE MAESTRÍA Y DOCTORADO EN CIENCIAS MATEMÁTICAS Y DE LA ESPECIALIZACIÓN EN ESTADÍSTICA APLICADA

# DEFORMATIONS ON ENTIRE TRANSCENDENTAL FUNCTIONS

THESIS
TO OBTAIN THE DOCTORAL DEGREE
IN MATHEMATICAL SCIENCES

PRESENTS:
RODRIGO ROBLES MONTERO

PhD. GUILLERMO JAVIER FRANCISCO SIENRA LOERA
FACULTAD DE CIENCIAS, UNAM.

PhD. SANTIAGO LÓPEZ DE MEDRANO SÁNCHEZ
INSTITUTO DE MATEMÁTICAS, UNAM.
PhD. PATRICIA DOMÍNGUEZ SOTO
FACULTAD DE CIENCIAS FÍSICO MATEMÁTICAS, BUAP

MEXICO CITY, MARCH 2023.

*Al rostro sabio y corazón firme de nuestra universidad,*
*tan hermosa y tan fuerte,*
*aquí he crecido, aprendido y construido,*
*donde hemos peleado, resistido, perdido y ganado,*
*en ella he reído, llorado y amado,*
*mis mejores amigos los he conocido aquí,*
*a mis mejores compañeros y maestros también.*

*Es uno de los pocos lugares en nuestro país*
*que tiene un corazón que piensa*
*y un cerebro que siente.*
*No es perfecta, pero nosotros, su comunidad,*
*la hacemos posible,*
*la hacemos real.*
*Una gran parte de lo que soy se lo debo a ella.*

*"Por mi raza hablará el espíritu."*

*A mi universidad, a la UNAM.*

*"He aquí a tu madre,*
*tu señora, de su vientre,*
*de su seno te desprendiste,*
*brotaste.*
*Como si fueras una yerbita,*
*una plantita, así brotaste.*
*Como sale la hoja,*
*así creciste, floreciste.*
*Como si hubieras estado dormido*
*y hubieras despertado."*

*Para cuidarse, para protegerse,*
*un pedacito de la vida ha tejido un lazo*
*en este lugar de luz y oscuridad,*
*es uno de los más fuertes que existen,*
*el amor de una madre y su hijo.*

*A Graciela Montero, mi mamá.*

*“Aquí están mis amores,*
*mi pareja y mi hijita,*
*mi collar de piedras finas,*
*mi plumaje de quetzal,*
*mi hechura humana,*
*las nacidas de mí,*
*y yo,*
*el nacido de ellas.*
*Ustedes son mi sangre,*
*mi color,*
*en sus almas está mi imagen.*
*Y en mi alma las suyas.”*

*A Mary, mi solecito.*
*A Mixi, mi estrellita.*

*"Para mí*
*sólo recorrer los caminos que tienen corazón,*
*cualquier camino que tenga corazón.*
*Por ahí yo recorro,*
*y la única prueba que vale*
*es atravesar todo su largo.*
*Y por ahí yo recorro mirando,*
*mirando, sin aliento."*

*A don Juan.*
*A sus aprendices.*

*Ilhuikatl & Koskayolotl*

# Contents

# Symbols

| | |
|---|---|
| $\mathbb{C}$ | The set of complex numbers, the complex plane. |
| $\overline{\mathbb{C}}$ | $\mathbb{C} \cup \{\infty\}$, the Riemann sphere. |
| $\mathbb{C}^*$ | $\mathbb{C} \setminus \{0\}$, the punctured plane. |
| $diam(A)$ | Euclidean diameter of $A \subseteq \mathbb{C}$. |
| $diam_s(A)$ | Spherical diameter of $A \subseteq \overline{\mathbb{C}}$. |
| $f^n$ | The $n^{th}$ iteration of a function $f$. |
| $f_t \rightrightarrows f$ | The uniform convergence of $f_t$ to $f$. |
| $\mathcal{F}$ | A collection of functions. |
| $F(f), F$ | The Fatou set of a function $f$. |
| $J(f), J$ | The Julia set of a function $f$. |
| $\mathcal{L}$ | The grand orbit of $\Lambda$ under $f$. |
| $P(f)$ | The postsingular set of $f$ defined as $\overline{\bigcup_{n=0}^{\infty} f^n\left(Sing\left(f^{-1}\right)\right)}$. |
| $\Lambda$ | A Baker lamination in a cycle of Baker domains $U$. |
| $(T, X, F)$ | A dynamical system with $T$ a monoid, $X$ a non-empty set and $F : T \times X \to X$ a function satisfying $F(0, x) = x$ and $F(t_2, F(t_1, x)) = F(t_1 + t_2, x)$. |
| $U/f$ | The quotient space identifying points under the grand orbit relation of $f$ in $U$. |
| $\mathbb{Z}^+$ | The set of positive integers. |

# Figures

Figure 1.1: Taken from [*Baransky & Fagella 2001*].

Figure 1.2: Made with Fractastream.

Figure 1.3: Made with Fractastream.

Figure 1.4: Made with Fractastream.

Figure 1.5: Taken from [*Fagella & Henriksen 2006*].

Figure 1.6: Made with Fractastream and edited with Inkscape.

Figure 1.7: Made with Inkscape.

Figure 2.1: Taken from [*Branner & Fagella 2014*].

Figure 2.2: Taken from [*Branner & Fagella 2014*].

Figure 2.3: Taken from [*Branner & Fagella 2014*].

Figure 2.4: Taken from [*Branner & Fagella 2014*].

Figure 2.5: Taken from [*Branner & Fagella 2014*].

Figure 2.6: Taken from [*Branner & Fagella 2014*].

Figure 2.7: Taken from [*Branner & Fagella 2014*].

Figure 2.8: Taken from [*Branner & Fagella 2014*].

Figure 2.9: Taken from [*Robles & Sienra 2022*] and made with Inkscape.

Figure 2.10: Made with Inkscape.

Figure 2.11: Taken from [*Haissinsky & Tan 2004*] and edited with Inkscape.

Figure 2.12: Made with Inkscape.

Figure 2.13: Taken from [*Robles & Sienra 2022*] and made with Inkscape.

Figure 3.1: Taken from [*Robles & Sienra 2022*] and made with Fractastream.

Figure 3.2: Taken from [*Robles & Sienra 2022*] and made with Inkscape.

Figure 4.1: Made with Inkscape.

Figure 4.2: Made with Inkscape.

Figure 4.3: Made with Inkscape.

Figure 4.4: Made with Inkscape.

Figure 4.5: Made with Inkscape.

Figure 5.1: Made with Fractastream.

Figure 5.2: Made with Fractastream and edited with Inkscape.

Figure 5.3: Made with Inkscape.

# 1 Introduction

Movement involves a great portion of the events taking place in the universe in the stream of time, and understanding it has been important for us to develop ourselves technologically, culturally, and intellectually.

Understanding the cycles of days, seasons, and years, has been fundamental to carrying out the agricultural revolution in our prehistory. Studying the motion of bodies and their thermodynamic and electromagnetic interactions in the XVII and XVIII centuries, gave us the basis for the development of the industrial revolution. The discovery of how our bodies function was and is essential to improve our quality of life and its expectancy.

It is important to realize that our human activities impact the biosphere, that we have extinguished and keep extinguishing hundreds of species of living beings, see [*Wilson, 1999*] and [*Bernstein & Chivian 2008*], and we need to solve adequately this interaction in order to keep us living in this planet. And much of these themes are, until a certain point, related to dynamical systems.

One may ask what is a dynamical system? In culture, this question may have many solid answers depending on which discipline answers it. In this thesis, we will approach it from the mathematical point of view by means of the next definition. For this, we will use the concept of *monoid*, i.e., a set $T$ with an associative operation $+ : T \times T \to T$ that has a neutral element 0 in $T$.

DEFINITION 1.1 *We say that a triad* $(T, X, F)$ *is a* ***dynamical system*** *if* $T$ *is a monoid,* $X$ *is a non-empty set, and* $F$ *is a function* $F : T \times X \to X$ *satisfying* $F(0, x) = x$ *and* $F(t_2, F(t_1, x)) = F(t_1 + t_2, x)$.

Notice that in the previous definition, we can use a semigroup, which does not have a neutral element, or a group instead of a monoid depending on the desired properties to work on the specific dynamical system.

The function $F(t, x)$ is called the *evolution of the dynamical system* which can usu-ally be interpreted as the movement of the system, and $x$ can be seen as the *object in movement.* $F$ links to each $x \in X$ only one image depending on $t$, understood as the *time.* In this context, $X$ is called the *phase space*, which intuitively is the space and the moving objects, and when we take the pair $(0, x)$ we have the *initial state of the system on* $x$.

In this thesis we will study a particular dynamical system where $T = \mathbb{Z}^+ \cup \{0\}$, the positive integer numbers with zero, which can be interpreted as discrete time; with $X = \mathbb{C}, \overline{\mathbb{C}}$ or $\mathbb{C}^*$, that is to say, the space can be the plane, the Riemann sphere, or the punctured plane; and $F(n, z) = f^n(z)$, with $f$ a holomorphic function in $X$ and $f^n(z)$ the $n^{th}$ iteration of $f$, which will regularly be an entire transcendental function or a rational function. This dynamical system is known as *Holomorphic Dynamics in one complex variable.* But where does this system come from?

Historically, the origins of this dynamical system come from two works studying Newton´s Method, this method finds the roots of a function and is an iterative algorithm. The first one consists of two papers by Ernst Schroeder, *"Ueber unendlich viele Algorithmen zur Auflosung de Gleichungen"* and *"Ueber iterite Functionen"*, see [*Schroeder 1870*] and [*Schroeder 1871*] respectively.

The second one, is by Arthur Cayley, *"Applications of the Newton-Fourier Method to an Imaginary Root of an Equation"*, see [*Cayley 1879*]. Even though the first known appearance of this method was implemented by the Babylons to approximate the square root of a number $a$, the formal studies started by Schroeder and Cayley would lead to Julia´s and Fatou´s analysis. See [*Alexander 1994*].

In 1915, the French *Académie de Sciences* announces that the research subject for its *Grand Prix des Sciences Mathemátiques* will be the iteration of holomorphic functions, promoting a global analysis. The works from Gaston Julia, *"Memoire sur l´itération des fonctions rationnelles"*, see [*Julia 1918*], and from Pierre Fatou, *"Sur les équations fonctionnelles"*, see [*Fatou 1919*], would excel by their analysis and would become the foundations of Holomorphic Dynamics, see [*Alexander 1994*], based in the theory of Normal Families from Paul Montel, see [*Montel 1927*].

Incapable of utterly classifying the Fatou components and proving the existence of Siegel disks (done until the 1940´s, see [*Siegel 1942*]) on one hand, and in the absence to visualize what was happening on the other hand, the field suffered big inactivity until the 1980´s when it mainly revitalized with the works from Douady and Hubbard, see [*Douady & Hubbard 1984-85*], Sullivan (who proved the non-existency of wandering domains for rational functions in [*Sullivan 1985*]), Milnor, Thurston, Baker, Lyubich, and Eremenko among other people.

Some of these theoretical developments were based on Quasiconformal theory and Teichmuller spaces, developed in the middle of the $XX$ century. The second aspect was discovered due to the growing calculus and graphing massive power of the computers, where it could be seen that related to the chaotic behavior existed some sets of impressive beauty and complexity, the fractals, autosimilar forms resembling natural objects like clouds and mountains, see [*Mandelbrot 1982*].

In these holomorphic dynamics, there is a dichotomy between completely invariant sets that dynamically behaves quite differently. Usually, it is introduced by means of normal families as was mentioned before:

DEFINITION 1.2 *Let $\mathcal{F}$ be a collection of functions from a Riemann surface $S$ to a Riemann surface $T$ where every infinite sequence of functions of $\mathcal{F}$ contains a subsequence which converges locally in compact sets and uniformly to a function $f$, not necessarily in $\mathcal{F}$, then $\mathcal{F}$ is a* ***normal family****.*

One of the most important results in the theory of normal families is the Montel Theorem. This states that if the functions belonging to $\mathcal{F}$ are holomorphic and $T = \overline{\mathbb{C}} \setminus \{a, b, c\}$, i. e., $T$ is a hyperbolic surface, then $\mathcal{F}$ is a normal family. With these assumptions we can resume the discussion about dynamic dichotomy with the following definition:

DEFINITION 1.3 *Let $f : S \to S$ where $S = \mathbb{C}, \overline{\mathbb{C}}$ or $\mathbb{C}^*$, we define the domain of normality of the collection of iterates $\{f^n\}$ as the* ***Fatou set*** *$F(f)$ or simply $F$. Additionally, its complement $S \setminus F$ is defined as the* ***Julia set*** *$J(f)$ or simply $J$.*

Some properties of these sets, when $f$ is non-linear, are that the Fatou set $F$ is open. While the Julia set $J$ is closed, perfect, non-numerable and it is the closure of the repelling periodic points. $F(f) = F(f^n)$, and $J(f) = J(f^n)$. If $F$ is empty, then $J = S$. For entire transcendental functions, both sets are unbounded, and it is common to define that $\infty \in J$. Furthermore, the dynamical system $\{\mathbb{Z}^+ \cup \{0\}, J(f), f^n\}$ is chaotic in Devaney´s sense (that is $f$ is topologically transitive, periodic points are dense in $J(f)$, and $f$ has sensitive dependence on initial conditions), for a deeper explanation see [*Devaney 1989*], [*Bergweiler 1993*] and [*Morosawa et al. 1998*].

The classification of the periodic components of the Fatou set $F(f)$ is the following. Let $U$ be a connected component of $F$, then $f^n(U) \subseteq F$, where proper inclusion is possible. For example, $f(z) = \lambda e^z$ with $\lambda \in (0, 1/e)$ has an attracting real fixed point and an attractive basin $U$ containing 0 but $0 \notin f(U)$. See [*Devaney 2010*].

We say that $U$ is *preperiodic* if there exist integer numbers $p > q \geq 0$ such that $f^p(U) = f^q(U)$. If $q = 0$, we say then that $U$ is *p-periodic*. If the component is not preperiodic we say that $U$ is a *wandering domain*.

THEOREM 1.1 [*Fatou 1919*], [*Cremer 1932*] & [*Siegel 1942*] *Let $U$ be a p-periodic component of the Fatou set, therefore only one of the next conditions is possible:*

- *$U$ is an immediate basin of attraction of an attracting p-periodic point $z_0 \in U$, and $\lim_{n\to\infty} f^{np} = z_0$ for every $z \in U$.*
- *$U$ is a parabolic basin of a parabolic p-periodic point $z_0 \in \partial U$, and $\lim_{n\to\infty} f^{np} = z_0$ for every $z \in U$.*
- *$U$ is a Siegel disk where $U$ is biholomorphic to $\mathbb{D}$, and $f^p |_U$ is analytically conjugated to an irrational rotation $z \mapsto e^{i2\pi\theta} z$, with $\theta \in \mathbb{R} \setminus \mathbb{Q}$, of the disk $\mathbb{D}$.*
- *$U$ is a Herman ring where $U$ is biholomorphic to $\mathbb{A}_r := \{z : 1 < |z| < r\}$, and $f^p |_U$ is analytically conjugated to an irrational rotation $z \mapsto e^{i2\pi\theta} z$, with $\theta \in \mathbb{R} \setminus \mathbb{Q}$, of the annulus $\mathbb{A}_r$.*
- *$U$ is a Baker domain where $\lim_{n\to\infty} f^{np} = z_0$ for every $z \in U$ and $f$ is not well defined in $z_0 \in \partial U$.*

For rational functions, there are no Baker domains, and for entire transcendental functions, the Baker domain is only possible when $z_0 = \infty$. Also, there are no Herman rings for entire functions, see [*Bergweiler 1993*].

In this thesis, we will see a special kind of Baker domain, the univalent one, because they are a fountain of examples that will be useful in Chapter 5. Some properties of Baker domains are the following. Let $Sing\,(f^{-1})$ be the *set of singularities of the inverse function* $f^{-1}$, i.e., the closure of critical and finite asymptotic values of $f$. Let $P(f)$ be the postsingular set of $f$ defined as

$$P(f) := \overline{\bigcup_{n=0}^{\infty} f^n\,(Sing\,(f^{-1}))}.$$

In [*Eremenko & Lyubich 1992*] it is proved that if $f$ has a Baker domain, then $Sing\,(f^{-1})$ is unbounded. Moreover, a Baker domain for entire transcendental functions is simply connected, see [*Baker 1975*]. The inverse image of a simply connected non-invariant domain $U$ of an entire transcendental function omitting a value is disconnected, particularly if $U$ is a Baker domain, see [*Bergweiler & Eremenko 2007*]. Furthermore, [*Stallard 1990*] shows that if $d\left(P(f), J(f)\right) > 0$ (i.e., $f$ is hyperbolic) then for every $z \in J(f)$, we have $\left|(f^n)'(z)\right| \to \infty$ as $n \to \infty$, i.e., $f$ is expansive in its Julia set.

On the other hand, if $U$ is a Baker domain of an entire transcendental function $f(z)$ since it is simply connected, by the Riemann Mapping Theorem there exists a biholomorphism $\psi : \mathbb{H} \to U$. Let $g = \psi^{-1} \circ f \circ \psi$ where $\underset{n\to\infty}{lim}\, g^n(w) = \infty$ for every $w \in \mathbb{H}$. Since $\psi$ is univalent too, $g \in PSL(2, \mathbb{R})$, and it is a polynomial conjugated to only two possible mappings:

$$g(w) = \begin{cases} aw \ (a > 1) & \text{Hyperbolic type,} \\ w + 1 & \text{Parabolic type.} \end{cases}$$

See Figure 1.1, where the lines and the arcs of circles are invariants, and the arrows show the direction of the dynamics under iteration, here $h(u)$ is the corresponding mapping to $g(w)$ onto the Poincaré disk $\mathbb{D}$.

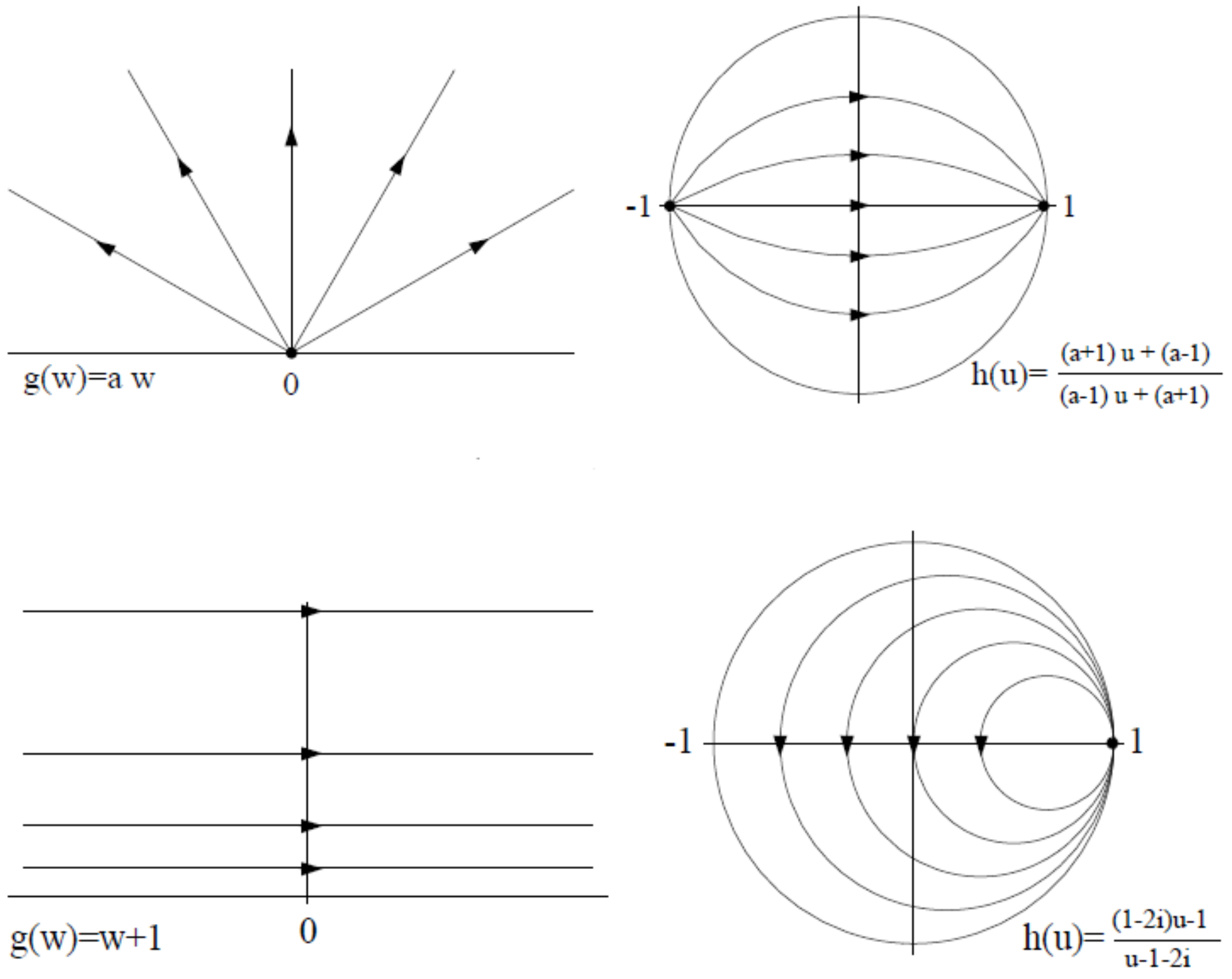

**Figure 1.1:** $g(w)$ on $\mathbb{H}$ and $h(u)$ on $\mathbb{D}$. Taken from [*Baransky & Fagella 2001*].

This motivates us to search for a classification of the univalent Baker domains for entire transcendental functions. So, we will need the next definition used in [*Baransky & Fagella 2001*]:

DEFINITION 1.4 *A point* $\zeta \in \overline{\mathbb{C}}$ *in the boundary of a simply connected domain* $U \subseteq \mathbb{C}$ *is* ***accessible from*** $U$ *if there is a curve* $\gamma : [0,\infty) \to U$ *landing at* $\zeta$*, i.e.,* $\lim_{t\to\infty} \gamma(t) = \zeta$*. And we say that two curves* $\gamma_1$ *and* $\gamma_2$ *have the* ***same access to*** $\zeta$ *if for every neighborhood* $V \subseteq \overline{\mathbb{C}}$ *of* $\zeta$ *there exists a curve* $\alpha : [0,1] \to U \cap V$ *such that* $\alpha(0) \in \gamma_1$ *and* $\alpha(1) \in \gamma_2$*.*

Let us note that if $\partial U$ is locally connected, thus all its boundary is accessible and equally, an access is a homotopy class within the collection of curves $\tilde{\gamma} : [0,1] \to \overline{\mathbb{C}}$, such that $\tilde{\gamma}((0,1)) \subset U$ and $\tilde{\gamma}(1) = \zeta$. Furthermore, [*Baker 1988*] demonstrates that $\infty$ is accessible from every Baker domain $U$ of an entire function. This implies that the iterates of every point in $U$ tend to $\infty$ through the same access, see Lemma A in [*Baransky & Fagella 2001*].

THEOREM 1.2 [*Baransky & Fagella 2001*] *Let* $f : \mathbb{C} \to \mathbb{C}$ *be an entire transcendental function and let* $U \subseteq \mathbb{C}$ *be an invariant and univalent Baker domain. Then, there is a point* $\zeta \in \overline{\mathbb{C}}$*, such that the backward iterates under* $(f|_U)^{-1}$ *of every point in* $U$ *tend to* $\zeta$*, which is attracting or parabolic, through the same access (which we call the* backward dynamical access *). Moreover, exactly one of the following occurs:*

- $U$ *is* HYPERBOLIC TYPE I: $\zeta \neq \infty$ *is an attracting or parabolic fixed point with multiplier* 1 *in* $\partial U$*. (See Figure 1.2)*
- $U$ *is* HYPERBOLIC TYPE II: $\zeta = \infty$*, the backward dynamical access is different from the forward one. In this case* $\partial U$ *is disconnected. (See Figure 1.3)*
- $U$ *is* PARABOLIC TYPE: $\zeta = \infty$*, the backward dynamical access is equal to the forward one. (See Figure 1.4)*

[*Baransky & Fagella 2001*] prove that for $U$ hyperbolic, the Riemann mapping $\psi$ has non-tangential limits in 0 (equivalent to $\zeta$) either in $\infty$(equivalent to $\infty$).

Furthermore, if $U$ is hyperbolic type and $\exists \lim_{w\to\varsigma} \psi(w)$, with $\varsigma \in \{0,\infty\}$ then $\zeta$ is the only periodic point in $\partial U$ for type I and there are no periodic points for type II. If $U$ is parabolic and $\exists \lim_{w\to\infty} \psi(w)$, there are no periodic points in $\partial U$.

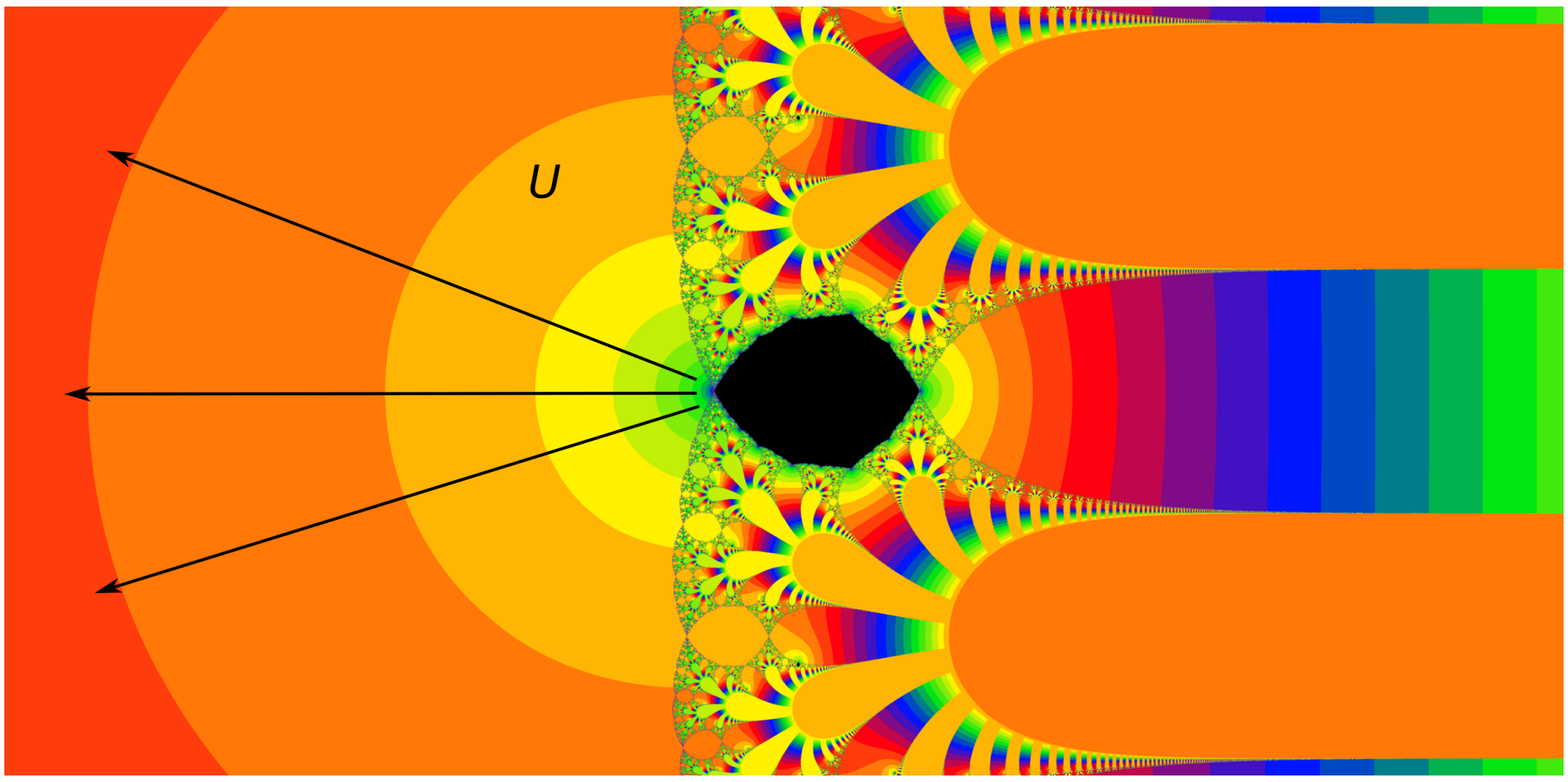

**Figure 1.2:** Hyperbolic Type I with Baker domain $U$, $f(z) = 2 - log(2) + 2z - e^z$.

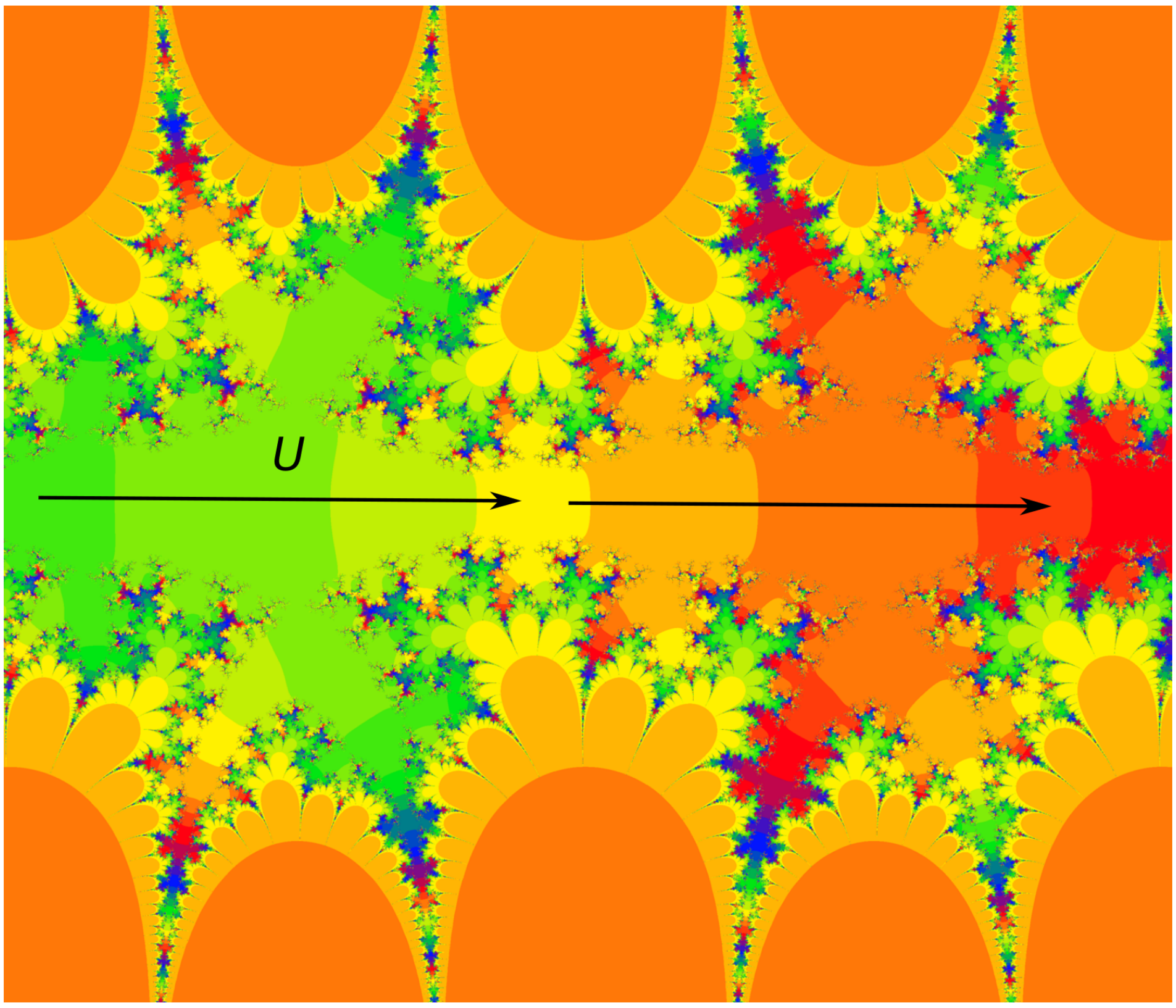

**Figure 1.3:** Hyperbolic Type II with Baker domain $U$, $f(z) = z + 1.8 + 0.6sin(z)$.

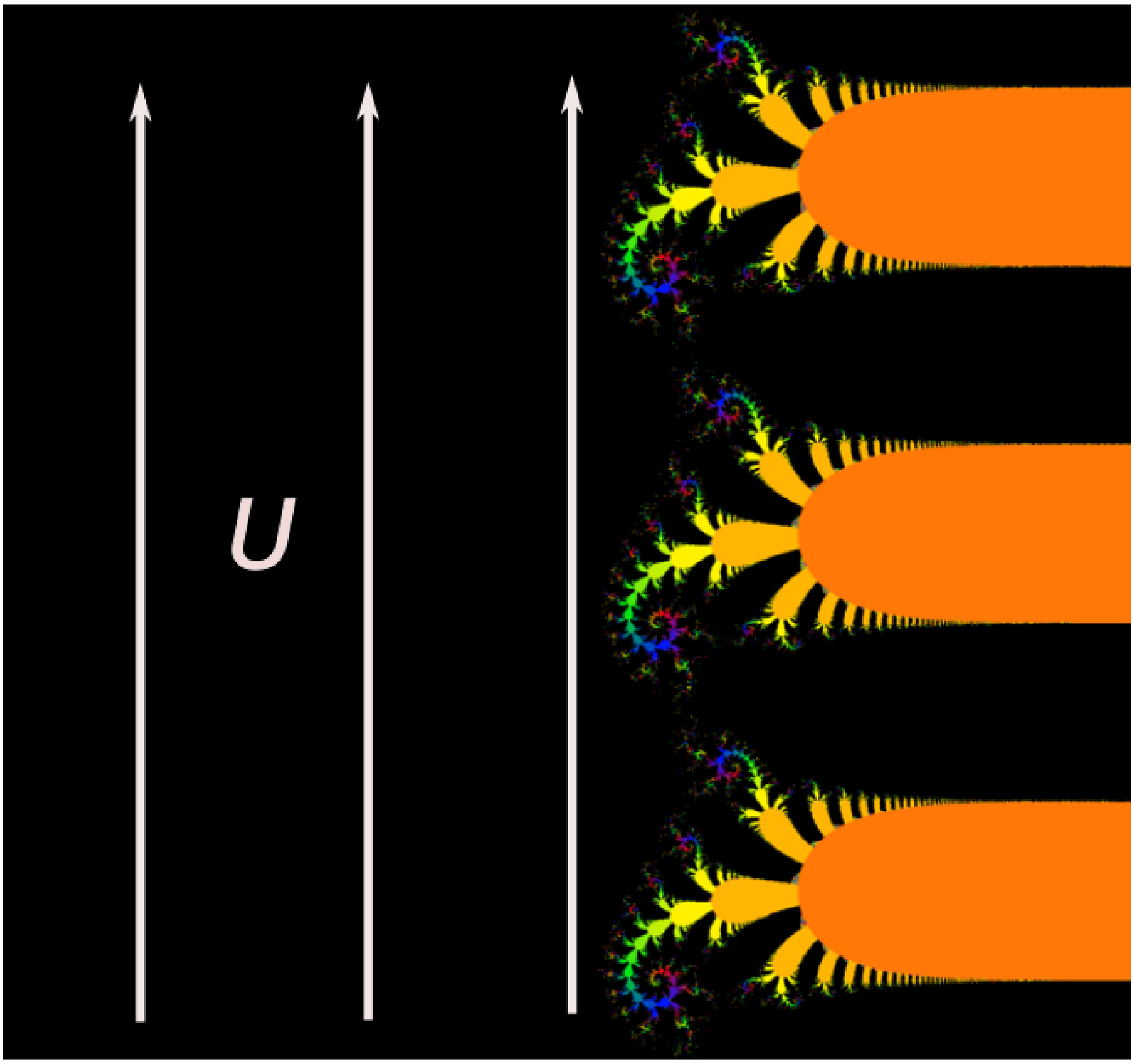

**Figure 1.4:** Parabolic Type with Baker domain $U$, $f(z) = z + log\left(\frac{\sqrt{5}-1}{2}\right) + e^z$.

Since we will be working in the boundary of Baker domains, it is important to know that for the non-univalent ones, [*Baker & Weinrich 1991*] showed $\partial U$ is not a Jordan curve and [*Baker & Domínguez 1999*] proved the infinite existence of different accesses to infinity.

Additionally, we define the *grand orbit* of $y \in U$ as the set $\{x \in U \mid f^n(x) = f^m(y)$ for some $n, m > 0\}$. The *grand orbit of a set* $A$ is the union of the grand orbits of the elements in $A$. The *grand orbit relation,* $x \sim y$, exists if and only if $x$ and $y$ have the same grand orbit, and it is an equivalence relation. We denote by $U/f$ the quotient space resulting from identifying points under the grand orbit relation of $f$ in $U$. In this context, we have the following result:

THEOREM 1.3 [*Fagella & Henriksen 2006*] *Let $f$ be an entire transcendental function with a Baker domain $U$. Therefore $U/f$ is a Riemann surface conformally isomorphic to only one of the following cylinders (See Figure 1.5):*

- $\{-s < Im(z) < s\}/\mathbb{Z}$ *for some $s > 0$ and we say that $U$ is* hyperbolic.
- $\{Im(z) > 0\}/\mathbb{Z}$ *and we say that $U$ is* simply parabolic.
- $\mathbb{C}/\mathbb{Z}$ *and we say that $U$ is* doubly parabolic. *In this case, $f : U \to U$ is not proper or has a degree greater than 1.*

In [*Fagella & Henriksen 2006*] and [*Bergweiler & Zheng 2012*] it is proved the relation among these surfaces with the univalency of the Baker domains. If $f$ is univalent the associated surface is hyperbolic or simply parabolic, and if $f$ is not univalent it can be of any type.

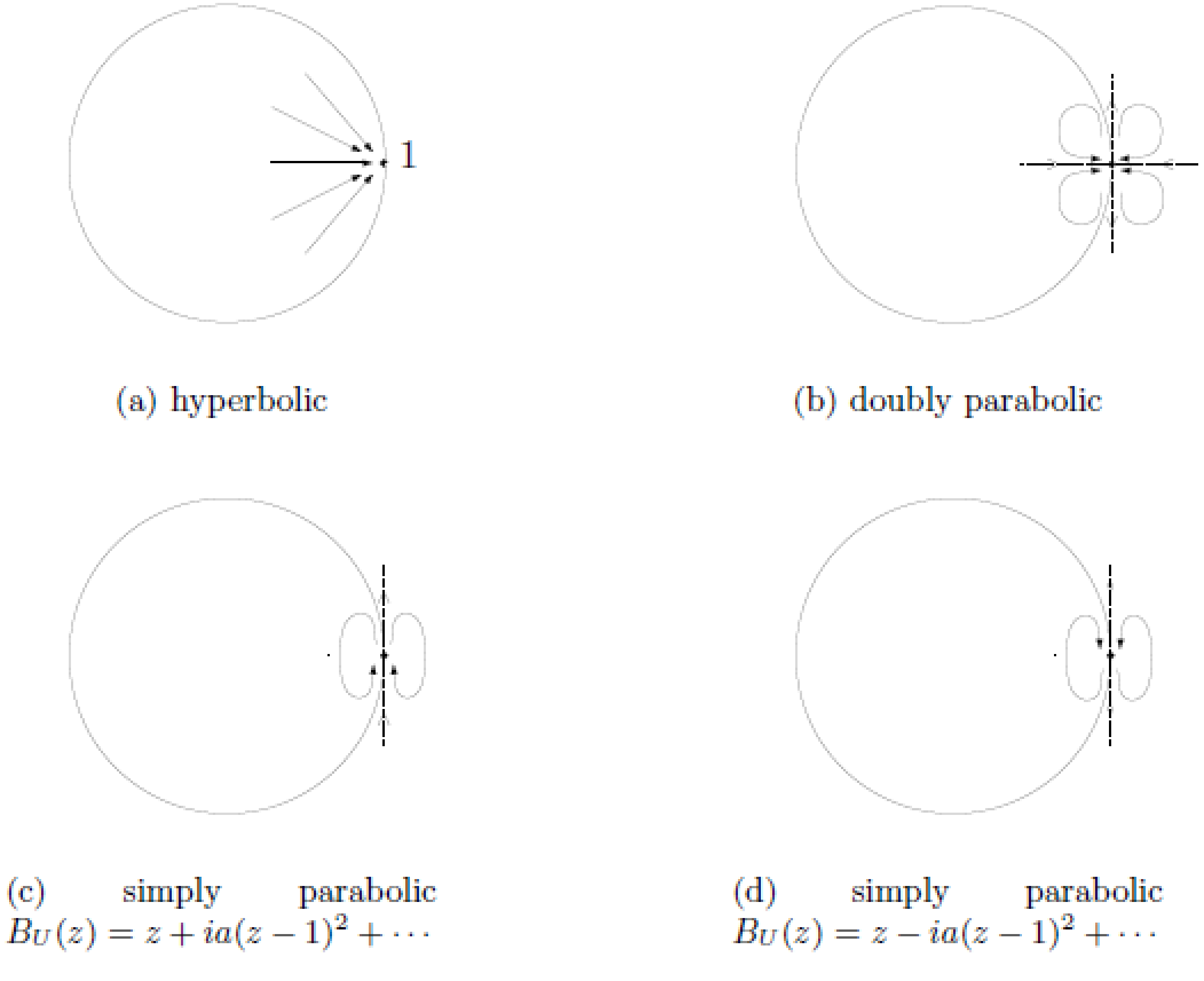

**Figure 1.5:** The orbits of $\mathbb{D}$ tending to 1, conjugated to the ones tending to $\infty$ in $U$. Taken from [*Fagella & Henriksen 2006*].

In this frame, we generalize the concept of *pinching deformation.* This was introduced by [*Makienko 2000*] as a tool to prove that a specific kind of *component of J-stability* is unbounded (i.e., the set of functions $g$ in $\mathbb{CP}^{2d+1}$ such that given a rational function $f$ with degree $d$ there exists a quasiconformal mapping conjugating $g$ with $f$ up to a Moebius transformation in given neighborhoods about $J(g)$ and $J(f)$).

The papers [*Haissinsky 2002*], [*Tan 2002*], and [*Haissinsky & Tan 2004*] redefined the concept of pinching deformation. This concept consists of the deformation of a rational function $f$ with an invariant periodic curve among a repulsor periodic point and attractor periodic points using quasiconformal conjugacies (given by the Ahlfors-Bers-Morrey theorem) and becomes in the limit into a function $F$ with a parabolic point.

In this work we develop another generalization of pinching by way of laminations:

DEFINITION [*Robles & Sienra 2022*] *Let* $f : \mathbb{C} \to \mathbb{C}$ *be an entire transcendental function with a Baker domain* $U$ *such that* $f(U) = U$ *provided with a hyperbolic metric, let* $\Lambda$ *be a set of complete geodesics in* $U$*. We say that* $\Lambda$ *is a* ***Baker lamination*** of $U$*, if the geodesics* $\lambda \in \Lambda$*, called* ***leaves*** *henceforth, satisfy:*

1. *The leaves of the lamination do not accumulate in* $U$*.*
2. *If* $\lambda \in \Lambda$ *then* $f^n(\lambda) \in \Lambda$*, with* $n \in \mathbb{N}$*. Also,* $\lambda \subset U$ *is in* $\Lambda$*, if* $f^n(\lambda) \in \Lambda$*, for some* $n > 0$*.*
3. *For any different leaves* $\lambda, \lambda' \in \Lambda$*,* $\lambda \cap \lambda' = \emptyset$*.*
4. *For any* $\lambda \in \Lambda$, *there exist* $\partial\lambda := \lim_{t\to\pm\infty} \lambda(t)$ *and* $\partial\lambda \subset \partial U \subset \overline{\mathbb{C}}$*.*

We are interested in the leaves $\lambda$ with $\infty$ as an endpoint of $\partial\lambda$, the leaves containing $\infty$ and a point $a \in \mathbb{C}$ as endpoints of $\partial\lambda$ are denoted by $\lambda_{a,\infty}$. Also, let $\mathcal{L} := \bigcup_{k\in\mathbb{N}} f^{-k}(\Lambda)$, i.e., the grand orbit of $\Lambda$.

Is in this lamination that we are going to perform the pinching deformation. The construction of the pinching deformation requires some nontrivial mathematical development, for this reason, we write down some related concepts and its definition, but further details are discussed in Chapter 2. The *almost complex structures* $\sigma_t$ are the union of the *linear conformal structures* $\sigma(u)$ which are vector spaces $\mathbb{C}(\mathbb{C}, +, \star)$ related to linear maps $L_u$, mapping certain ellipses (real scaled) on the tangent space $T_u$ to circles on $u \in \overline{\mathbb{C}}$. And $\mathcal{V} := \bigcup_{k\in\mathbb{N}} f^{-k}\left(V_\delta\left(\Lambda\right)\right)$ is the grand orbit of disjoint *good neighborhoods (of thickness* $\delta$*) for the leaves* $\lambda \in \Lambda$*,* $V_\delta\left(\Lambda\right)$*.*

DEFINITION [*Robles & Sienra 2022*] *Let $f$ be an entire transcendental function with at least one cycle of periodic Baker domains $U = \{U_0, U_1, \ldots, U_{p-1}\}$ with a Baker lamination $\Lambda$ in $U$. The family of almost complex structures $(\sigma_t)_{t\in[0,1)}$ defines a* ***pinching deformation of*** *$f$ with support in $\mathcal{V}$. These structures come with quasiconformal maps $h_t : \overline{\mathbb{C}} \to \overline{\mathbb{C}}$ via integration by the* Ahlfors-Bers-Morrey *Theorem, that we can normalize assuming $h_t$ fixes $\infty$ and two points $p, q \in J(f)$. In addition, the* function *$f_t := h_t \circ f \circ h_t^{-1}$ is holomorphic for every $t \in [0, 1)$.*

*We say that a* ***pinching deformation converges uniformly*** *if $h_t \rightrightarrows H$ (i.e., $h_t$ converges uniformly on the euclidean metric to a function $H$) and the nontrivial fibers of $H$ is the grand orbit $\mathcal{L}$. In the sense that $diam_s (h_t (\bar{\gamma})) \to 0$ as $t \to 1$, for each $\gamma \in \mathcal{L}$, where $diam_s(A)$ is the spherical diameter of a set $A \subset \overline{\mathbb{C}}$.*

An example of convergent pinching along $\lambda_{a,\infty}$ in non-univalent invariant Baker domains is presented in [*Domínguez & Sienra 2015*], deforming the Fatou function $f(z) = z + 1 + e^{-z}$ to $F_{p/q}(z) = z + e^{-z} + 2\pi ip/q$. This process deforms a doubly parabolic completely invariant domain into infinitely many doubly periodic domains or a wandering domain. In contrast to this example, we proved in [*Robles & Sienra 2022*] that the pinching along $\lambda_{a,\infty}$ in non-completely invariant Baker domains does not converge.

THEOREM A [*Robles & Sienra 2022*] *Let $f : \mathbb{C} \to \mathbb{C}$ be an entire transcendental function with a non-completely invariant Baker domain $U$. Consider a Baker lamination $\Lambda$ in $U$ with a leaf of type $\lambda_{z_0,\infty}$ having endpoints at $z_0 \in \mathbb{C}$ and $\infty$, with $z_0$ a non-exceptional point in $\partial U$. Thus, the pinching deformation along the grand orbit $\mathcal{L}$ does not converge uniformly.*

Hence, in this particular case, the boundary of the quasiconformal deformations is incomplete.

As an immediate consequence of this theorem, we have the next corollary for the possibility of asymptotic values:

COROLLARY A [*Robles & Sienra 2022*] *Let $f : \mathbb{C} \to \mathbb{C}$ be an entire transcendental function with a non-completely invariant Baker domain $U$. Consider a Baker lamination $\Lambda$, with a leaf $\lambda_{a,b}$ having endpoints at non-exceptional points $a, b \in \mathbb{C}$. If $\lambda_{a,b}$ intersects the set of asymptotic values of $f$, then the pinching deformation along the grand orbit of the lamination, $\mathcal{L}$, does not converge uniformly.*

On the other hand, we have the problem of finding sufficient conditions so that the pinching deformation is uniformly convergent. In this sense, [*Haissinsky & Tan*

*2004*] proved a theorem for rational functions. In this paper, Hassinsky and Tan Lei showed some of the richness of the pinching deformation as a mathematical tool. They gave some generalizations of theorems by Rees, Tan, and Shishikura for geometrically finite polynomials, see [*Rees 1986*], [*Tan 1990*], and [*Shishikura 2000*]: two Poscritically finite quadratic polynomials $f_c$ and $f_{c'}$ are matable if and only if $c$ and $c'$ do not belong to conjugate limbs of the Mandelbrot set.

Based on [*Haissinsky & Tan 2004*], in this thesis, we develop analogous results for uniform convergent pinching deformations of entire transcendental functions on laminations in their Baker domains. Some of the results of Haissinsky and Tan Lei are already valid in the case of entire transcendental functions and the rest are proved here.

As we will see in the proof of Theorem C, we need a theorem on rigidity for entire transcendental functions proved in [*Skorulski & Urbanski 2012*]. Equally, Theorem C needs the concept of semihyperbolicity and Theorem B.

DEFINITION *An entire transcendental function $f$ is* ***semihyperbolic at*** $a \in J(f)$*, if there exist $r > 0$ and $N \in \mathbb{N}$ such that for all $n \in \mathbb{N}$ and for all components $V$ of $f^{-n}(D_r(a)) := \{z \in \mathbb{C} \mid f^n(z) \in D_r(a)\}$ the function $f^n : V \to D_r(a)$ is a proper* function *of degree at most $N$, where $D_r(a) := \{z \in \mathbb{C} : |z - a| < r\}$. If this is satisfied for every $a \in J(f)$, we say that $f$ is semihyperbolic.*

We note that the concept of semihyperbolicity is introduced to entire transcendental functions as weak hyperbolicity is introduced to rational functions in the following sense. Semihyperbolicity will be used in Section 2.2 to guarantee some control on the inverse images of neighborhoods about the Julia set to shrink the inverse images of the lamination, id est, we need a version of Theorem II of [*Mañe 1993*] for entire transcendental functions, which was proved in [*Bergweiler & Morosawa, 2002*]:

THEOREM [*Bergweiler & Morosawa 2002*] *Let $f$ be an entire transcendental function and suppose that $f$ is semihyperbolic at $a \in J(f)$. Then there exists $r > 0$ with the following property: for all $\varepsilon > 0$, there exists $M \in \mathbb{N}$ such that if $n \geq M$ and $V$ is a component of $f^{-n}(D_r(a))$, then diam$(V) < \varepsilon$.*

The structure to prove Theorem C is the same as in [*Haissinsky & Tan 2004*], the same division by lemmas and propositions, but the order is different, is straightforward and is separated by sections. We need the next theorem to prove it:

THEOREM B *Let $f : \mathbb{C} \to \mathbb{C}$ be a semihyperbolic entire transcendental function with a Baker domain $U$. Let $\mathcal{L}$ be the grand orbit of a Baker lamination $\Lambda$ in $U$ that* *does not contain a leaf of the type* $\lambda_{a,\infty}$ *and sing* $(f^{-1}) \cap \mathcal{L} = \emptyset$*, then the maps* $\{h_t\}$ *that integrate the family of almost complex structures* $(\sigma_t)_{t\in[0,1)}$ *are equicontinuous in* $\overline{\mathbb{C}}$*. Furthermore, for any subsequence* $\{h_{t_k}\}$ *converging to a map* $H$*, the nontrivial fibers of* $H$ *are exactly the components of* $\mathcal{L}$*.*

To avoid the possibility of closed curves in $\mathcal{L}$ we ask that $f$ satisfies property $P$, described in Definition 2.3. In addition, we will need that the Lebesgue measure of $J(f)$ to be zero, thus we will ask that $J(f)$ be thin at $\infty$, see Section 4.8 for further detail.

THEOREM C *Let* $f : \mathbb{C} \to \mathbb{C}$ *be a semihyperbolic entire transcendental function with a Baker domain* $U$ *that satisfies property* $P$ *and with* $J(f)$ *thin at* $\infty$*. Let* $\mathcal{L}$ *be a grand orbit of a Baker lamination* $\Lambda$ *in* $U$ *that does not contain a leaf of the type* $\lambda_{a,\infty}$ *and sing* $(f^{-1}) \cap \mathcal{L} = \emptyset$*, then there exists a uniformly convergent continuous pinching deformation* $f_t = h_t \circ f \circ h_t^{-1}$ *to an entire function* $F$*. The mappings* $h_t$ *are quasiconformal mappings that converge uniformly to a map* $H$*, whose non-trivial fibers are the* $\mathcal{L}-$*components.*

As a consequence, Theorem C solves the following open problem. There are known examples of wandering domains with or without postsingular set within, i.e., the forward iterations of the critical points of $f$, but the existence of a wandering domain $W$ where the distance between the Poscritical set and $W$ is positive is unknown. The function $f(z) = 2 - log(2) + 2z - exp(z)$ is studied in [*Bergweiler 1995*], which has a Poscritical set $P(f)$ in the interior of a wandering domain $W$, whose distance to the Baker domain $U$ is positive. To this function, we will apply Theorem C for an adequate lamination $\mathcal{L}$ in the fundamental domain of $f$ in $U$, and we will obtain a positive answer to the problem. See Figure 1.6 and Figure 1.7. For further details see Chapter 5.

THEOREM D *There exists an entire function* $F$ *with a wandering domain* $W$ *such that* $d\,(P(F), W) > 0$*.*

It is worth mentioning that Theorem D is in accordance with Corollary C by [*Baransky et al., 2020*] stating:

COROLLARY C [*Baransky et al., 2020*]

*Let* $f$ *be a topologically hyperbolic meromorphic map* (i.e., $dist(P(f), J(f) \cap \mathbb{C}) > 0$) *and* $U$ *be a Fatou component of* $f$*. Denote by* $U_n$ *the Fatou component such that* $f^n(U) \subset U_n$ *and suppose that* $U_n \cap P(f) = \emptyset$ *for* $n > 0$*. Then, for every compact* *set $K \subset U$, every $z \in K$ and every $r > 0$ there exists $n_0$ such that for all $n \geq n_0$, $D_r(f^n(z)) \subset U_n$. In particular, diam$(U_n) \to \infty$ and dist$(f^n(z), \partial U_n) \to \infty$, for every $z \in U$, as $n \to \infty$.*

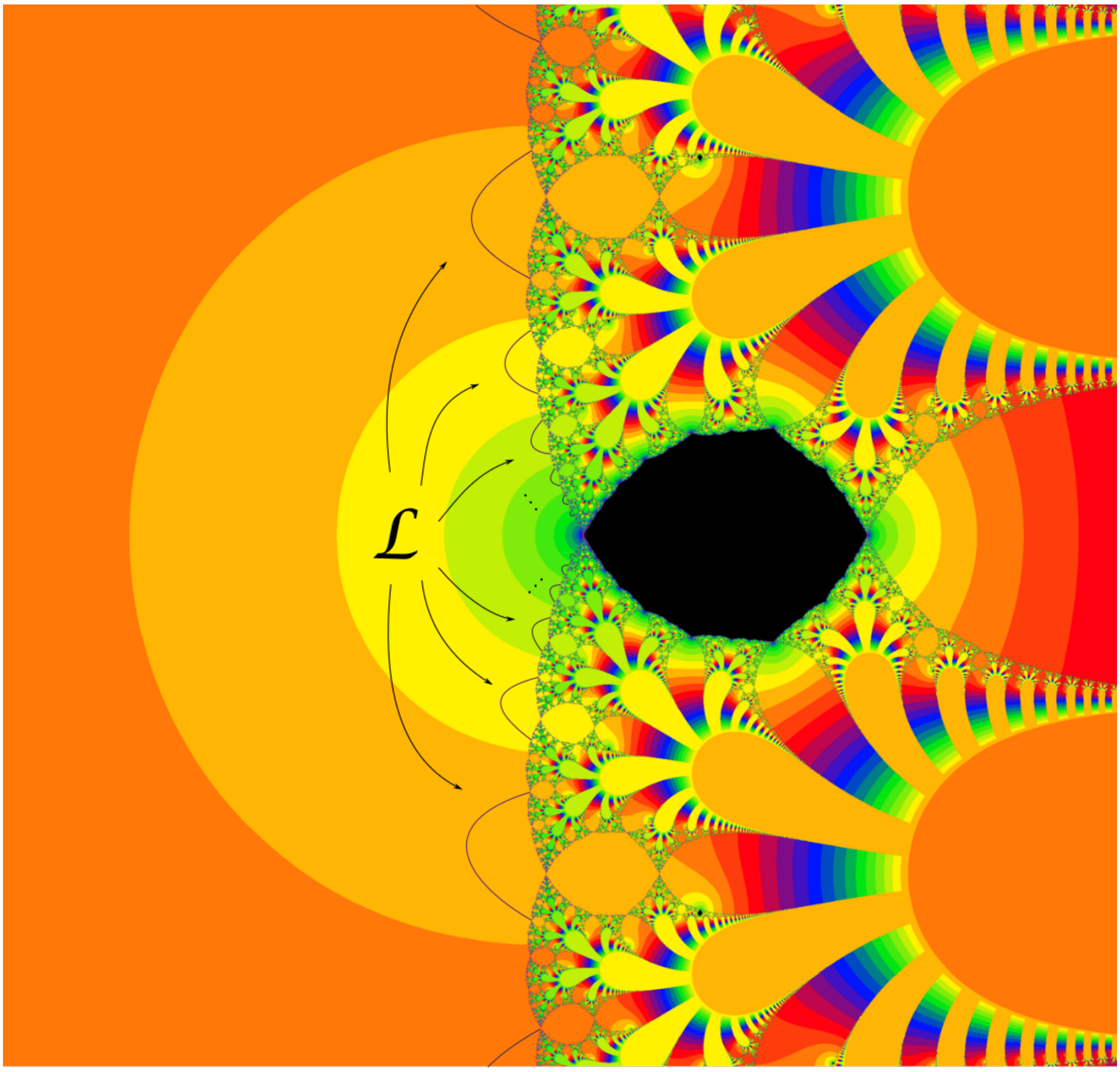

**Figure 1.6:** Pinching deformation of $f(z)$ along $\mathcal{L}$ for Theorem D.

In particular, the wandering domain $W_1$ of $F$, resulting from the pinching deformation, satisfies this corollary, because the pinching deformation only deforms the grand orbit $\mathcal{L}$ to point components, a neighborhood $\mathcal{Y} \setminus \mathcal{L}$ around the lamination is bounded deformed, and the complement $\mathbb{C} \setminus \mathcal{Y}^*$ is not deformed (see section 2.3 for definitions). This implies that the Baker domain $U_1$ of $F$ is practically the same as the Baker domain $U$ for $f$, which both are hyperbolic type I. Thus their fundamental domains grow conjugated to $z \mapsto az$, where in the case of $f$, $a = 2$. As the leaves of the lamination $\mathcal{L}$ are in the fundamental domains of $U$, when the pinching deformation is realized the wandering domain $W_1$ inherits this growth, satisfying corollary C from [*Baransky et al., 2020*]. See Figure 1.7.

To the knowledge of the author, the known examples of wandering domains do not grow without limit. It is relevant to ask if the return of this Corollary C by [*Baransky et al., 2020*] is valid. In the sense that the set of singularities of the inverse function, $sing\,(f^{-1})$, and the postsingular set $P(f)$ are related with many considerable results in holomorphic dynamics, for instance, in any attracting or parabolic basin, there exists a critical point.

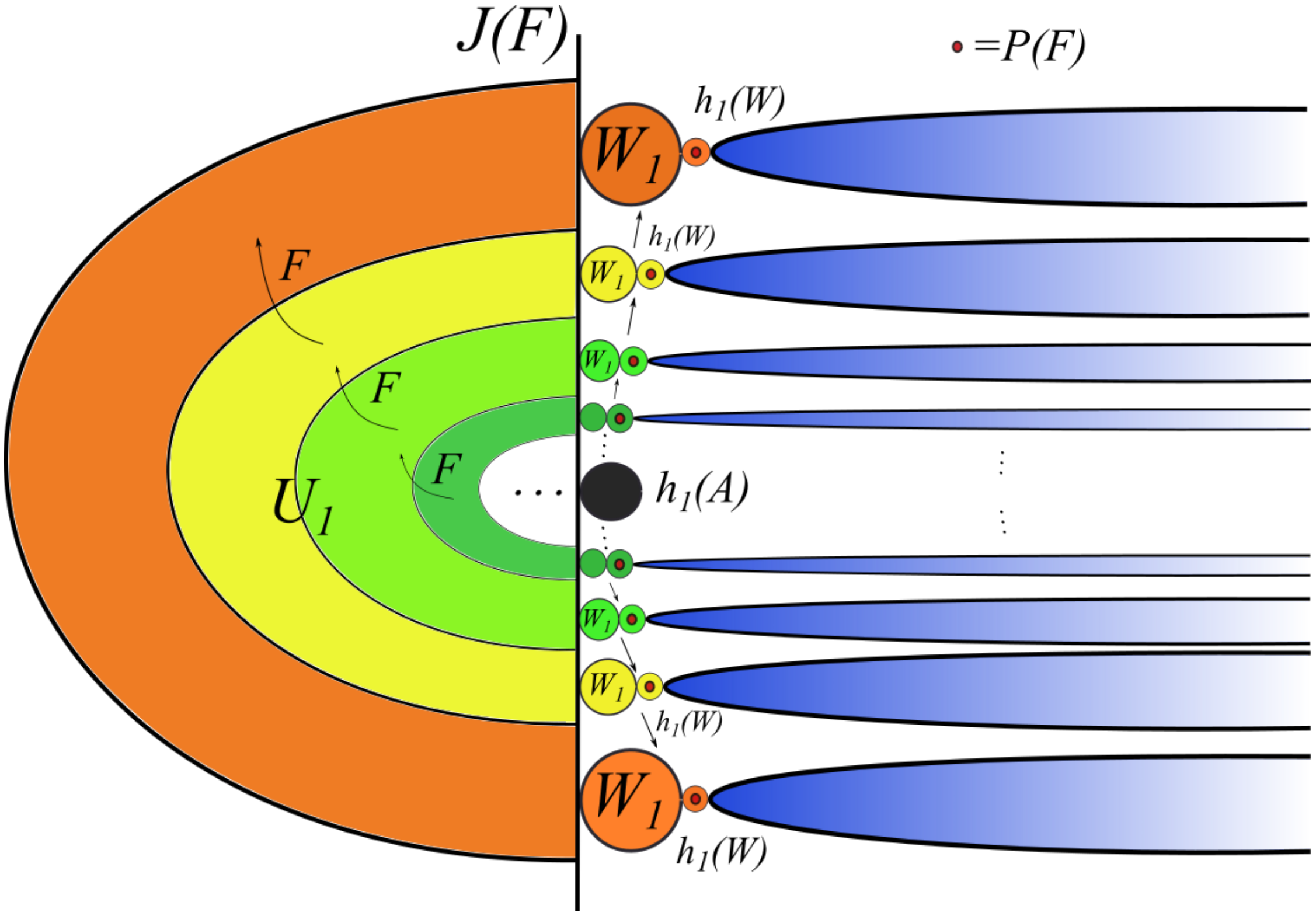

**Figure 1.7:** The wandering domain $W_1$ lies at a positive distance of $P(F)$ in red. With Baker domain $U_1$ and its preimages in blue.

# 2 Preliminaries

## 2.1 Quasiconformal Theory

In this section, we develop basic concepts of quasiconformal theory and its fundamental theorem. These will be essential for the realization of the pinching deformation. For further study see [*Ahlfors 1966*], [*Lehto 1987*], [*Gardiner 2000*], [*Zakeri & Zeinalian 1996*]. In particular, for holomorphic dynamics, we follow [*Branner & Fagella 2014*] approach closely.

There are different ways to define a quasiconformal map. Some of them are more intuitive, using modules of quadrilaterals or rings, or more analytical, using distributional derivatives. Here, we give a midway definition and mention the geometric one. We say that $f : I \to \mathbb{R}$ is *absolutely continuous on* $I$ if for every $\varepsilon > 0$, there exists $\delta > 0$ such that for finite intervals $(x_k, y_k) \subseteq I$, satisfying $\sum_k |x_k - y_k| < \delta$, this implies $\sum_k |f(x_k) - f(y_k)| < \varepsilon$. As a consequence, we have that absolute continuity implies uniform continuity. The next set of relations gives a reasonable idea where these definitions stand analytically:

$$\{\text{Continuously differentiable}\} \subseteq \{\text{Lipschitz continuous}\} \subseteq$$

$$\subseteq \{\text{Absolutely continuous}\} \subseteq \{\text{Bounded Variation}\} \subseteq$$

$$\subseteq \{\text{Differentiable almost everywhere}\}.$$

As an example, the square root is absolutely continuous but it is not Lipschitz continuous, and the devil´s staircase is not absolutely continuous but it is of bounded variation.

Now, a continuous real-valued function $u$ is said to be *absolutely continuous on lines* (ACL) in a domain $U \subseteq \mathbb{C}$ if for each closed rectangle $\{x + iy \mid a \leq x \leq b,\ c \leq y \leq d\} \subset U$, the function $x \mapsto u(x + iy)$ is absolutely continuous in $[a, b]$ for almost all $y \in [c,d]$ and the function $y \mapsto u(x+iy)$ is absolutely continuous in $[c,d]$ for almost all $x \in [a,b]$. A complex function is *absolutely continuous in $U$* (ACL) if its real and imaginary parts are ACL in $U$.

For example, let $f : \mathbb{R} \to \mathbb{R}$ be defined as

$$f(x) = \begin{cases} 0 & x < 0, \\ \text{devil´s staircase} & x \in [0,1]\,, \\ 1 & x > 1. \end{cases}$$

Let $F : \mathbb{C} \to \mathbb{C}$ be defined as $F(x+iy) = x + i\,(y + f(x))$, then $F$ is not ACL in $\mathbb{C}$. See Figure 2.1.

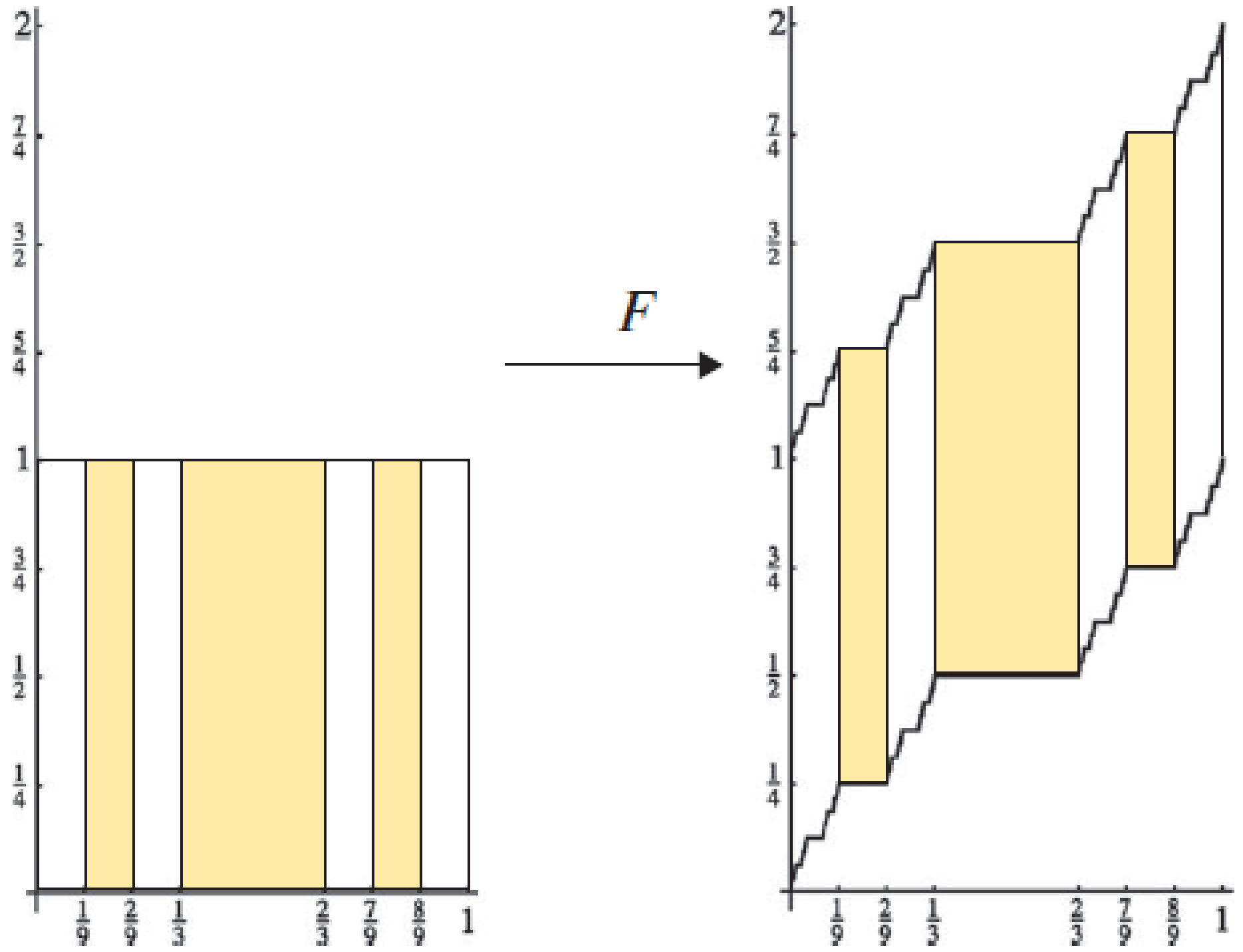

**Figure 2.1:** Illustrating the function $F$. Taken from [*Branner & Fagella 2014*].

A mapping $h : U \to V$, with $U, V \subseteq \mathbb{C}$, is $K-quasiconformal$ if and only if $h$ is a homeomorphism, $h$ is $ACL$ in $U$, and $|\partial_{\bar{z}} h| \leq k\,|\partial_z h|$ almost every where, with $\partial_z h := \frac{1}{2}\left(\frac{\partial h}{\partial x} - i\frac{\partial h}{\partial y}\right)$, $\partial_{\bar{z}} h := \frac{1}{2}\left(\frac{\partial h}{\partial x} + i\frac{\partial h}{\partial y}\right)$, $K := \frac{1+|\mu|}{1-|\mu|}$ defined as the *quasiconformal dilatation*, and the *complex dilatation* or *Beltrami coefficient of $h$* is defined as $\mu = \mu_h(z) := \frac{\partial_{\bar{z}} h(z)}{\partial_z h(z)}$, and $k := \frac{K-1}{K+1} < 1$.

Conversely, let $\mu(z)$ be a measurable complex-valued function defined on $U$ for which $\|\mu\|_\infty = k < 1$ almost everywhere, then we say that $\mu$ is a *k-Beltrami coefficient of* $U$. As a consequence, one may ask if there is a quasiconformal map $h$ satisfying the *Beltrami equation* $\partial_{\bar{z}} h(z) = \mu(z)\partial_z h(z)$. The answer is the next theorem:

THEOREM 2.1 *Measurable Riemann Mapping Theorem,* [*Ahlfors, 1966*]. *The Beltrami equation gives a one-to-one correspondence between the set of quasiconformal homeomorphisms of* $\overline{\mathbb{C}}$ *that fix the points* $0, 1, \infty$*, and the set of measurable complex-valued functions* $\mu$ *on* $\overline{\mathbb{C}}$ *for which* $\|\mu\|_\infty < 1$.

A geometric property of a $K - quasiconformal$ mapping is the following one. A domain whose boundary is a Jordan curve is called a *Jordan domain.* We define a *quadrilateral* $Q\,(z_1, z_2, z_3, z_4)$ as a Jordan domain and a sequence of four points $z_1, z_2, z_3, z_4 \in \partial Q$ following each other in positive orientation respect to $Q$.

Let $R = \{x + iy | 0 < x < a, 0 < y < b\}$, and let $g : Q \to R$ be an onto conformal mapping in such a way that $z_1, z_2, z_3, z_4$ are mapped to the vertices of $R$ under $g$, where the *arc* $[z_1, z_2] \subseteq \partial Q$ is mapped to the side $[0, a] \subseteq \partial R$. The number $a/b$ is called the *module* of the quadrilateral $Q$ and is denoted as

$$M(Q) = M\left(Q\left(z_1, z_2, z_3, z_4\right)\right) := \frac{a}{b}.$$

If $h : U \to V$ is a $K - quasiconformal$ mapping with $U, V \subseteq \mathbb{C}$, then

$$K(h) := \sup_{Q \subseteq U} \frac{M(h(Q))}{M(Q)} < \infty.$$

The value $K(h)$ coincides with the quasiconformal dilatation of $h$ defined above. An orientation-preserving homeomorphism $h : U \to V$ satisfying this inequality is a $K - quasiconformal$ mapping. See [*Lehto 1987*] and [*Lehto & Virtanen 1973*]. The equivalency is formally proven in [*Bers 1962*].

Now, suppose that $U \subseteq \mathbb{C}$, we define a *measurable field of ellipses* $\mathcal{E} \subseteq TU$ as a set of origin centered ellipses $\{E_u \subseteq T_u U \mid u \in U \text{ and } E_u \text{ is defined up to real scaling}\}$, where $TU$ is the tangent bundle of $U$ and $T_u U$ is its element for $u \in U$, see Figure 2.2. The map $u \mapsto \mu(u)$ from $U$ to $\mathbb{D}$ is Lebesgue measurable, with $\mu(u) := \dfrac{M - m}{M + m} e^{i2\theta}$, $M :=$ major axis of $E_u$, $m :=$ minor axis of $E_u$, $\theta$ is the argument of the direction of the minor axis in $[0, \pi)$. We denote $\mu(u)$ as the *Beltrami coefficient of* $E_u$, see Figure 2.3.

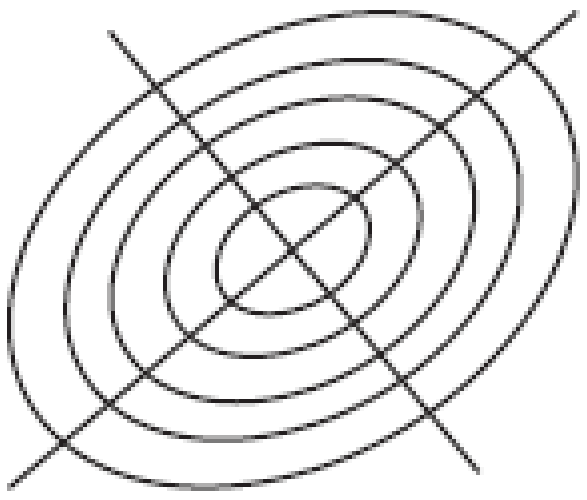

**Figure 2.2:** An element $E_u$ of $\mathcal{E}$. Taken from [*Branner & Fagella 2014*].

In this way, for every element of $E_u$ we can build a linear isomorphism $L_u(z) = az + b\bar{z} = a(z + \mu(u)\bar{z})$, with $\mu(u) = b/a = \dfrac{M-m}{M+m}e^{i2\theta}$, such that $L_u(E_u)$ is a circle. To the usual vector space $\mathbb{C}(\mathbb{C}, +, \cdot)$ with the standard complex sum and scalar multiplication, we call it the *linear standard conformal structure* $\sigma_0$. Additionally, $L_u$ induces a *linear conformal structure* $\sigma(u)$, which is a new vector space $\mathbb{C}(\mathbb{C}, +, \star)$ with new scalar multiplication $\star$, and $\mu(u)$ is again, defined as the *Beltrami coefficient of* $\sigma(u)$. And $\sigma(u)$ is related to $L_u$ and $E_u$ in the following sense.

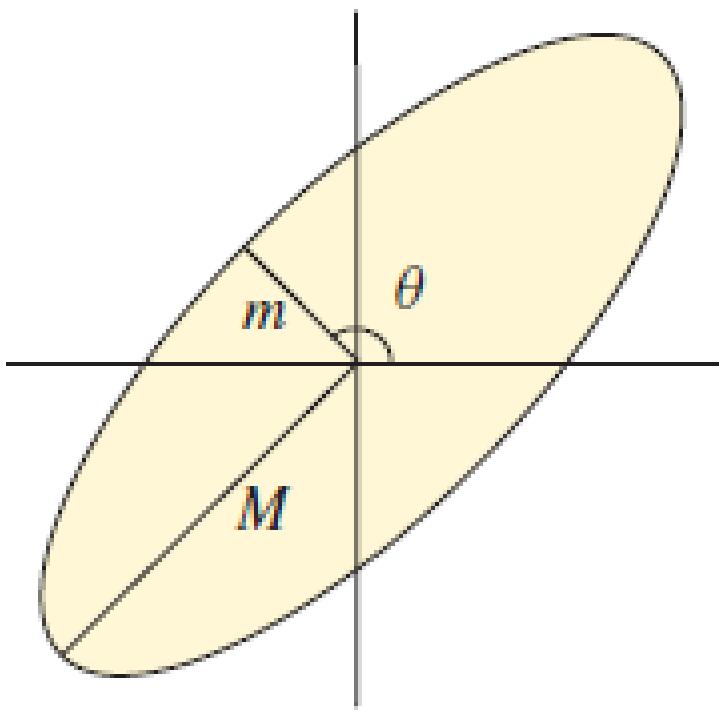

**Figure 2.3:** A representative ellipse of $E_u$ in $\mathcal{E}$. Taken from [*Branner & Fagella 2014*].

The standard complex sum, the real scalar multiplication, and the usual distribution remain. We define the complex scalar multiplication $c \star z := Re(c)z + Im(c)(i \star z)$, for $c \in \mathbb{C}$. Hence the problem is reduced to find an $\mathbb{R}$-linear map $J$ that multiplies by $i$ under $\star$, i.e., $J(z) = i \star z$. Thus $J(J(z)) = -z$ if $i \star i = -1$, and the linear conformal structure $\sigma(u)$ induced by $L_u$ is induced by $J = L_u^{-1} \circ (z \mapsto iz) \circ L_u$, making the tangent space $T_uU$ a $\mathbb{C}\,(\mathbb{C}, +, \star)$ linear vector space. See Figure 2.4.

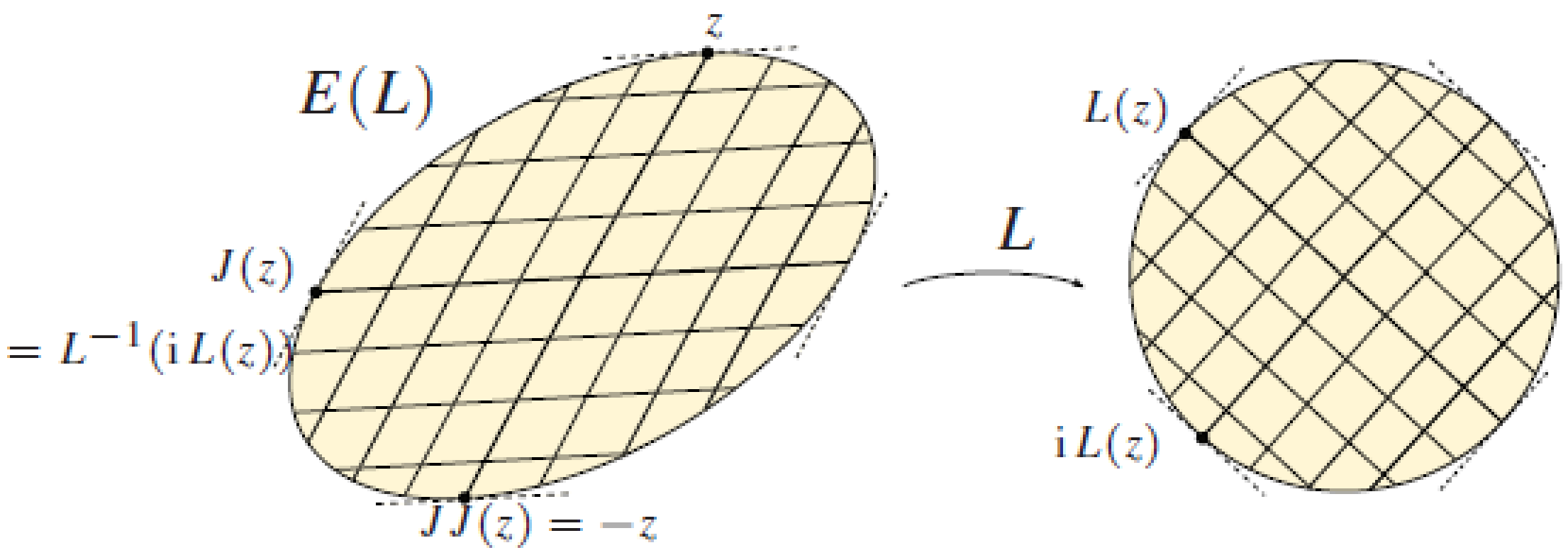

**Figure 2.4:** The map $L$ inducing a Linear Conformal Structure. Taken from [*Branner & Fagella 2014*].

Now, we define an *almost complex structure* $\sigma$ as $\bigcup_{u\in U} \sigma(u)$, and the *dilatation of* $\sigma$ *as* $K(\sigma) := \operatorname{ess\,sup}_{u\in U} K(u)$, where $K(u) := \dfrac{1+|\mu(u)|}{1-|\mu(u)|}$ is the *dilatation* of $E_u$. Note that the notions of a measurable field of ellipses $\mathcal{E}$ and of an almost complex structure $\sigma$, are equivalent, therefore the corresponding Beltrami coefficients and dilatations are equivalent as well.

In this context, let $\mathcal{D}^+(U,V)$ be the class of continuous orientation-preserving functions $h$ from $U$ to $V$, open subsets of $\mathbb{C}$, which are $\mathbb{R}$-differentiable almost everywhere with a non-singular differential $D_u h : T_u U \to T_{f(u)} V$ almost everywhere, and $D_u h = \partial_z h(u) dz + \partial_{\bar{z}} h(u) d\bar{z}$ depending measurably on $u$. With the previous discussion, $D_u h$ defines ellipses $E_u \subseteq T_u U$ via the inverse image of origin centered circles under $D_u h$, with Beltrami coefficient $\mu_h(z) := \dfrac{\partial_{\bar{z}} h(z)}{\partial_z h(z)}$ or, equivalently, a linear conformal structure induced under $D_u h$ on $T_u U$. See Figure 2.5.

In addition, we do this in every $u \in U$ where $h$ is differentiable. Hence, we obtain a measurable field of ellipses $\mathcal{E}_h$ and an almost complex structure $\sigma_h$ on $U$ with Beltrami coefficient $\mu_h$. We define $\mathcal{E}_h$ as the *pullback of* $\mathcal{E}_0$, the measurable field of circles, $\sigma_h$ as the *pullback of* $\sigma_0$, and $\mu_h$ as the *pullback of* $\mu_0 \equiv 0$. We denote them as

$$\mathcal{E}_h = h^*\mathcal{E}_0 \qquad \sigma_h(u) = h^*\sigma_0(u) \qquad \mu_h(u) = h^*\mu_0(u).$$

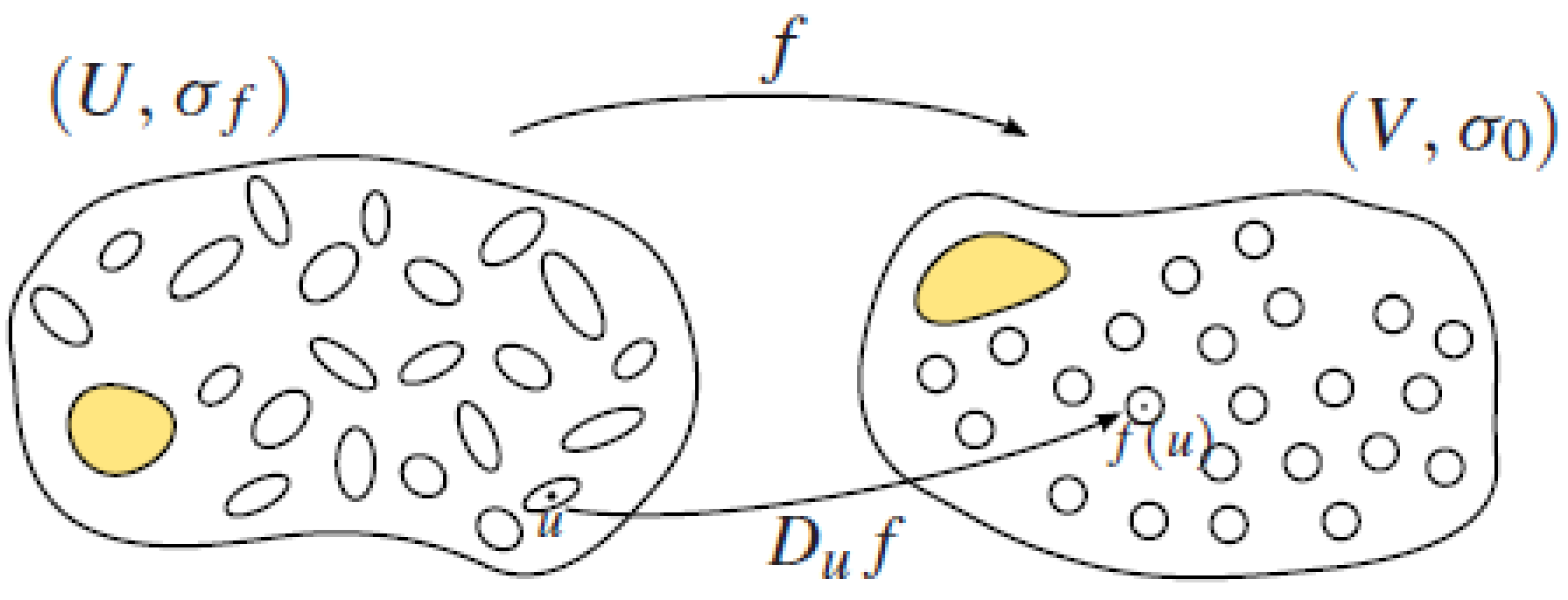

**Figure 2.5:** The pullback $\mathcal{E}_0$ under $f$. Taken from [*Branner & Fagella 2014*].

Now, we can pullback any almost complex structure $\sigma$ under a map $h$ if the inverse image of any set of measure zero is of measure zero, i.e., is *absolutely continuous with respect to the Lebesgue measure.* We write $\mathcal{D}_0^+(U,V)$ as the subclass of $\mathcal{D}^+(U,V)$ with this property. If $h \in \mathcal{D}_0^+(U,V)$ and $\mathcal{E} \subseteq TV$ is a measurable field of ellipses with Beltrami coefficient $\mu(v) = \mu$, *the pullback of* $\mathcal{E}$ *under* $h$ is the measurable field of ellipses $\mathcal{E}'$ with elements $E_u = (D_u h)^{-1}\left(E_{h(u)}\right)$ for almost every $u \in U$. See Figure 2.6.

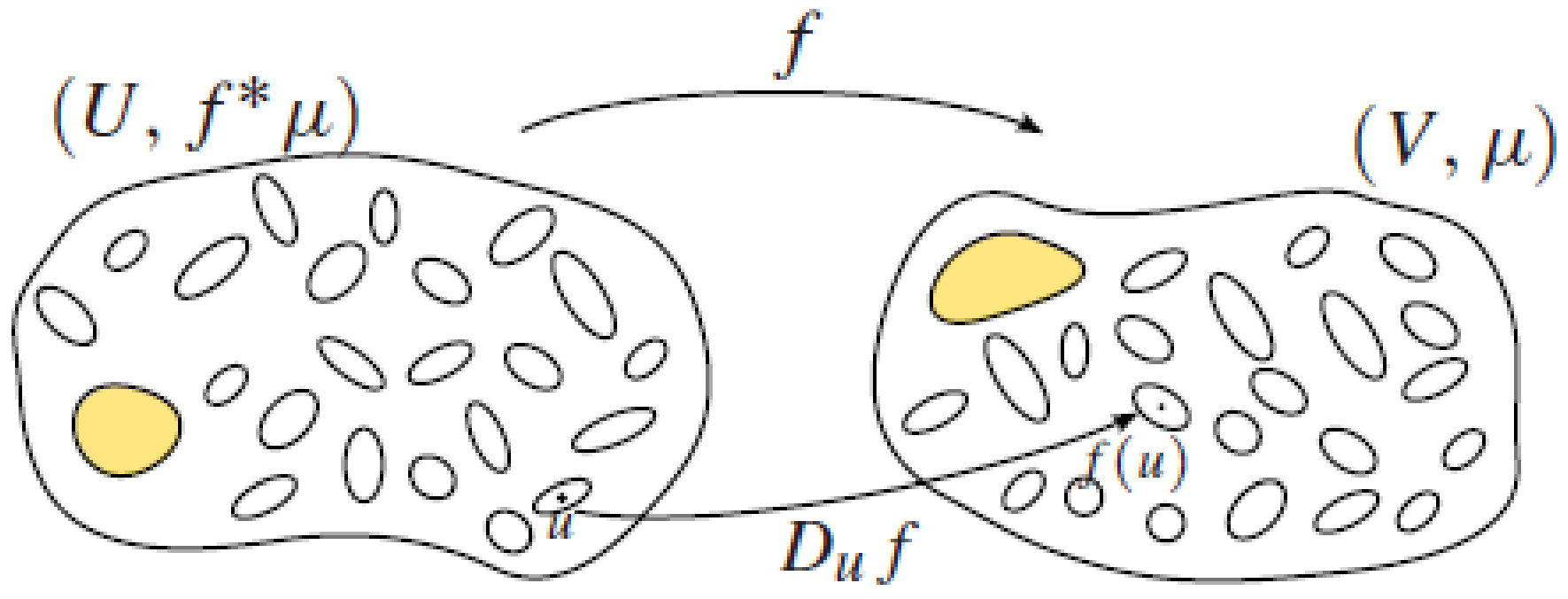

**Figure 2.6:** The pullback of $\mathcal{E}$ under $f$. Taken from [*Branner & Fagella 2014*].

In the case where $\mu$ is the Beltrami coefficient of a map $g : V \to W$ belonging to $\mathcal{D}^+(U,V)$, we can write $\mu = u_g$. Hence

$$h^*\mu_g = h^*\left(g^*\mu_0\right) = (g \circ h)^*\,\mu_0.$$

Observe that composing two linear isomorphisms $L_i(z) = a_i z + b_i \bar{z}, i = 1,2$ we obtain

$$L_1 \circ L_2(z) = \left(a_1 a_2 + b_1 \overline{b_2}\right) z + \left(a_1 b_2 + b_1 \overline{a_2}\right) \bar{z},$$

with Beltrami coefficient

$$\mu_{L_1 \circ L_2} = \frac{b_2 + \mu_{L_1} \overline{a_2}}{a_2 + \mu_{L_1} \overline{b_2}}.$$

Which implies

$$h^* \mu(u) = \frac{\partial_{\bar{z}} h(u) + \mu\left(h(u)\right) \overline{\partial_z h(u)}}{\partial_z h(u) + \mu\left(h(u)\right) \overline{\partial_{\bar{z}} h(u)}}.$$

In the case that $h$ is holomorphic, $\partial_{\bar{z}} h \equiv 0$ and $h^* \mu(u) = \mu\left(h(u)\right) \dfrac{\overline{\partial_z h(u)}}{\partial_z h(u)}$.

If $f \in \mathcal{D}^+(U,U)$ and $\sigma$ is an almost complex structure with Beltrami coefficient $\mu$ such that $f^* \mu(u) = \mu(u)$ for almost every $u \in U$, then we say that $\mu$ is $f-$*invariant*, and so is $\sigma$. Equivalently, we say that $f$ is *holomorphic* (in fact *conformal*) with respect to $\mu$ or $\sigma$.

Let $F \in \mathcal{D}_0^+(V,V), f \in \mathcal{D}_0^+(U,U), h \in \mathcal{D}_0^+(U,V)$ and let $\sigma$ be an almost complex structure $F$-invariant with Beltrami coefficient $\mu$. If

$$\begin{array}{ccccc} & & f & & \\ & U & \longrightarrow & U & \\ h & \downarrow & & \downarrow & h \\ & V & \longrightarrow & V & \\ & & F & & \end{array}$$

commutes then

$$f^* \left(h^* \mu\right) = \left(h \circ f\right)^* \mu = \left(F \circ h\right)^* \mu = h^* F^* \mu = h^* \mu.$$

In other words, if $F$ is $\sigma$-invariant, then $f$ is $(h^* \sigma)$-invariant.

Furthermore, Weyl´s Lemma states that if $h$ is quasiconformal in $U$ and $\partial_{\bar{z}} h \equiv 0$ almost everywhere, then there exists a holomorphic function $\tilde{h}$ in $U$ such that $\tilde{h}(z) = h(z)$ almost everywhere.

We define the *pushforward* of $h \in \mathcal{D}_0^+(U,V)$ as $h_* \mu := \left(h^{-1}\right)^* \mu$.

LEMMA 2.1 *If $F = h \circ f \circ h^{-1}$, such that $f$ is holomorphic and $h^*\mu_0 = \mu$, thus $F$ is holomorphic.*

Proof.

Note that:

$$F^*\mu_0 = \left(h \circ f \circ h^{-1}\right)^* \mu_0 = h_* f^* h^* \mu_0 = \left(h^{-1}\right)^* f^* \mu = \left(h^{-1}\right)^* \mu = \mu_0.$$

By Weyl's Lemma, $F$ is holomorphic.

■

Now, we can extend this theory to Riemann surfaces $S, S'$ as follows. Let $\phi : S \to S'$ and homeomorphism, if there exists a $K \geq 1$ so that $\phi$ is locally $K$-quasiconformal when is expressed in all the charts, then $\phi$ is *quasiconformal*. In the Figure 2.7, the composition $\varphi' \circ \phi \circ \varphi^{-1} : U \subseteq \mathbb{C} \to U' \subseteq \mathbb{C}$ is $K$-quasiconformal, while the transition mappings $\varphi \circ \tilde{\varphi}^{-1}$ and $\varphi' \circ \tilde{\varphi'}^{-1}$ are conformal.

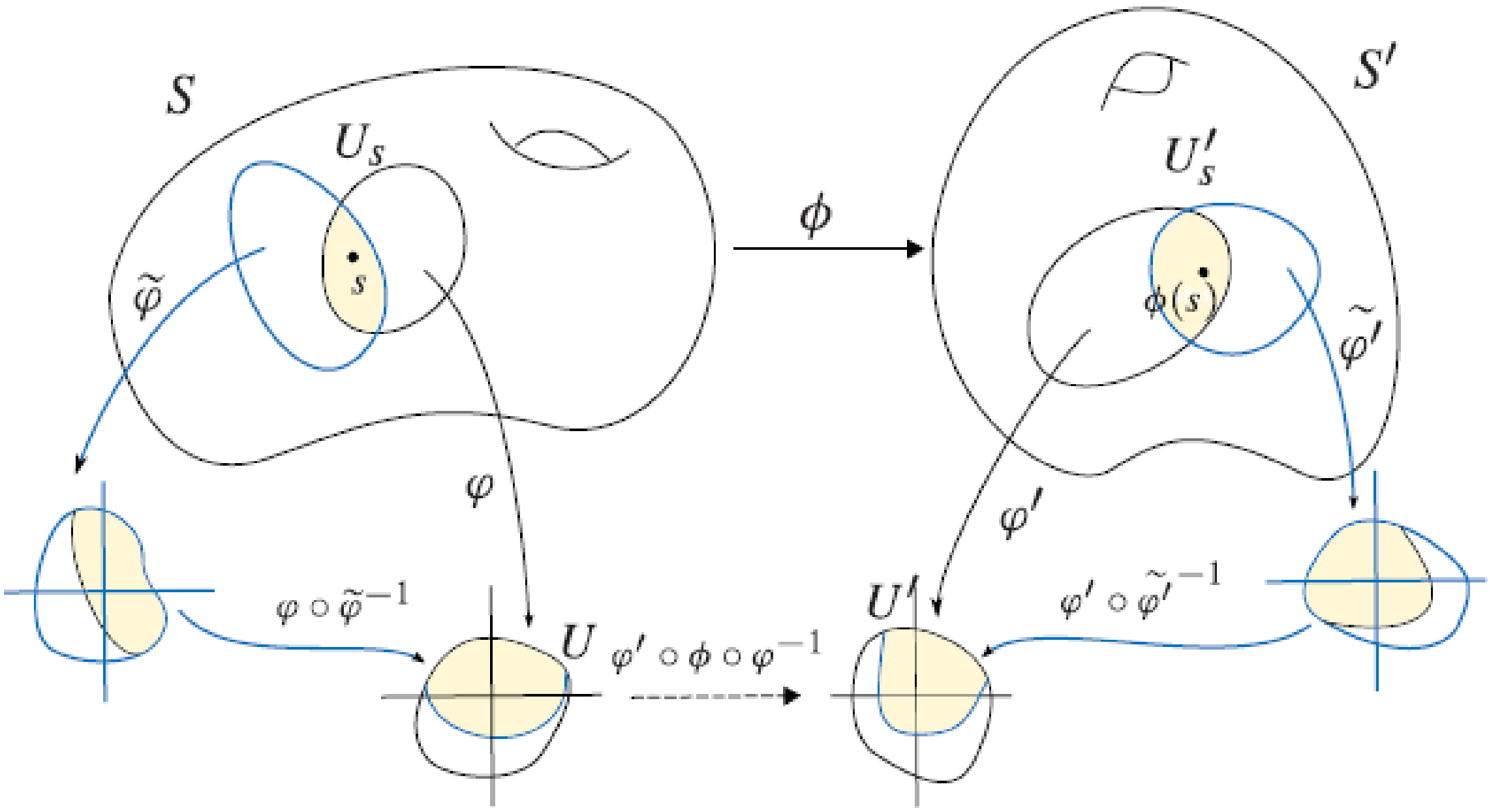

**Figure 2.7:** A $K$-quasiconformal $\phi$ between Riemann surfaces. Taken from [*Branner & Fagella 2014*].

A *Beltrami form* or a *Beltrami differential* $\mu$ on a Riemann surface $S$ is a $(-1, 1)$ differential on $S$, which is expressed as $\mu(z)d\bar{z}/dz$. This implies that if $\varphi : U_s \to U$ and $\tilde{\varphi} : \tilde{U}_s \to \tilde{U}$ are two overlapping charts on $S$ with holomorphic transition map $h = \tilde{\varphi} \circ \varphi^{-1}$, where $z = \varphi(s)$ and $\tilde{z} = \tilde{\varphi}(s)$, then the Beltrami form satisfies

$$\mu_{\varphi}(z) = \mu_{\tilde{\varphi}}(\tilde{z}) \frac{\overline{h'(z)}}{h'(z)}. \qquad (2.1)$$

Compare this with the above formula for a pullback of a Beltrami coefficient for a holomorphic function.

Now, we will study how the Beltrami form defines a measurable field of ellipses on the tangent bundle $TS$, and by equivalence, an almost complex structure. The Beltrami coefficient $\mu_{\varphi}(z)$ defines an ellipse up to real scaling at $T_zU$ with dilatation $K = \dfrac{1 + |\mu_{\varphi}(z)|}{1 - |\mu_{\varphi}(z)|}$.

An ellipse $E$ is mapped by the transition map $D_zh : T_zU \to T_{\tilde{z}}\tilde{U}'$ to an ellipse $\tilde{E}$, where $\tilde{E}$ is $E$ scaled by $|h'(z)|$ and rotated by the argument of $h'(z)$. See Figure 2.8.

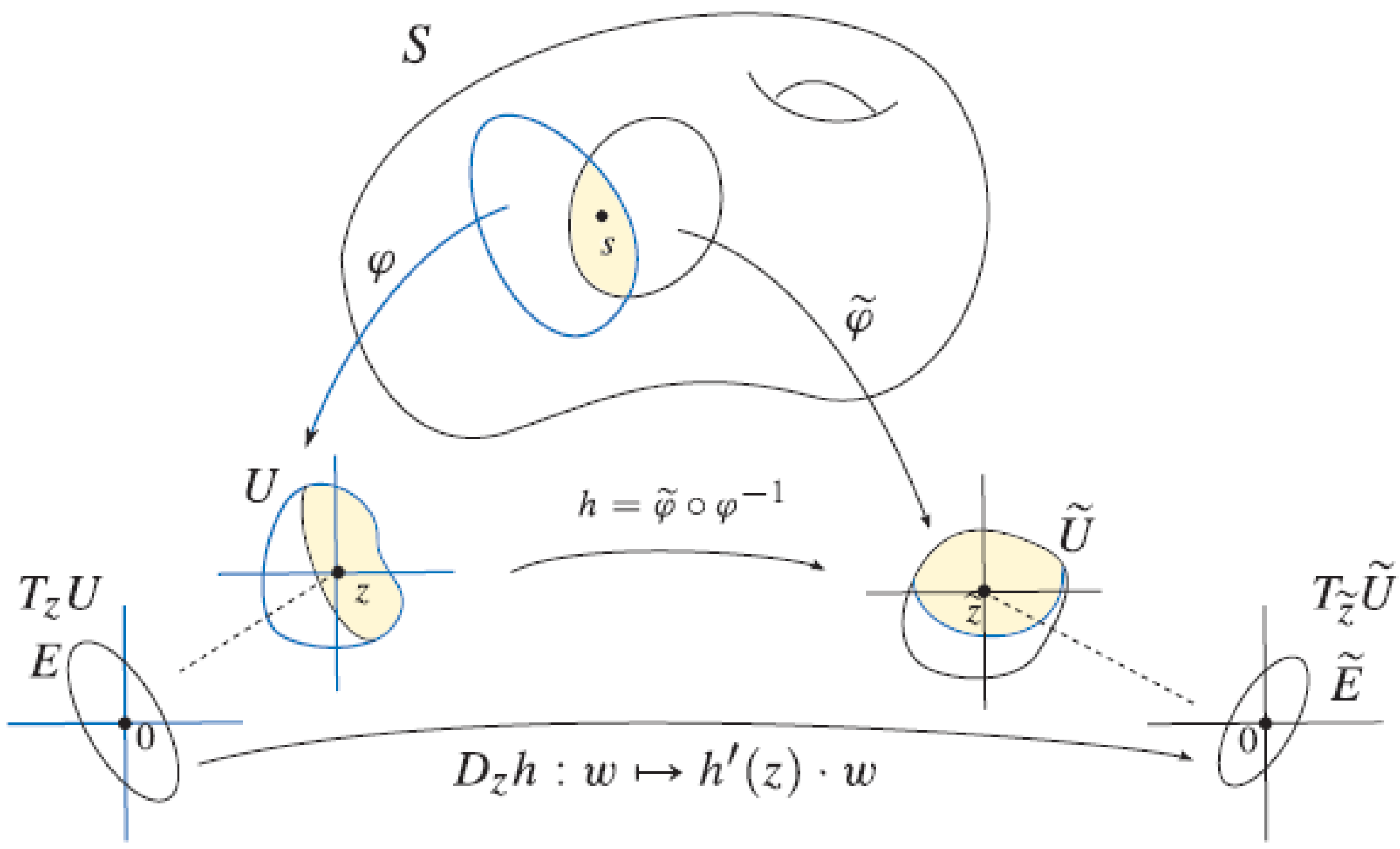

**Figure 2.8:** A field of ellipses on a Riemann surface $S$. Taken from [*Branner & Fagella 2014*].

Note that the equation (2.1) implies that $|\mu_{\varphi}(z)| = |\mu_{\tilde{\varphi}}(\tilde{z})|$. Therefore, the two ellipses have the same dilatation:

$$K(E) = K\left(\tilde{E}\right),$$

and again equation (2.1) implies that $arg\left(\mu_{\tilde{\varphi}}\left(\tilde{z}\right)\right) = arg\left(\mu_{\varphi}(z)\right) + 2arg\left(h'(z)\right)$. It follows that $\tilde{E}$ represents the ellipse up real scaling defined by the Beltrami coefficient $\mu_{\tilde{\varphi}}\left(\tilde{z}\right)$. Thus, the Beltrami form defines a measurable field of ellipses well defined in $TS$ and an almost complex structure on $S$.

Let $h : S \to S'$ be a quasiconformal mapping, with arbitrary charts $\varphi : U_S \to U$ and $\varphi' : U_{S'} \to U'$ on points $s \in S$ and $h(s) \in S'$ where $z = \varphi(s)$ and $z' = \varphi'\left(h(s)\right)$ respectively. If $\mu'$ is a Beltrami form on $S'$ then the *pullback* $h^*\mu'$ is defined as the Beltrami form on $S$ for which the chart $\varphi$ has the Beltrami coefficient

$$\left(h^*\mu'\right)_{\varphi}(z) = \frac{\partial_{\bar{z}} f(z) + \mu'_{\varphi'}\left(f(z)\right)\overline{\partial_z f(z)}}{\partial_z f(z) + \mu'_{\varphi'}\left(f(z)\right)\overline{\partial_{\bar{z}} f(z)}},$$

where $\mu'_{\varphi'}\left(z'\right)$ is the Beltrami coefficient of $\mu'$ in the chart $\varphi'$ and $z' = f(z) = \varphi' \circ h \circ \varphi^{-1}(z)$.

Now, let $f$ be a holomorphic function $f$. We want to deform $f$ via quasiconformal conjugation in such a way the deformations will be holomorphic. We can do this via a quasiconformal map $h$ such that $h^*\mu_0 = \mu$, and the Ahlfors-Bers-Morrey theorem guarantees its existence, see [*Branner & Fagella 2014*] and [*Zakeri & Zeinalian 1996*]:

THEOREM 2.2 (Ahlfors-Bers-Morrey) *Let $S$ be a simply connected Riemann surface isomorphic to $\mathbb{C}, \overline{\mathbb{C}}$ or $\mathbb{D}$, and let $\sigma$ be an almost complex structure on $S$ with measurable Beltrami form $\mu$. If the dilatation $K(\sigma)$ of $\sigma$ is uniformly bounded, then $\mu$ is integrable, i.e., there exists a quasiconformal homeomorphism $h$ from $S$ to its uniformization $\mathbb{C}, \overline{\mathbb{C}}$ or $\mathbb{D}$ which satisfies*

$$h^*\mu_0 = \mu.$$

*And $h$ is unique up to composition with automorphisms of the uniformization $\mathbb{C}, \overline{\mathbb{C}}$ or $\mathbb{D}$. Furthermore, in the case of $\overline{\mathbb{C}}$, $h$ fixes at least three points.*

## 2.2 Baker Laminations

In this thesis, the pinching deformations will be done on a Baker lamination, a set of geodesics in a Baker domain $U$ of an entire transcendental function $f$. Since the deformations will be done with quasiconformal theory, we work on a band whose central line is mapped to a geodesic in $U$.

Let $\alpha := \mathbb{R} \times \left\{\frac{\pi}{2}\right\}$ and let $B_\delta := \mathbb{R} \times \left(\frac{\pi}{2} - \delta, \frac{\pi}{2} + \delta\right) \subset \mathbb{C}$ with $0 < \delta < \frac{\pi}{2}$. Applying the exponential function, we define $\beta := exp(\alpha) = \{it \mid t \in \mathbb{R}^+\}$ and $V'_\delta(\beta) := exp\,(B_\delta) \subset \mathbb{H}$, where $V'_\delta(\beta)$ is called a *good neighborhood of thickness* $\delta$ of the complete geodesic $\beta$. If $\gamma$ is any other complete geodesic in $\mathbb{H}$, there is a unique oriented isometry $M \in PSL(2, \mathbb{R})$ of $\mathbb{H}$ such that $M(\beta) = \gamma$ and we say that $V'_\delta(\gamma) := M\,(V'_\delta(\beta))$ is a *good neighborhood (of thickness* $\delta$*) for* $\gamma$. See Figure 2.9.

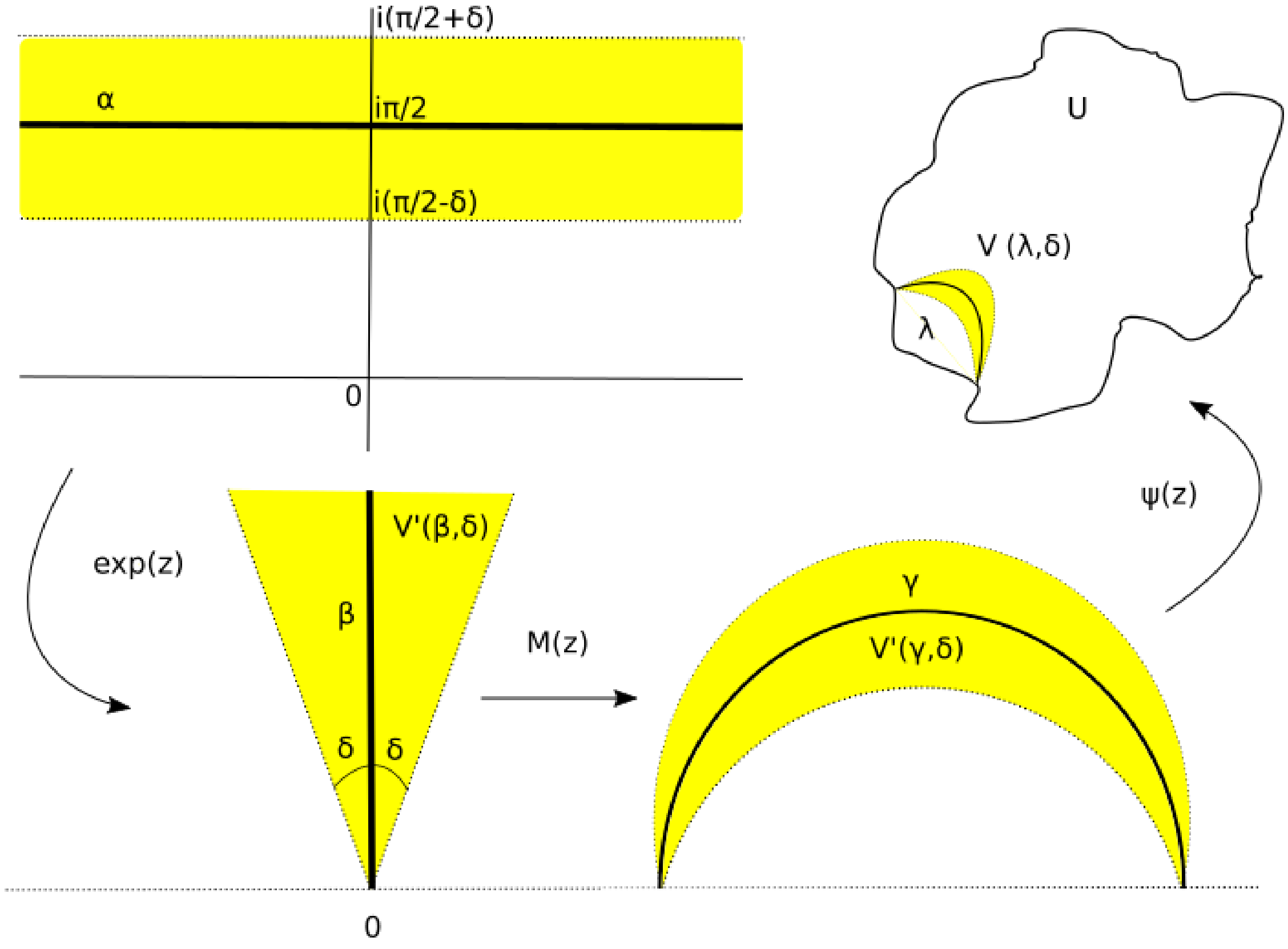

**Figure 2.9:** A good neighborhood $V_\delta(\lambda) = V(\lambda, \delta)$.

Now, suppose that $\lambda \in U$ is a complete geodesic, so $\psi^{-1}(\lambda) = \gamma$ is a geodesic in $\mathbb{H}$. Furthermore, assign a good neighborhood $V'_\delta(\gamma)$ to $\gamma$, and define $V_\delta(\lambda) := \psi\,(V'_\delta(\gamma))$ *good neighborhood (of thickness* $\delta$*) for* $\lambda$, where $\psi : \mathbb{H} \to U$ is the Riemann mapping.

On the other hand, we need to guarantee some conditions on the inverse images of neighborhoods around the Julia set in order to shrink the inverse images of the lamination. For this, we need a version of Theorem II of [*Mañé, 1993*] for entire transcendental functions. [*Bergweiler & Morosawa, 2002*] proved this case with the condition that $f$ *is semihyperbolic at* $a \in J(f)$, see Definition 2.1 below:

DEFINITION 2.1 *An entire transcendental function $f$ is* ***semihyperbolic at*** $a \in J(f)$*, if there exist $r > 0$ and $N \in \mathbb{N}$ such that for all $n \in \mathbb{N}$ and for all components $V$ of $f^{-n}(D_r(a)) := \{z \in \mathbb{C} \mid f^n(z) \in D_r(a)\}$ the function $f^n : V \to D_r(a)$ is a proper function of degree at most $N$, where $D_r(a) := \{z \in \mathbb{C} : |z - a| < r\}$. If this is satisfied for every $a \in J(f)$, we say that $f$ is* ***semihyperbolic****.*

Notice that an entire function $f$ is not semihyperbolic at a parabolic point, at a recurrent critical point nor at an asymptotic value, see [*Bergweiler & Morosawa, 2002*]. This means that $f$ is not semihyperbolic in $Sing(f^{-1})$, $P(f)$ (see proof of Lemma 4.2 below) nor at $\infty$.

THEOREM 2.3 [*Bergweiler & Morosawa, 2002*] *Let $f$ be an entire transcendental function and suppose that $f$ is semihyperbolic at $a \in J(f)$. Then there exists $r > 0$ with the following property: for all $\varepsilon > 0$, there exists $M \in \mathbb{N}$ such that if $n \geq M$ and $V$ is a component of $f^{-n}(D_r(a))$, then $diam(V) < \varepsilon$.*

In this context, we introduce the *Baker Laminations*, which is a natural setting for our pinching deformation.

DEFINITION 2.2 [*Robles & Sienra 2022*] *Let $f : \mathbb{C} \to \mathbb{C}$ be an entire transcendental function with a Baker domain $U$ such that $f(U) = U$ and provided with the hyperbolic metric. Let $\Lambda$ be a set of complete geodesics in $U$ with disjoint good neighborhoods. We say that $\Lambda$ is a* **Baker lamination of** *$U$, if the geodesics $\lambda \in \Lambda$, called* **leaves** *henceforth, satisfy:*

1. *The leaves of the lamination do not accumulate in $U$.*
2. *If $\lambda \in \Lambda$ then $f^n(\lambda) \in \Lambda$, with $n \in \mathbb{N}$. Also, $\lambda \subset U$ is in $\Lambda$, if $f^n(\lambda) \in \Lambda$, for some $n > 0$.*
3. *For any different leaves $\lambda, \lambda' \in \Lambda, \lambda \cap \lambda' = \emptyset$.*
4. *For any $\lambda \in \Lambda$, there exist different $\partial\lambda := \lim_{t \to \pm\infty} \lambda(t)$ and $\partial\lambda \subset \partial U \subset \overline{\mathbb{C}}$.*

The elements of the boundary $\partial\lambda$ are called *endpoints* and we denote $\overline{\lambda} := \lambda \cup \partial\lambda$. Set $\mathcal{L} := \bigcup_{k \in \mathbb{N}} f^{-k}(\Lambda)$, i.e., the grand orbit of $\Lambda$, and the connected $\mathcal{L} -$ *components* will be called the *grand leaves* as well. If $U$ is a univalent Baker domain, and $\Lambda$ is a Baker lamination in $U$, we say that $\Lambda$ is a *univalent Baker lamination.*

If $\bar{\Lambda} := \{\lambda \cup \partial\lambda \mid \lambda \in \Lambda\} \subseteq \bar{U}$, we define $\bar{\mathcal{L}} := \bigcup_{k\in\mathbb{N}} f^{-k}\left(\bar{\Lambda}\right)$. In the same way, let us denote $\mathcal{V} := \bigcup_{k\in\mathbb{N}} f^{-k}\left(V_\delta\left(\Lambda\right)\right)$, the grand orbit of disjoint good neighborhoods, and $\mathcal{U} := \bigcup_{k\in\mathbb{N}} f^{-k}(U)$, the grand orbit of the cycle of Baker domains.

*Remark:* It is important to note that for the convergence on Theorem C, we need necessarily that the components of $\mathcal{V}$ contain disjoint grand leaves of $\mathcal{L}$.

### 2.2.1 Constructing Baker Laminations

In this subsection, we will construct different examples of Baker laminations.

a) We say that $U$ is a *cycle of p-periodic components* $U_i$ if $U = \bigcup_{i=0}^{p-1} U_i$ with $U_i$ a p-periodic component, $f\left(U_i\right) \subseteq U_{i+1}$ if $i \neq p-1$, and $f\left(U_{p-1}\right) \subseteq U_0$. If $\Lambda_0$ is a Baker lamination in a p-periodic Baker domain $U_0$ under $f^p$, and $\partial U_i \cap \partial U_j = \emptyset$, we can induce a Baker lamination in all the components $U_i$ and in the cycle of Baker domains $U$ under iteration. Thus, $\Lambda_k = f^k\left(\Lambda_0\right)$ is a Baker lamination in $U_k$ under $f^p$, and $\Lambda = \bigcup_{i=0}^{p-1} \Lambda_i$ is a Baker lamination in $U$ under $f$. Good neighborhoods for each leaf in every lamination are constructed by iteration of the good neighborhoods of the leaves of $\Lambda_0$.

b) Let $U_i$ be a univalent Baker domain of period $p$ under the function $f$, with two points $u, v \in \partial U_i \cap I_f$ where $I_f := \left\{z \in \mathbb{C} \mid \lim_{n\to\infty} f^n(z) = \infty\right\}$ and $\partial U_i \subseteq J(f) = \partial I_f$, see [Morosawa, Nishimura, Taniguchi & Ueda, 1998].

We say that an almost open set $R \subset U_i$ (i.e., $R$ is the disjunctive union of an open set with a set of Lebesgue measure zero) is a *fundamental domain of* $f$ in $U_i$, if any different points $u, v \in R$ are in different grand orbits under $f$, and if $w \in U_i$ there exists $z \in R$ and $n \in \mathbb{Z}$ such that $f^n(z) = w$.

If $u$ and $v$ belong to a fundamental domain $R$ of $f$ in $U_i$, we can construct a geodesic $\gamma_{u,v}$ with endpoints $u$ and $v$ non-exceptional points, such that $\gamma_{u,v} \cap f^p\left(\gamma_{u,v}\right) = \emptyset$, is disjoint to $Sing(f^{-1})$ and to $P(f)$, then the forward and backward orbits of $\gamma_{u,v}$ under $f^p$ make a Baker lamination $\Lambda$ in $U_i$. Note that since $R$ contains an open set and $f^n(R) \neq R$ for every $n \in \mathbb{Z}$, for every $\gamma \in \Lambda$ there is a unique $n_0 \in \mathbb{Z}$ such that $\gamma \in f^{n_0 p}(R)$, therefore the leaves of $\Lambda$ do not accumulate in $U$.

c) As defined above, we map a vertical geodesic as $\beta = \{it \mid t \in \mathbb{R}\} \subset \mathbb{H}$ under the Riemann mapping $\psi : \mathbb{H} \to U$ and suppose that we can extend $\psi(\infty) = \infty$. If $U$ is of period $p$, we define $\lambda_{\zeta,\infty} := \{\psi(\beta), f(\psi(\beta)), ..., f^{p-1}(\psi(\beta))\}$, where $\zeta = \psi(0)$. We are interested in the case when $\lambda_{\zeta,\infty}$ is a Baker lamination $\Lambda$. Henceforth we will consider $\lambda_{\zeta,\infty}$ satisfying Definition 2.2.

Let us interpret this geometrically. If $U$ is a Riemann surface, and $f$ is an endomorphism of $U$, a univalent Baker domain of hyperbolic type I, $f|_U$ is conjugated to $g : \mathbb{H} \to \mathbb{H}$ with $g(w) = aw$. $\mathbb{H}/g$ is an annulus $\mathcal{A}$ with *core geodesic* $\tilde{\beta} := \beta/g$, defined as the unique closed geodesic in the hyperbolic annulus. Since $g(1) = a$, thus *length* $\left(\tilde{\beta}\right) = log(a) - log(1) = 2\pi/mod(\mathcal{A})$. See [*Mc Mullen 1994*].

If we take a geodesic $\lambda$ in $U$, with $\infty$ as one of its endpoints, one may ask if $\lambda/f$ is the core geodesic or not. When $\lambda/f$ is the core geodesic, then we will denote it by $\lambda_{\zeta,\infty}$ and it will be the only invariant leaf containing $\infty$ in a univalent Baker lamination. If not, $\{f^n(\lambda)\}_{n\in\mathbb{Z}}$ is a set of infinite geodesics that share $\infty$ as an endpoint by uniformization. In other words, if $b \neq \zeta$, there would be preimages of the geodesic $\lambda_{b,\infty}$ accumulating on $\lambda_{\zeta,\infty}$ in contradiction with condition (1) of Definition 2.2, thus if $\lambda_{\zeta,\infty} \in \Lambda$ there are no other leaves with one endpoint at $\infty$. See Figure 2.10.

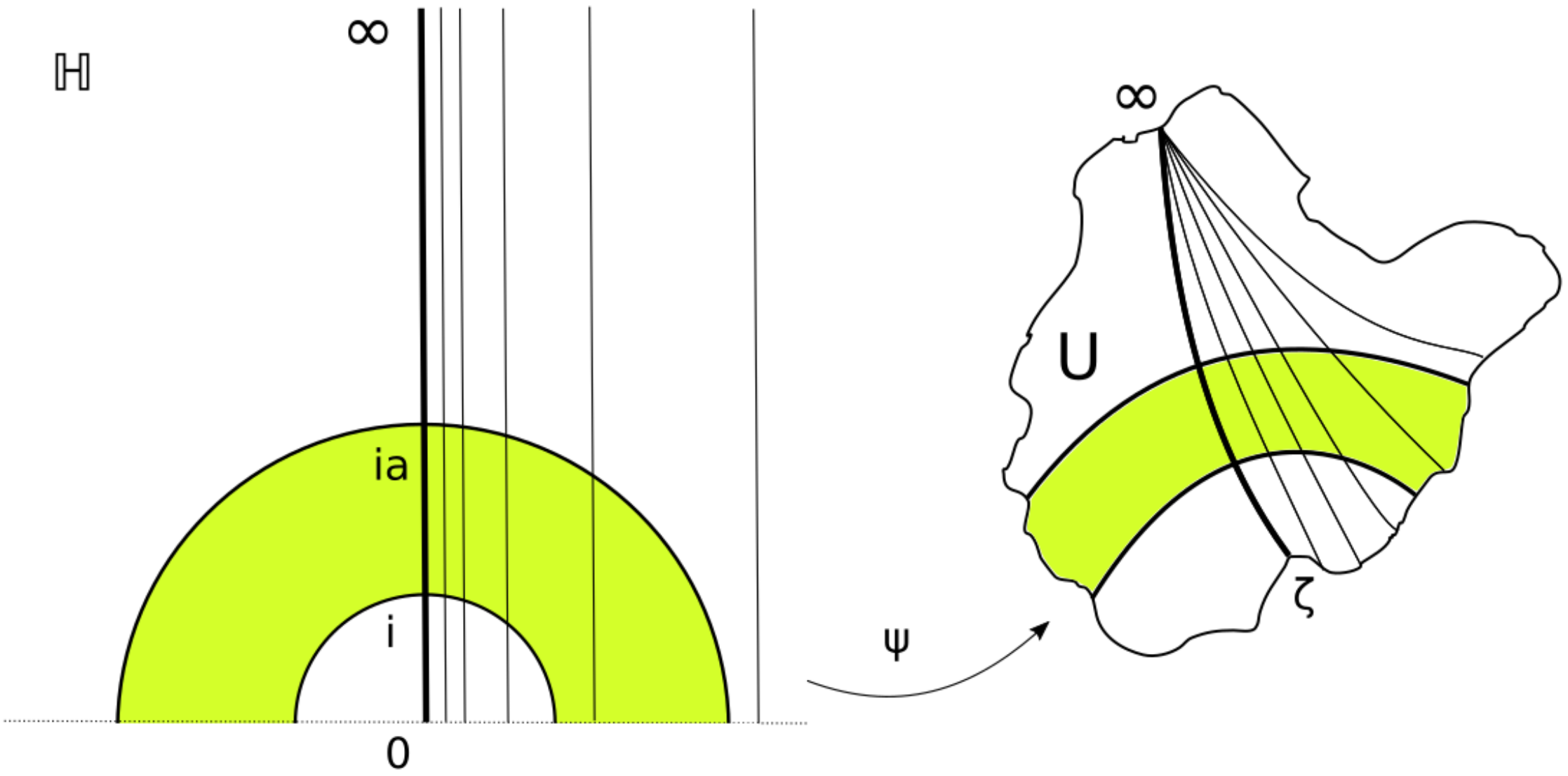

**Figure 2.10:** The leaf $\lambda_{a,\infty}$ in a univalent Baker domain $U$ of Hyp. T. I.

### 2.2.2 Closed curves in Baker Laminations

In this subsection, we discuss the non-existence of a closed curve built from the closure of leaves in a Baker lamination of an entire transcendental function where finite Baker domains intersect in their boundaries, and some of them are not completely invariant. The existence of such a curve would not allow the uniform convergence of the pinching deformation, because this closed curve would be pinched to a point under the pinching deformation and $\overline{\mathbb{C}}$ would become the wedge product of two spheres of dimension two at each point of the Julia set.

First, since the Baker domains of a transcendental entire function are simply connected, there exists no closed complete geodesic as a leaf for a Baker lamination.

Second, since the endpoints of leaves in a Baker lamination are in the boundary of a Baker domain, we want to know the circumstances where the components of the boundaries of Baker domains intersect each other (besides $\infty$). Because there could be two leaves of a Baker lamination conforming a cycle, as the example given by [*Sienra 2006*] where two or more univalent Baker domains intersect among them in their boundary in a non-singular fixed point.

PROPOSITION 2.1 *If $f$ is an entire transcendental function with $U$ and $V$ disjoint periodic univalent Baker domains of $f$ such that $f^n(U) = U$, $f^m(V) = V$, and $U$ is not completely invariant, where $\{p, q\} \subset \partial U \cap \partial V$, then $p = \infty$ and $f^{nm}(q) = q$.*

Proof.

Since $f^{nm}(U) = U$ and $f^{nm}(V) = V$, we define $g := f^{nm}$. Thus $g(U) = U$ and $g(V) = V$.

Let $\{p, q\} \subset (\partial U \cap \partial V) \setminus \{\infty\}$, let $\gamma_U$ be a simple curve in $U$ such that $\partial\gamma_U = \{p, q\}$ and let $\gamma_V$ be a simple curve in $V$ such that $\partial\gamma_V = \{p, q\}$. We define $\Gamma := \gamma_U \cup \gamma_V$, then $\overline{\Gamma}$ is a closed simple curve, boundary of the bounded disk $D_\Gamma$.

By the maximum principle,

$$max \left\{ |g^k(z)| : z \in \overline{D_\Gamma} \right\} \in \left| \partial \left( g^k \left( D_\Gamma \right) \right) \right|,$$

and since the Fatou set and the Julia set are completely invariant

$$\partial \left( g^k \left( D_\Gamma \right) \right) \subset \left( U \cup V \cup \left\{ g^k(p), g^k(q) | k \in \mathbb{N} \right\} \right),$$

then for $U_0 := g^{-1}(U) \setminus \overline{U}$ an open preimage of $U$ where $g^{-1}(U) \neq U$, since $U$ is not completely invariant,

$$\lim_{k\to\infty} g^k\left(D_\Gamma\right) \subset \overline{\mathbb{C}} \setminus U_0.$$

But this is a contradiction with the fact that for every open set $V$ containing a point of $J(g)$, $\overline{\mathbb{C}} \setminus E \subset \bigcup_{n=1}^{\infty} g^n(V)$, where $E$ is the set of all the exceptional points with $\#E \leq 2$ (see [*Morosawa et al. 1998*]), since $\#U_0 > 2$, $U_0$ is an open set. Therefore $p = \infty$.

Furthermore, $g(q) = q$. In the case that $g(q) \neq q$, if $q, g(q) \in \mathbb{C}$ the precedent argument implies that $\infty \in \{q, g(q)\}$. But $q$ can not be $\infty$ since it is an essential singularity of the entire function $g$, and if $g(q) = \infty$, $q$ would be a pole which is not possible since $g$ is an entire function.

■

As an immediate consequence, we have the next corollary:

COROLLARY 2.1 *If $f$ is an entire transcendental function with $U$ and $V$ disjoint periodic Baker domains of $f$ such that $f^n(U) = U$ and $f^m(V) = V$, $U$ is not completely invariant and $\Lambda$ a Baker lamination on $U$ which does not contain a leaf of type $\lambda_{a,\infty}$ where $f^{nm}(a) = a$. If $\lambda_1 \in \mathcal{L} \cap U$, $\lambda_2 \in \mathcal{L} \cap V$ and $\partial\lambda_i \cap \{a\} \neq \emptyset$, then $\overline{\lambda_1} \cap \overline{\lambda_2} = \emptyset$.*

Corollary 2.1 implies that the only way to make a closed curve with the closure of some leaves of a Baker lamination is with their endpoints as fixed points of the function in the boundary of some Baker domains.

PROPOSITION 2.2 *If $f$ is an entire transcendental function with $U_i$ disjoint periodic Baker domains of $f$ such that $f^{m_i}(U_i) = U_i$ and $\lambda_i$ are leaves in $\mathbb{C}$ on a grand orbit $\mathcal{L}$ of a Baker lamination $\Lambda$, with $\lambda_i \subset U_i$, $i \in A$, $A$ a finite set of indexes, such that $\overline{U_i} \cap \overline{U_k} \neq \emptyset$ for some $k \in A \setminus \{i\}$, and some $U_{i_j}$ are not completely invariant. Then $\bigcup_{i\in A} \overline{\lambda_i}$ is not a closed curve.*

Proof.

The proof is analogous to the one of Proposition 2.1.

Since $f^{\prod_{i\in A} m_i}(U_i) = U_i$, we define $g := f^{\prod_{i\in A} m_i}$. Thus $g(U_i) = U_i$ for every $i \in A$.

Let us suppose that $\bigcup_{i\in A} \overline{\lambda_i}$ is a closed curve denoted as $\Gamma$, hence is the boundary of the bounded disk $D_\Gamma$.

By the maximum principle,

$$max\left\{|g^n(z)| : z \in \overline{D_\Gamma}\right\} \in |\partial\left(g^n\left(D_\Gamma\right)\right)|,$$

and since the Fatou set and the Julia set are completely invariant

$$\partial\left(g^n\left(D_\Gamma\right)\right) \subset \left(\left(\bigcup_{i\in A} U_i\right) \cup \left(\bigcup_{i\in A} g^n\left(\partial\lambda_i\right)\right)\right),$$

then for $U_0 := g^{-1}(U_{i_j}) \setminus \overline{U_{i_j}}$ an open preimage of $U_{i_j}$ where $g^{-1}\left(U_{i_j}\right) \neq U_{i_j}$, since $U_{i_j}$ is not completely invariant,

$$\lim_{n\to\infty} \bigcup_{n\in\mathbb{N}} g^n\left(D_\Gamma\right) \subset \overline{\mathbb{C}} \setminus U_0.$$

But this is a contradiction with the fact that for every open set $V$ containing a point of $J(g)$, $\overline{\mathbb{C}} \setminus E \subset \bigcup_{n=1}^{\infty} g^n(V)$, where $E$ is the set of all the exceptional points with $\#E \leq 2$, (see [*Morosawa et al. 1998*]), since $\#U_0 > 2$, because $U_0$ is an open set.

■

Therefore, in entire transcendental functions with a finite number of Baker domains intersecting each other and with some Baker domains not completely invariant, we can not have laminations whose leaves conform a closed curve. This proof is still valid with an infinite number of Baker domains where there exists the least common multiple of their periods.

To my knowledge, there is no explicit example where the case with an infinite number of Baker domains intersecting each other and some Baker domains not completely invariant could occur. However, examples of entire transcendental functions with infinitely many $p$-cycles of Baker domains for each $p \in \mathbb{N}$ was proved by [*Rippon & Stallard 1999*]. For this reason, we introduce the following definition and we will add it as a hypothesis in Theorem C.

DEFINITION 2.3 *Let $f$ be an entire transcendental function with a Baker domain $U$, and a Baker lamination $\Lambda \subseteq U$. We say that $f$ satisfies the* ***property P*** *if $f$ does not have an infinite number of Baker domains intersecting each other with no least common multiple of their periods and the grand orbit $\mathcal{L}$ does not contain a closed curve.*

In this context, we remark on the following. Let $\mathcal{L}$ be a grand orbit for $f$, and let

$\gamma$ be a connected path in $\overline{\mathcal{L}}$ in such a way that between any two grand leaves of $\mathcal{L}$ there is another one (except perhaps at the extremes of $\gamma$). This path is made up of a numerable set of grand leaves in $\mathcal{L}$, and their endpoints in $J(f)$.

This set of endpoints is a Cantor set. In the process of the pinching deformation, the grand leaves will tend continuously to a point. At the end of the pinching deformation, the path $\gamma$ is shrunk into a path but not to a point. See the devil´s staircase as an example of a continuous function mapping a cantor set into an interval.

## 2.3 Pinching Deformations on Baker Laminations

It is common to study the theory of holomorphic dynamics introducing deformations of a function via conjugation classes, i.e., analyzing a certain space of functions. One of these tools is the pinching deformation introduced by [*Makienko 2000*] to prove that the component of J-stability is unbounded in $\mathbb{CP}^{2d+1}$ for rational functions with disconnected Julia sets as well as with connected Julia sets with some restrictions on accesses. [*Tan 2002*] generalized this concept and gave it a different approach.

In general terms, it is to take an invariant curve $\gamma$ with an attractor fixed point and a repulsor fixed point as boundaries of a function $f$ and deform $f$ via quasiconformal conjugations, to shrink $\gamma$ to a point. An intuitive explanation of this is to fuse the attractor point with the repulsor point creating a parabolic fixed point. This is exactly the behavior that we want to generalize to entire transcendental functions with Baker domains. We will follow [*Haissinsky & Tan 2004*] very closely to build the deformations in this particular case.

Let $f$ be an entire transcendental function with at least one cycle of periodic Baker domains $U = \{U_0, U_1, \ldots, U_{p-1}\}$ with a Baker lamination $\Lambda$ in $U$.

Take $L_b, L_y, L_r \in \mathbb{R}$ such that $0 < L_b < L_y < L_r$ and a function $\tau : [0, 1) \to [L_r, \infty)$ such that $\tau \in C^1[0, 1)$ is an increasing function. Using $\tau$ we build a closed set $M \subset \mathbb{R}^2$ bounded by

$$([0,1] \times \{L_b\}) \cup (\{0\} \times [L_b, L_r]) \cup (\{1\} \times [L_b, \infty)) \cup \{(t, \tau(t)) \mid t \in [0,1)\}.$$

Now choose $v_t(y)$ such that $v_t(y) = y$ for $y \in [L_b, L_y]$ and $(t, y) \mapsto (t, v_t(y))$ is a $C^1$-diffeomorphism from $[0,1] \times [L_b, L_r] \setminus \{(1, L_r)\}$ onto $M$. See Figure 2.11.

Here we introduce a technical assumption that will be used in the proof of Lemma 4.8:

*Technical Assumption.* For any $L' < L_r$, there is $t(L') \in (0,1)$ with $t(L') \to 1$ as $L' \to L_r$ such that for any $(s,y) \in (t(L'), 1] \times [L_b, L']$, we have $v_s(y) = v_{t(L')}(y)$. It means that the deformation of the yellow zone after the point $t(L')$ is constant, as can be seen in Figure 2.11 below from $L'$.

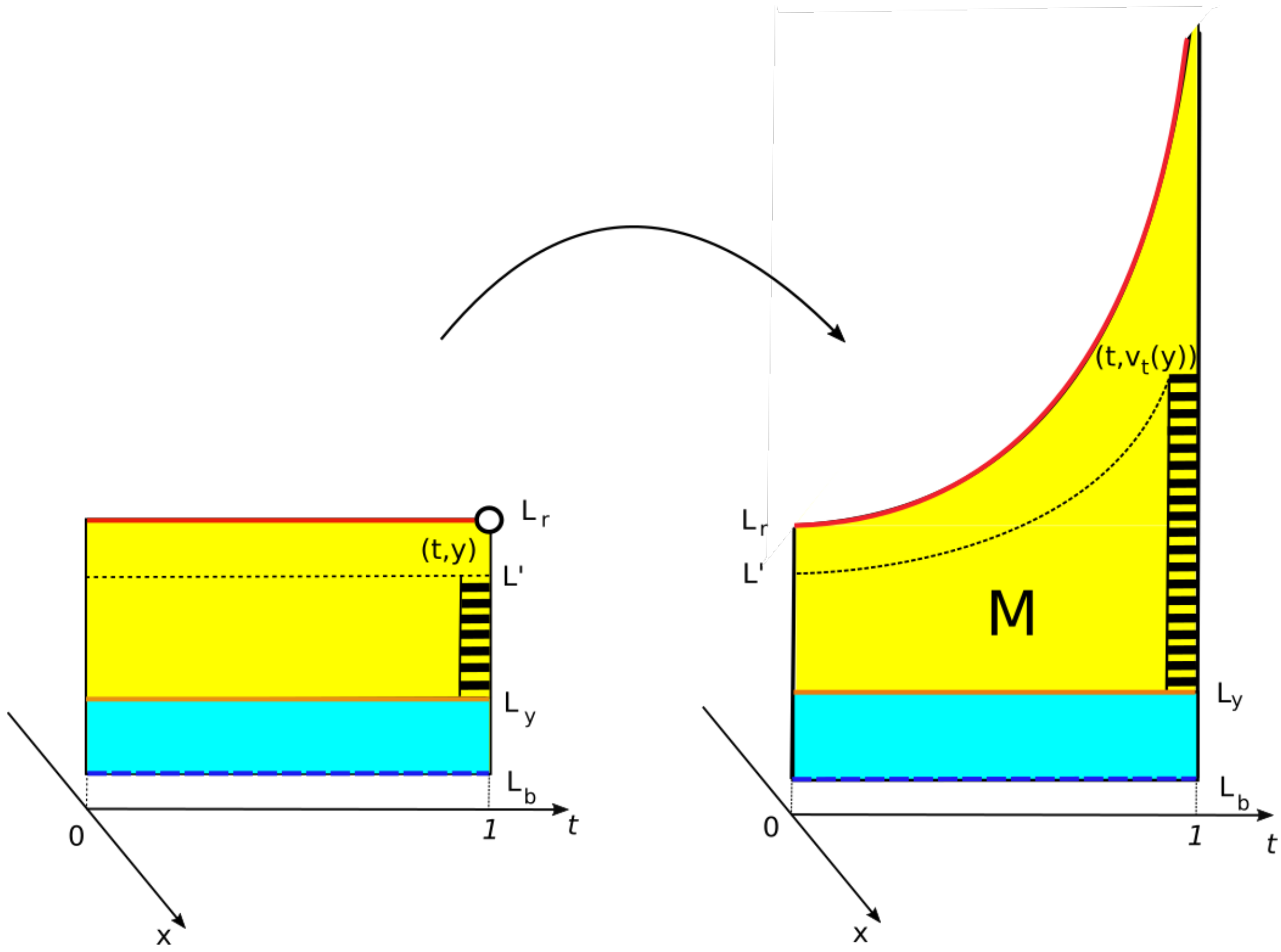

**Figure 2.11:** The diffeomorphism $(t, y) \to (t, v_t(y))$.

With $v_t$, we build a map $\widetilde{P}_t$ defined on the strip $\{x + iy \mid x \in \mathbb{R}, y \in [L_b, L_r]\}$ with $t \in [0,1]$ where

$$\widetilde{P}_t(x + iy) = x + iv_t(y)$$

and it has the following properties:

1. It commutes with any real translation.
2. It is the identity map on $\mathbb{R} \times [L_b, L_y]$.

3. The coefficient of its Beltrami form is

$$\frac{\partial_{\bar{z}}\widetilde{P}_t}{\partial_z\widetilde{P}_t}(x+iy)=\frac{1-\partial_y v_t(y)}{1+\partial_y v_t(y)},$$

which is continuous on $(t,x+iy)\in[0,1]\times\mathbb{R}\times[L_b,L_r]\setminus\{(1,x,L_r)\}$, its norm is locally uniformly bounded away from 1 and tends to 1 as $(t,x,y)\to(1,x,L_r)$, for every $t\in[0,1]$.

4. $\widetilde{P}_t(z)=-1/\widetilde{P}_t(-1/z)$ is continuous in $(t,z)\in[0,1]\times\mathbb{R}\times[L_b,L_r]\setminus\{(1,x,L_r)\}$. For $t<1$, $P_t$ is injective.

Now we have to connect these bands with those in the Baker lamination. Let $B_\delta^+:=\mathbb{R}\times\left[\frac{\pi}{2},\frac{\pi}{2}+\delta\right]$, i.e., the upper part of $B_\delta$, we define the holomorphic map $S_+:\mathbb{R}\times[L_b,L_r]\to B_\delta^+$ as

$$S_+(z)=\frac{\delta L_r L_b}{L_r-L_b}\left(-\frac{1}{z}\right)+i\left(\frac{\pi}{2}-\frac{\delta L_b}{L_r-L_b}\right).$$

Also, we define the map

$$\phi_+:=\Psi\circ M\circ exp\circ S_+:\mathbb{R}\times[L_b,L_r]\to V^+\subseteq U_i$$

where $V^+\subseteq V$, with $V$ a good neighborhood of $\lambda\in\Lambda$, and with a well defined inverse branch $\psi_+:V^+\to\mathbb{R}\times[L_b,L_r]$.

For $t\in[0,1)$, set $(\sigma_t')_+:=\left(\widetilde{P}_t\circ\psi_+\right)^*(\sigma_0)$ on $V^+$ to be the pullback of the standard almost complex structure on $\{x+iv_t(y)\mid t\in[0,1)\}$.

Analogously, we do the same proceeding for $V^-\subseteq V$. Let $B_\delta^-:=\mathbb{R}\times\left[\frac{\pi}{2}-\delta,\frac{\pi}{2}\right]$, we define the holomorphic map $S_-:\mathbb{R}\times[L_b,L_r]\to B_\delta^-$ as

$$S_-(z)=\frac{\delta}{L_r-L_b}(z-iL_r)+i\left(\frac{\pi}{2}\right).$$

We define the map

$$\phi_-:=\Psi\circ M\circ exp\circ S_-:\mathbb{R}\times[L_b,L_r]\to V^-\subseteq U_i,$$

where $V^-\subseteq V$,with inverse branch $\psi_-:V^-\to\mathbb{R}\times[L_b,L_r]$. For $t\in[0,1)$, let us set $(\sigma_t')_-:=\left(\widetilde{P}_t\circ\psi_-\right)^*(\sigma_0)$ on $V^-$ be the pullback of the standard almost complex structure on $\{x+iv_t(y)\mid t\in[0,1)\}$.

Then, we spread $(\sigma'_t)_+$ and $(\sigma'_t)_-$ to the grand orbit $\mathcal{V}$ by defining

$$\sigma_t := \bigcup_n \left( (f^n)^* \left( (\sigma'_t)_+ \right) \cup (f^n)^* \left( (\sigma'_t)_- \right) \right),$$

and we define $\sigma_t$ outside of $\mathcal{V}$ on the Riemann sphere by setting $\sigma_t := \sigma_0$.

Even though this "yellow zone" is where the deformation will happen, we will need another "blue" zone around the "yellow zone" that will be used in the proof of Theorem C. For example, in Lemma 4.3 we will use this blue zone to use quasiconformal theory in moduli of annuli.

DEFINITION 2.4 *We define the* ***yellow strip neighborhood*** *as*

$$Y(\lambda) := \phi_+ (\mathbb{R} \times [L_y, L_r]) \cup \phi_- (\mathbb{R} \times [L_y, L_r]).$$

*In the same way, we define the* ***blue strip neighborhood*** *as*

$$B'(\lambda) := \phi_+ (\mathbb{R} \times (L_b, L_r]) \cup \phi_- (\mathbb{R} \times (L_b, L_r]).$$

*Additionally, we define $Y^*(\lambda) := Y(\lambda) \cup \partial\lambda$.*

*We will join to $B'(\lambda)$ two suitable neighborhoods of radius lesser to the $r$ from semi-hyperbolicity (see Definition 2.1), $\Delta_\alpha$ and $\Delta_\beta$ of the two point set $\partial\lambda = \{\alpha, \beta\}$. We will do this in such a way that we do not intersect any other leaf $\lambda'$ of the lamination $\Lambda$ nor others neighborhoods $\Delta_{a'}$ or $\Delta_{\beta'}$ of $\partial\lambda'$, except $\infty$ possibly. Additionally, consider the radii of these discs to be between $r_\beta/2$ and $r'_\beta/6$ where $r'_\beta$ is the radius $r$ of Theorem 2.3 and $r_\beta$ will be defined at the beginning of the proof of Proposition 4.7. Hence, we define the* ***blue neighborhood*** *$B(\lambda) := B'(\lambda) \cup \Delta_\alpha \cup \Delta_\beta$. Also, we define $B^* := B'(\lambda) \cup \partial\lambda$.*

*The grand orbit of $\bigcup_{\lambda \in \Lambda} Y(\lambda)$ is denoted as $\mathcal{Y}$, and it is where the deformation occurs. Also we define $\mathcal{Y}^* := \mathcal{Y} \cup \mathcal{L} \cup \partial\mathcal{L}$.*

*We define $\mathcal{B}$ to be an open neighborhood of $\mathcal{Y}$ with the property that on $\mathcal{B} \setminus \mathcal{Y}$, $\sigma_t$ is conformal, $f(\mathcal{B}) \subseteq \mathcal{B}$, and such that it contains the grand orbit of $B(\lambda)$. Similarly, we define $\mathcal{B}^*$ as the grand orbit of $B'(\lambda)$.*

Note that $\mathcal{Y} = \mathcal{V}$ as sets, but they are images of different functions from where are built. We remark that the components of $\mathcal{B}$ contain disjoint grand leaves of $\mathcal{L}$. See Figure 2.12 to see these neighborhoods.

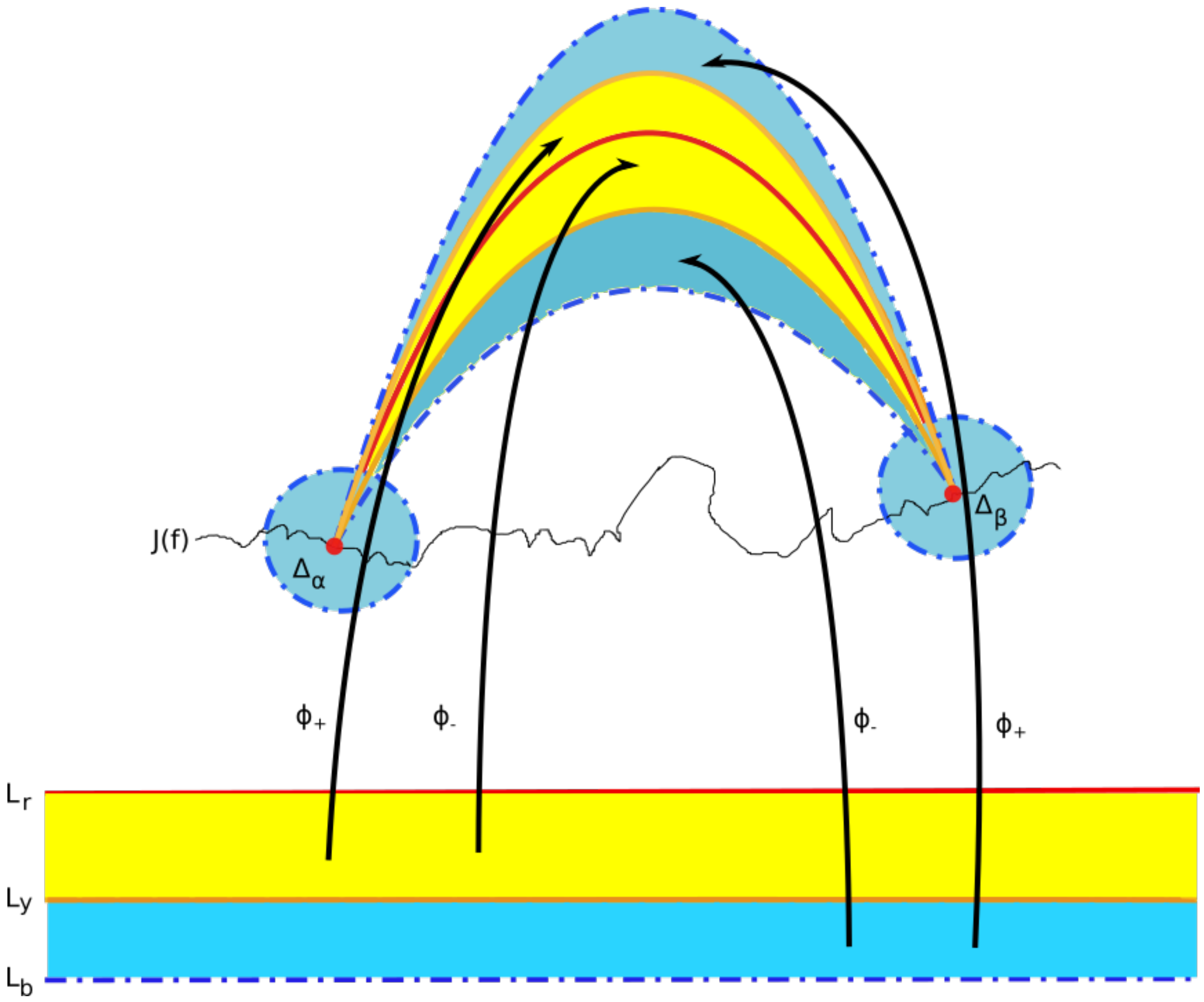

**Figure 2.12:** The yellow and the blue strip neighborhoods.

Now we have all the elements to define a pinching deformation formally. See [*Robles & Sienra 2022*].

DEFINITION 2.5 [*Robles & Sienra 2022*] *Let $f$ be an entire transcendental function with at least one cycle of periodic Baker domains $U = \{U_0, U_1, \ldots, U_{p-1}\}$ with a Baker lamination $\Lambda$ in $U$. The family of almost complex structures $(\sigma_t)_{t\in[0,1)}$ defines a* ***pinching deformation of*** *$f$ with support in $\mathcal{V}$ contained in $\mathcal{B}$. These structures come with quasiconformal maps $h_t : \overline{\mathbb{C}} \to \overline{\mathbb{C}}$ via integration by Theorem 2.2 (Ahlfors-Bers-Morrey), that can be normalized assuming $h_t$ fixes $\infty$ and two points $p, q \in J(f)$. The map $f_t := h_t \circ f \circ h_t^{-1}$ is holomorphic for every $t \in [0, 1)$.*

*We say that a* ***pinching deformation converges uniformly*** *if $h_t \rightrightarrows H$ (double arrow means uniform convergence on euclidean metric) and the nontrivial fibers of $H$ are the grand orbit $\mathcal{L}$. In the sense that $diam_s(h_t(\bar{\gamma})) \to 0$ as $t \to 1$, for each*

$\gamma \in \mathcal{L}$.

See Figure 2.13 for a visualization of the pinching deformation at time $t$. Consider that we have put two bands together in the domain and codomain of $\tilde{P}_t$ where, in the mathematical sense, the band in the lower part is not well drawn, but they help visualize the situation.

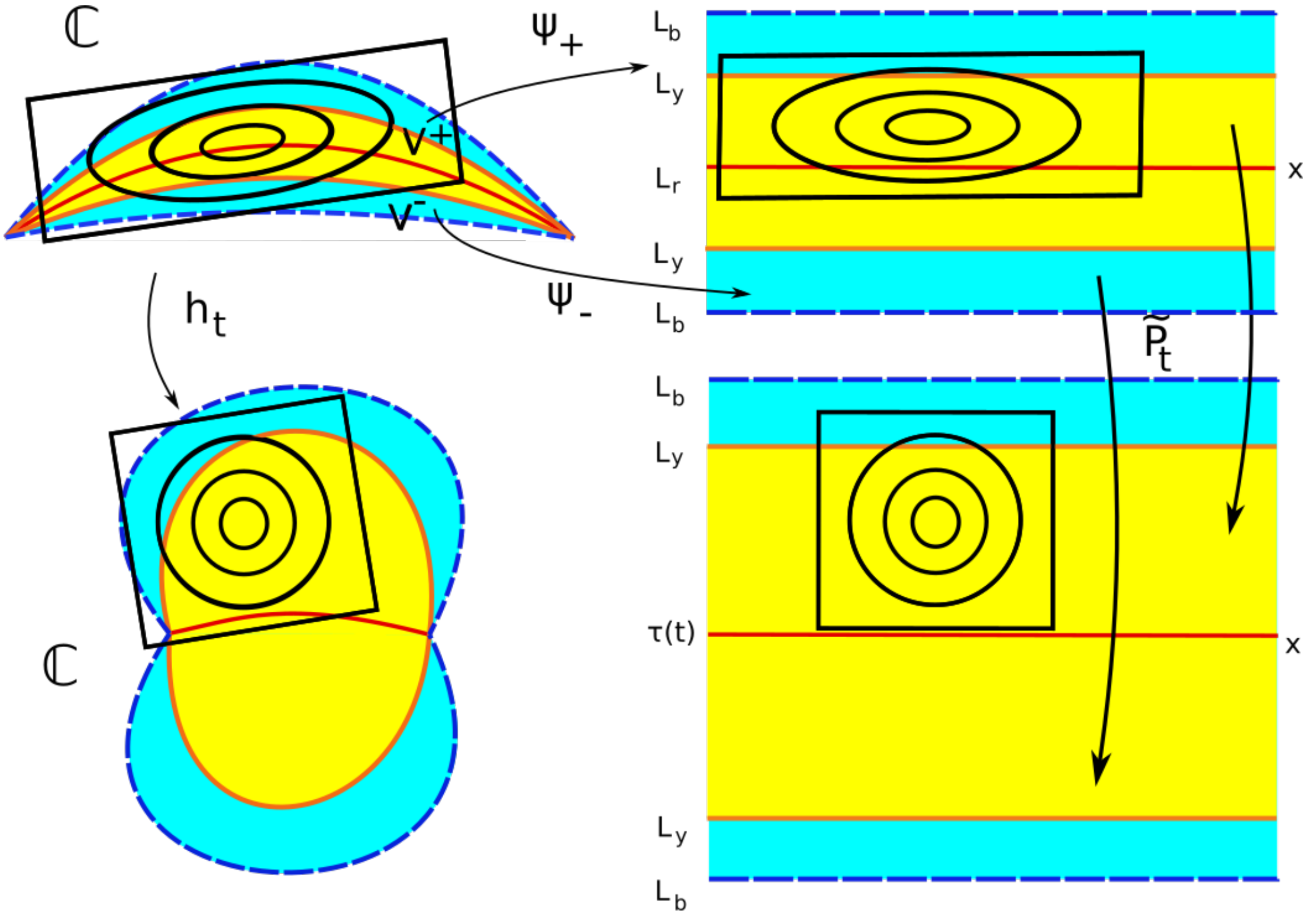

**Figure 2.13:** The pinching deformation at time $t$.

By Lemma 2.1, the composition $h_t \circ f \circ h_t^{-1}$ defines an entire transcendental function. We are interested in showing the circumstances when pinching Baker laminations converge uniformly via $f_t$, i.e. when $f_t \rightrightarrows F$, where $F$ is an entire transcendental function. In this context we have the next lemma, proved in Appendix A of [*Haissinsky & Tan, 2004*]:

LEMMA 2.2 *Let $f : \mathbb{C} \to \mathbb{C}$ be a continuous surjective map. For $t \in [0, 1)$, let $G_t, R_t : \overline{\mathbb{C}} \to \overline{\mathbb{C}}$ be two families of homeomorphisms of $\overline{\mathbb{C}}$. Assume that, as $t \to 1$, $G_t$ and $H_t$ converge uniformly to continuous maps $G, R$, respectively, and $f$ maps each fiber of $G$ into a fiber of $R$. Then $f_t := R_t \circ f \circ \left(G_t^{-1}\right) : \mathbb{C} \to \mathbb{C}$ converges uniformly to a continuous map $F$ and $F \circ G = R \circ f$.*

Taking $F_t = G_t := h_t$ in Lemma 2.2, it implies that if $h_t \rightrightarrows H$, then $f_t \rightrightarrows F$, and this is the reason why in Definition 2.5 we are interested only in the uniform convergence of $h_t$. It guarantees the uniform convergence of $f_t$.

Finally, when we have a uniform convergent pinching deformation the next two lemmas establish the shrinking of the grand lamination $\mathcal{L}$ in a Baker domain $U$ under the following circumstances.

The following lemma is a slight generalization of [*Dominguez & Sienra 201*5] of Lemma 2.1 taken from [*Tan 2002*] and shows that, in the limit, the spherical diameter of the pinched curves shrinks to a point.

*LEMMA 2.3 Let $f$ be an entire transcendental function with a Baker lamination $\Lambda$ in a Baker domain $U$ of $f$. Let $(\sigma_t)_{t\in[0,1)}$ be a uniform convergent pinching deformation of $f$. Then, for any grand leaf $\gamma \in \mathcal{L}$,*

$$\lim_{t\to 1} diam_s h_t\left(\bar{\gamma}\right) = 0.$$

Proof.

This lemma is a direct consequence of Definition 2.5 and formalizes the comment at the end of it.

Let $h_t$ be the integrating maps of Definition 2.5 such that $h_t \rightrightarrows H$ when $t \to 1$. By Theorem 2.2 (Ahlfors-Bers-Morrey) the pullback of the measurable field of circles under the mappings $\widetilde{P}_t \circ \psi_\pm$ is the same as the pullback of the measurable field of circles under $h_t$.

By definition, $\phi_\pm := \Psi \circ M \circ exp \circ S_\pm$, where $\Psi, M, exp$, and $S_\pm$ are conformal mappings. Hence, the pullback of the measurable field of circles under their inverse branches $\psi_\pm$ is a measurable field of circles. And the actual deformation is provided by the pullback under $\widetilde{P}_t$.

Now we recall Remark 9 by [*Robles & Sienra 2022*]. From property 3 in the definition of $\widetilde{P}_t$, its Beltrami coefficient is

$$\frac{\partial_{\bar{z}}\widetilde{P}_t}{\partial_z\widetilde{P}_t}(x+iy) = \frac{1-\partial_y v_t(y)}{1+\partial_y v_t(y)}.$$

By definition of quasiconformal dilatation in Section 2.2, we have that $K\left(\widetilde{P}_t\right) = \frac{1}{\partial_y v_t(y)}$. Thus $K\left(h_t\right) = K\left(\widetilde{P}_t \circ \psi_\pm\right) = K\left(\widetilde{P}_t\right) = \frac{1}{\partial_y v_t(y)}$, which has support in the grand orbit of $\mathcal{Y} = \bigcup_{\lambda \in \Lambda} Y(\lambda)$. Since $\partial_y v_t(y) \to 0$ as $y \to L_r$ and $t \to 1$, then $K(h_t) \to \infty$.

From the definition of quasiconformal maps $h_t$ by the almost complex structures $\sigma_t$, and for a quadrilateral $Q$ that intersects a grand leaf $\gamma \in \mathcal{L}$ we obtain that $K(h_t(Q)) \to \infty$. Therefore, for any $z_1, z_2 \in \gamma$, by taking a quadrilateral with two opposite sides intersecting the segment $\lambda$ in $z_1$ and $z_2$, we have that $d_s(h_t(z_1), h_t(z_2)) \to 0$ as $t \to 1$. See the left-hand side of Figure 2.13. The uniform convergence of $h_t$ uniformizes this convergence to 0.

Then, for any grand leaf $\gamma \in \mathcal{L}$, $\lim_{t \to 1} diam_s h_t(\bar{\gamma}) = 0$.

■

As observed by [*Robles & Sienra 2022*] in Remark 10, from property 3 in the definition of $\widetilde{P}_t$, we have that the coefficient of Beltrami of $\widetilde{P}_t$ is locally uniformly bounded away from 1 at any point $(t, x, y) \neq (1, x, L_r)$, therefore following the argument of the proof of Lemma 2.3, we have for any quadrilateral $Q$ disjoint of the grand orbit $\mathcal{L}$ for a Baker lamination $\mathit{\Lambda}$ in a Baker domain $U$, we have that $K(h_t(Q))$ is bounded away from $\infty$ and so the quadrilateral $h_t(Q)$ does not shrink to a point in the process of a pinching deformation. This implies that if $D$ is any disc in $\mathbb{C} \setminus \mathcal{L}$, then $h_t(D)$ is homeomorphic to $D$ for $t \in [0, 1)$. When $t = 1$, some problems may appear (see the proof and the remark of Theorem A).

As a direct consequence of the theorem of Mañé due to semihyperbolicity on an entire transcendental function $f$, i.e. see Theorem 2.3, we have the next lemma, see [*Haissinsky & Tan, 2004*].

LEMMA 2.4 *Let $f$ be an entire transcendental semihyperbolic function with a Baker domain $U$ and a Baker lamination $\Lambda$ such that $sing(f^{-1}) \cap \mathcal{L} = \emptyset$, where $\mathcal{L}$ is the grand orbit of $\Lambda$. Thus the diameter of any sequence of elements of $\mathcal{V} := \bigcup_{k \in \mathbb{N}} f^{-k}(V_\delta(\Lambda))$ tends to 0 when $k \to \infty$ if $\lambda_{a,\infty} \notin \Lambda$.*

Proof.

This lemma is a direct consequence of the Shrinking Lemma for rational maps proved by [*Lyubich & Minsky 1997*], pp. 86. Their proof is based on the fact that a rational map has at most $k$ singular points. Since in entire transcendental functions the set $sing(f^{-1})$ could be infinite is sufficient to ask for $sing(f^{-1}) \cap \mathcal{L} = \emptyset$. With this hypothesis, we can build a good neighborhood for $\mathcal{L}$ disjoint to $sing(f^{-1})$ and the same proof by [*Lyubich & Minsky 1997*] is valid for $\mathcal{V}$.

■

Note that we are interested to do a pinching deformation of an entire transcendental function $f$ on a Baker lamination $\Lambda$ with support in $\mathcal{V}$, and by Definitions 2.4 and 2.5, we can build good neighborhoods on $V_{\delta}(\Lambda)$ such that $\mathcal{V} \cap sing(f^{-1}) = \emptyset$.

# 3 Baker Domains and Divergent Deformations

Pinching deformation of some completely invariant Baker domains was studied in [*Dominguez & Sienra 2015*]. In particular, it is proven that the Fatou function $e^{-z} + z + 1$ can be pinched to the Baker-Dominguez function $e^{-z} + z$. The Fatou function has a completely invariant doubly parabolic Baker domain that contains the right half plane. When it is pinched to the Baker-Dominguez function, the pinching deformation $e^{-z} + z$ has infinite invariant doubly parabolic Baker domains (see Figure 3.1).

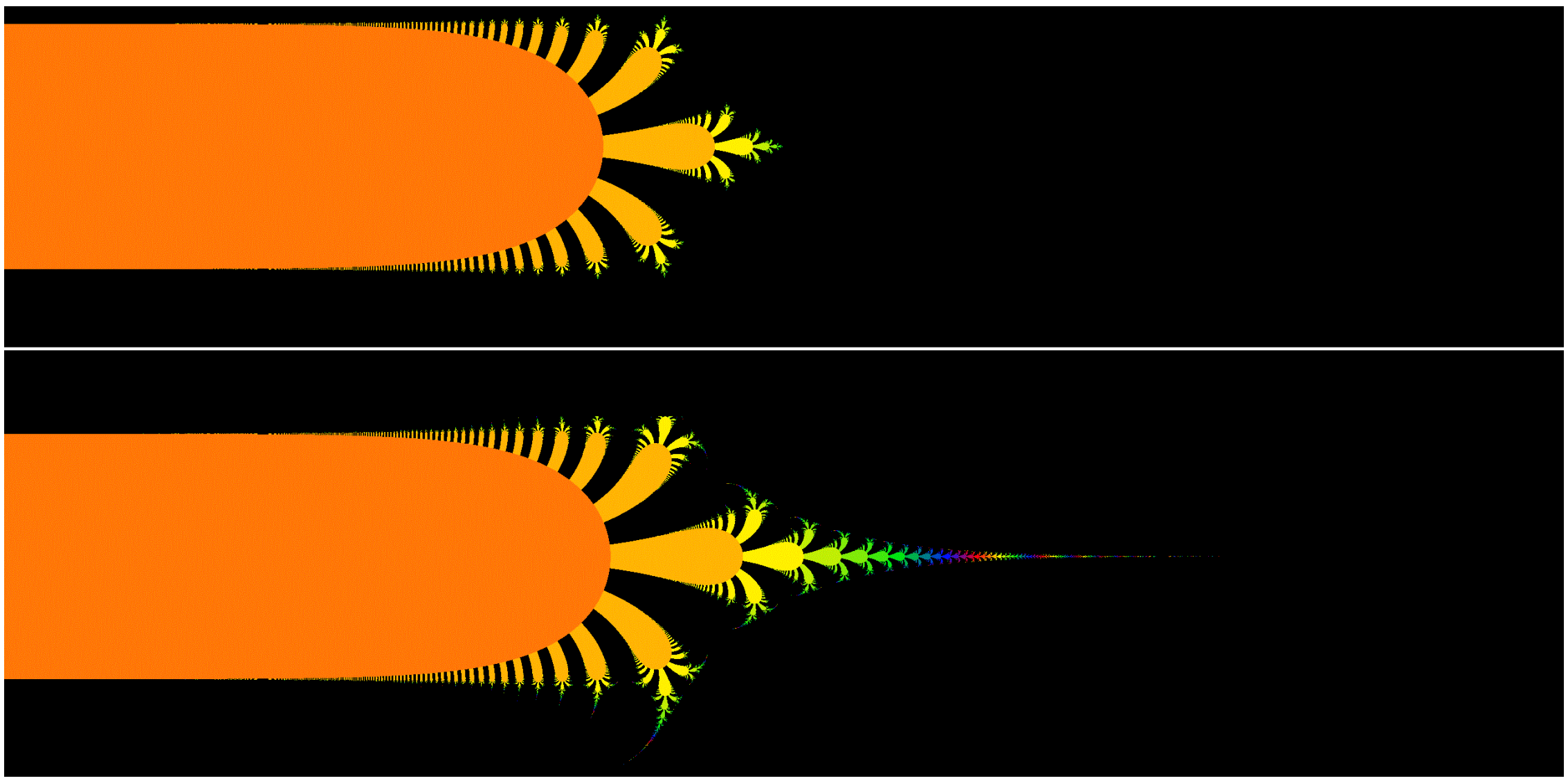

**Figure 3.1:** Pinching $e^{-z} + z + 1$ (above) to $e^{-z} + z$ (below).

In the next theorem, we consider the case when the Baker domain is not completely invariant. We prove that for certain curves the pinching process along such curves is divergent, see [*Robles & Sienra 2022*]:

THEOREM A [*Robles & Sienra 2022*] *Let $f : \mathbb{C} \to \mathbb{C}$ be an entire transcendental function with a non-completely invariant Baker domain $U$. Consider a Baker lamination $\Lambda$ in $U$ with a leaf $\lambda_{z_0,\infty}$ having endpoints at $z_0 \in \mathbb{C}$ and $\infty$, with $z_0$ a non-exceptional point in $\partial U$. Thus, the pinching deformation along the* grand orbit *$\mathcal{L}$ does not converge uniformly.*

Proof.

As $z_0$ is a non-exceptional point, $\overline{\bigcup_{n=1}^{\infty} f^{-n}(z_0)} = J(f)$ and so there is a subsequence $\{z_{n_k}\} := \{f^{-n_k}(z_0)\} \to z_0$ as $n_k \to \infty$, this convergence is understood in the spherical metric. Then, we have a family of curves

$$\{\gamma_{n_k}\} := \left\{f^{-n_k}(\lambda_{z_0,\infty})\right\} \subset \mathcal{L} \subset \left\{f^{-n_k}(U)\right\}$$

with $\{z_{n_k}\}$ and $\infty$ as endpoints of each one.

Since $U$ is not completely invariant we have that the curves $\{\gamma_{n_k}\} \subset \overline{\mathbb{C}} \setminus U$ are disjoint in $\mathbb{C}$, see [*Bergweiler & Eremenko 2007*]. Thus, for $\varepsilon_0$ there is a natural number $N_{\varepsilon_0}$ such that for every $n_k > N_{\varepsilon_0}$ the subsequence $\{z_{n_k}\} \subset D_{\varepsilon_0}(z_0) \cap \overline{\mathbb{C}} \setminus U$. Notice that $\infty$ is accessible from $\overline{\mathbb{C}} \setminus h_t(U)$ by the curves $h_t(\{\gamma_{n_k}\})$, (all these curves are attached to $\infty$), see Figure 3.2.

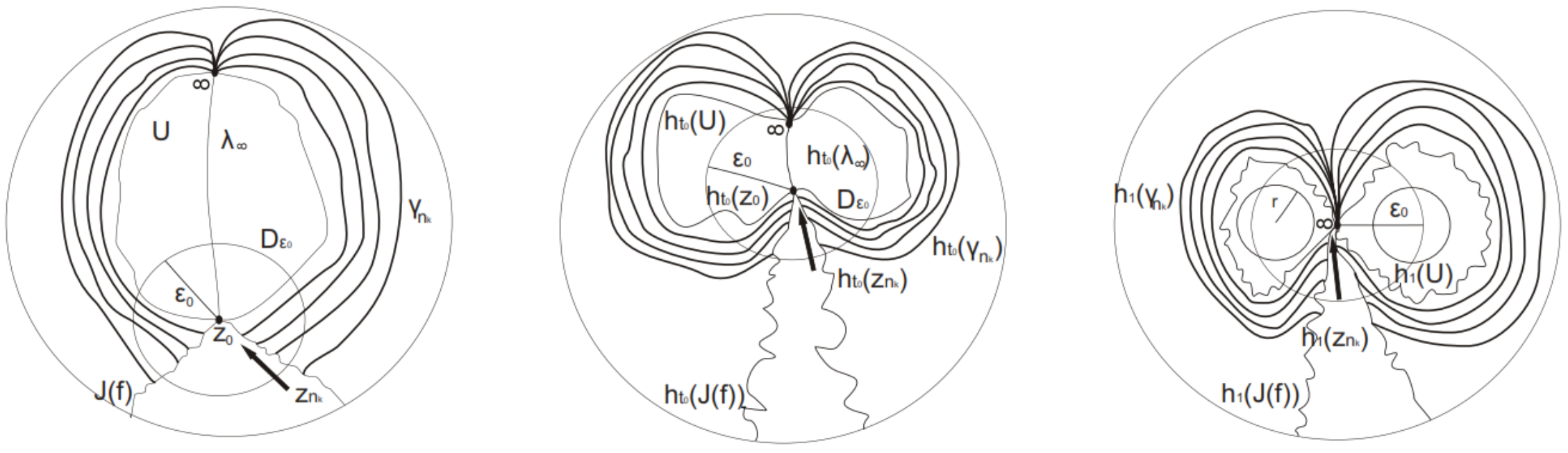

**Figure 3.2:** The pinching process in the proof of Lemma 3.1.

Notice that the duality between $\lambda_{z_0,\infty}$ being in $U$ and $\{\gamma_{n_k}\}$ not being in $U$ is the heart of the problem of convergence.

As a consequence of Lemma 2.3, we can assume that the pinching deformation along $\mathcal{L}$ converges uniformly via the quasiconformal maps $h_t$. This implies $h_t \rightrightarrows H$ as in Definition 2.5. Then, for $\gamma \in \mathcal{L}$, $diam_s(h_t(\bar{\gamma})) \to 0$ as $t \to 1$. Since $h_t$ fixes infinity for all $t$, then $h_t(\gamma)$ tends to infinity. In particular $h_t(\gamma_{n_k}) \to \infty$ as $t \to 1$.

Notice that the set $C_t := \bigcup_{k=0}^{\infty} \left(h_t\left(\gamma_{n_k}\right) \cup h_t\left(\lambda_{z_0,\infty}\right)\right)$ disconnects the complex plane in at least two connected components $\Omega_i$, such that $\Omega_i \cap U \neq \emptyset$, with $i = 1, 2$. We have two cases: either some region $h_t\left(\Omega_i\right)$ collapses to $\infty$ as $t \to 1$, i.e., $h_1\left(\Omega_i\right) = \infty$ or none of the regions collapses. By hypothesis we have that $f_t$ converges uniformly to an entire transcendental function $g$ and the non-trivial fibers of the pinching deformation are the grand leaves of the grand orbit $\mathcal{L}$. Hence, the family $h_t\left(\Omega_i\right)$ can not collapse to a point when $t \to 1$, with $i = 1, 2$. In the Remark below, we explain why the requirement for the fibers, in this situation.

By continuity of $h_t$, $\lim_{n_k \to \infty} h_t\left(z_{n_k}\right) = h_t\left(z_0\right)$ for every $t \in [0, 1]$. It follows that $\left\{h_t\left(z_{n_k}\right)\right\} \subset D_{\varepsilon_0}\left(h_t\left(z_0\right)\right) \cap \left(\overline{\mathbb{C}} \setminus h_t\left(U\right)\right)$, and since the pinching deformation is convergent, $\lim_{t \to 1} h_t\left(z_0\right) = h_1\left(z_0\right) = \infty$.

Observe that the regions $h_t\left(\Omega_i\right)$ contain an open set for $t \in [0, 1]$ and do not collapse, then there are two open discs, one at each side of $h_t\left(\lambda_{z_0,\infty}\right)$ and contained in $h_t\left(U \setminus \Lambda\right)$ of radius $r_t > 0$, such that $0 < r_t < diam_s\left(\left\{h_t\left(\gamma_{n_k}\right)\right\}\right)$ for $n_k > N\left(\varepsilon_0\right)$, for all $0 \leq t \leq 1$.

Observe that there exists $t_0$ such that $d_s\left(\infty, h_{t_0}\left(z_0\right)\right) < \varepsilon_0$, therefore $D_{\varepsilon_0}\left(h_t\left(z_0\right)\right) \cap \overline{\mathbb{C}} \setminus h_t\left(U\right)$ has two components for $t > t_0$. One component contains the endpoints $\left\{h_t\left(z_{n_k}\right)\right\}$, the other component contains the access to $\infty$ from $\overline{\mathbb{C}} \setminus h_t\left(U\right)$.

This implies that for every curve in $\left\{h_t\left(\gamma_{n_k}\right)\right\}$, its intersection with $D_{\varepsilon_0}\left(h_t\left(z_0\right)\right)$ has two components. But this is a contradiction because the convergence of the pinching deformation implies that $diam_s\left(\left\{h_t\left(\gamma_{n_k}\right)\right\}\right) \to 0$ when $t \to 1$. Thus the pinching deformation along $\mathcal{L}$ does not converge uniformly.

■

*Remark on the proof.* Observe that for every $t \in [0, 1]$, we have $h_t\left(\Omega_i\right) \cap J\left(f_t\right) \neq \emptyset$ with $i = 1, 2$. If $p \in J(f)$, then $p_t := h_t(p) \in J\left(f_t\right)$. Additionally, by Montel´s theorem, there is a $m \geq 0$ such that for $V_{p_t}$ any neighborhood of $p_t$, $f_t^m\left(V_{p_t}\right) \cap h_t\left(\Omega_i\right) \neq \emptyset$, where the integer $m$ depends on $V_{p_t}$. Therefore

$$V_{p_t} \cap f_t^{-m}\left(h_t\left(\Omega_i\right)\right) \neq \emptyset.$$

Assume that the functions $f_t$ converge uniformly to an entire function $g$ and $h_t\left(\Omega_i\right)$ collapses to $\infty$ as $t \to 1$. Then, there exists $p \in J(f) \setminus \Omega_i$, such that $p_t \in J\left(f_t\right) \setminus h_t\left(\Omega_i\right)$ for $t \in [0, 1]$, otherwise $J(g) = \infty$. Therefore, $p_1 \in J(g)$ and for any neighborhood $V_{p_1}$ of $p_1$, there is an inverse branch of $\infty$ in $V_{p_1}$. This implies that $p_1$ is either a prepole or the accumulation point of different preimages of $\infty$. Moreover, $p_1$ is an essential singularity, so $g$ is not an entire function, contradicting the hypothesis.

□

Now, as a non-trivial example consider the Bergweiler function $f(z) = 2 - log(2) + 2z - e^z$. This has a Baker domain $U$ of hyperbolic type I not completely invariant and contains the half plane $\{z \mid Re(z) < 2\}$. A lamination $\Lambda$ on $U$ with a leaf $\lambda_{z_0,\infty} \in \Lambda$ consists of the set $(-\infty, \zeta) \times \{0\} \subseteq \mathbb{C}$, where $z_0 = \zeta$ is a repelling fixed point of $f$ in $J(f)$ and the grand orbit $\mathcal{L} = \bigcup\limits_{n\in\mathbb{N}} f^{-n}(\lambda_{z_0,\infty})$. From Lemma 3.1, the pinching deformation of $f$ along $\mathcal{L}$ does not converge uniformly.

Note that in the previous example, the core geodesic of the cylinder $U/f$ is pinched and the limit surface exists, as a noded surface, see [*Bers 1974*], but the limit function does not.

On the other hand, there is a possibility that a Baker lamination intersects the set of asymptotic values of $f$, as a consequence we have the following Corollary.

COROLLARY A [*Robles & Sienra 2022*] *Let $f : \mathbb{C} \to \mathbb{C}$ be an entire transcendental function with a non-completely invariant Baker domain $U$. Consider a Baker lamination $\Lambda$, with a leaf $\lambda_{a,b}$ having endpoints at non-exceptional points $a, b \in \mathbb{C}$. If $\lambda_{a,b}$ intersects the set of asymptotic values of $f$, then the pinching deformation along the grand orbit of the lamination, $\mathcal{L}$, does not converge uniformly.*

Proof.

Let $\lambda_a \in \mathcal{L}$ be a grand leaf that intersects the set of asymptotic values of $f$. Then there is a grand leaf $\sigma \in \mathcal{L}$ with $f(\sigma) = \lambda_a$ such that $\sigma$ is in some component of the inverse image of $U$. Moreover, $\sigma$ has one of its extreme points at $\infty$. Hence, by Theorem $A$, the pinching deformation along $\mathcal{L}$ does not converge uniformly.

■

Observe that Theorem A and Corollary A can be generalized to a cycle of non-completely invariant Baker domains $U = \{U_0, U_1, \ldots, U_{p-1}\}$ with $f^p(U_i) \subseteq U_i$ for each $i \in \{0, 1, \ldots, p-1\}$. If for some $i$, we have a leaf $\lambda_{a,\infty} \in \Lambda \cap U_i$, with $z_0$ a non-exceptional point, or a leaf $\lambda_{a,b} \in \Lambda \cap U_i$ intersecting the asymptotic values of $f$, thus we define $g := f^p$, where $g(\lambda_{a,\infty}) = \lambda_{g(a),\infty} \in \Lambda \cap U_i$, and $g(\lambda_{a,b}) = \lambda_{g(a),g(b)} \in \Lambda \cap U_i$ by Definition 2.2, and the latter intersects the asymptotic values of $g$. Here, $U_i$ is a non-completely invariant Baker domain for $g$. Thus we apply Theorem A or Corollary A to the function $g$ in the respective case for each leaf.

Note that quasiconformal deformations of a map are a class of quasiconformal maps, this is to avoid trivial situations. For instance, in case that $h_t(z)$ converges, the quasiconformal maps $\tilde{h}_t(z) := \dfrac{h_t(z)}{1-t}$ does not converge uniformly, even though it integrates the same structure. However, we showed in Theorem A and Corollary A, that the pinching deformation does not converge uniformly, no matter which integrating map $h_t$ is chosen for each $t$.

# 4 Convergent Deformations on Baker Domains

In this chapter, we will prove the next theorem.

THEOREM C *Let* $f : \mathbb{C} \to \mathbb{C}$ *be a semihyperbolic entire transcendental function with a Baker domain* $U$ *that satisfies property* $P$ *and with* $J(f)$ *thin at* $\infty$. *Let* $\mathcal{L}$ *be a grand orbit of a Baker lamination* $\Lambda$ *in* $U$ *that does not contain a leaf of the type* $\lambda_{a,\infty}$ *and* $sing\,(f^{-1}) \cap \mathcal{L} = \emptyset$, *then there exists a uniformly convergent continuous pinching deformation* $f_t = h_t \circ f \circ h_t^{-1}$ *to an entire function* $F$. *The mappings* $h_t$ *are quasiconformal mappings that converge uniformly to a map* $H$, *whose non-trivial fibers are the* $\mathcal{L}-$ *components.*

Due to the extensive length of the proof, this will be divided into sections and we will follow [*Haissinsky & Tan, 2004*] closely because it is the same proof in many points. The idea of the proof is the following, we will prove that the family of quasiconformal mappings $h_t$ is equicontinuous, therefore is normal and there exists a subsequence converging uniformly to a map $H$, from there we will see that $H$ is the unique map to which $\{h_t\}$ converge, and thus all the $\{h_t\}$ converge uniformly to $H$, and so $\{f_t\}$ converges uniformly to $F$ where $H \circ f = F \circ H$.

A family of functions on metric spaces $\mathcal{F} = \{f : A \to (Y, d_Y) \mid A \subseteq (X, d_X)\}$ is *equicontinuous in* $z_0 \in A$, if for any $\varepsilon > 0$ there exists $\delta(\varepsilon) > 0$ such that if $z \in A$ satisfies $d_X(z, z_0) < \delta$, then $d_Y(f(z), f(z_0)) < \varepsilon$ for any $f \in \mathcal{F}$.

To realize the proof of Theorem C we prove the next theorem:

THEOREM B *Let* $f : \mathbb{C} \to \mathbb{C}$ *be a semihyperbolic entire transcendental function with a cycle of Baker domains* $U$. *Let* $\mathcal{L}$ *be a Baker lamination in* $U$ *that does not contain a leaf of the type* $\lambda_{a,\infty}$ *and* $sing\,(f^{-1}) \cap \mathcal{L} = \emptyset$, *then the maps* $\{h_t\}$ *that integrate the family of almost complex structures* $(\sigma_t)_{t\in[0,1)}$ *are equicontinuous in* $\overline{\mathbb{C}}$. *Furthermore, for any subsequence* $\{h_{t_k}\}$ *converging to a map* $H$, *the nontrivial fibers of* $H$ *are exactly the components of* $\mathcal{L}$.

Due to the length of the proof of this theorem, we prove it in the next five subsections:

- 4.1 Equicontinuity of $\{h_t\}$ in $\overline{\mathbb{C}} \setminus (J(f) \cup \mathcal{L})$.
- 4.2 Equicontinuity of $\{h_t\}$ in $J(f) \setminus \{\infty\}$.
- 4.3 Equicontinuity of $\{h_t\}$ at $\infty$.
- 4.4 Equicontinuity of $\{h_t\}$ in $\mathcal{L}$.
- 4.5 The nontrivial fibers of $H$ are the grand leaves in $\mathcal{L}$.

## 4.1 Equicontinuity of $\{h_t\}$ in $\overline{\mathbb{C}} \setminus (J(f) \cup \mathcal{L})$

PROPOSITION 4.1. *The maps $\{h_t\}$ are equicontinuous $\overline{\mathbb{C}} \setminus (J(f) \cup \mathcal{L})$, where $\lambda_{a,\infty} \notin \mathcal{L}$.*

Proof.

As explained in the proof of Lemma 2.3, by Theorem 2.2 (Ahlfors-Bers-Morrey), the Beltrami coefficient of $h_t$ is the same of $\widetilde{P}_t$. And we have two cases:

a) Let $z_0 \in \overline{\mathbb{C}} \setminus (J(f) \cup \mathcal{L})$ such that $\psi_{\pm}(z_0) = x_0 + iy_0$, with $y_0 \in [L_b, L_y]$.

By construction in Section 2.3, $\widetilde{P}_t(x_0 + iy_0) = x_0 + iv_t(y_0) = x_0 + iy_0$. By Section 2.3, we have that

$$\mu_{h_t}(z_0) = \frac{\partial_{\bar{z}} h_t}{\partial_z h_t}(z_0) = \frac{\partial_{\bar{z}} \widetilde{P}_t}{\partial_z \widetilde{P}_t}(x_0 + iy_0) = \frac{1 - \partial_y v_t(y_0)}{1 + \partial_y v_t(y_0)} = 0.$$

Then by Weyl´s lemma, $h_t$ is holomorphic in $z_0$ for every $t \in [0, 1)$. Since

$$\left|\partial_z \left(\widetilde{P}_t \circ \psi_{\pm}\right)(z_0)\right| = |\partial_z \psi_{\pm}(z_0)| \leq M\text{'},$$

this implies that $|\partial_z h_t(z_0)| \leq M$ for each $t \in [0, 1)$.

Let $\varepsilon > 0$, we choose $\delta = \varepsilon/M$, for every $z \in B_\delta(z_0)$ we have that integrating along the line segment from $z_0$ to $z$:

$$|h_t(z) - h_t(z_0)| = \left|\int_{z_0}^{z} \partial_z h_t(w)dw\right| \leq M|z - z_0| < M\delta = \varepsilon.$$

Therefore, $h_t$ is equicontinuous in $z_0$.

b) Let $z_0 \in \overline{\mathbb{C}} \setminus (J(f) \cup \mathcal{L})$ such that $\psi_\pm(z_0) = x_0 + iy_0$ with $y_0 \in [L_y, L_r)$. By the technical assumption in Section 2.3, for any $L' < L_r$, there is $t(L') \in (0,1)$ with $t(L') \to 1$ as $L' \to L_r$ such that for any $(s,y) \in (t(L'), 1] \times [L_b, L']$, we have $v_s(y) = v_{t(L')}(y)$. See Figure 2.11. Thus for every $t > t(L')$

$$\mu_{h_t}(z_0) = \frac{\partial_{\bar{z}} h_t}{\partial_z h_t}(z_0) = \frac{\partial_{\bar{z}} \widetilde{P}_t}{\partial_z \widetilde{P}_t}(x_0 + iy_0) = \frac{1 - \partial_y v_t(y_0)}{1 + \partial_y v_t(y_0)} = \frac{1 - \partial_y v_{t(L')}(y_0)}{1 + \partial_y v_{t(L')}(y_0)} = k_0,$$

and the Beltrami coefficient of $h_t$ is bounded by some $k_0 < 1$. As above,

$$\left|\partial_z \left(\widetilde{P}_t \circ \psi_\pm\right)(z_0)\right| \leq M",$$

which implies that $|\partial_z h_t(z_0)| \leq M^*$ for each $t \in [0,1)$. Since $h_t$ is real-differentiable:

$$h_t(z) = D_{h_t}(z_0)(z - z_0) + h_t(z_0) + e_t(z),$$

where $D_{h_t}(z_0)(z - z_0) = \partial_z h_t(z_0)(z - z_0) + \partial_{\bar{z}} h_t(z_0)(\bar{z} - \overline{z_0})$ and $\lim_{z \to z_0} |e_t(z)| = 0$.

Thus

$$|h_t(z) - h_t(z_0)| = |D_{h_t}(z_0)(z - z_0) + e_t(z)| =$$

$$\leq |D_{h_t}(z_0)(z - z_0)| + |e_t(z)| =$$
$$= \left|\partial_z h_t(z_0)\left[(z - z_0) + \frac{\partial_{\bar{z}} h_t(z_0)}{\partial_z h_t(z_0)}(\bar{z} - \overline{z_0})\right]\right| + |e_t(z)| \leq$$
$$\leq |M^*[(z - z_0) + k_0(\bar{z} - \overline{z_0})]| + |e_t(z)| \leq$$
$$\leq M^*|z - z_0| + M^* k_0 |\bar{z} - \overline{z_0}| + |e_t(z)| =$$
$$= M^*(1 + k_0)|z - z_0| + |e_t(z)|.$$

Now, let $\varepsilon > 0$. Since $\lim_{z \to z_0} |e_t(z)| = 0$, for $\frac{\varepsilon}{2}$ there exists $\delta_{e_t} > 0$ such that if $0 < |z - z_0| < \delta_{e_t}$ it holds that $|e_t(z)| < \frac{\varepsilon}{2}$.

Let $\delta = inf\left\{\frac{\varepsilon}{2M^*(1 + k_0)}, \delta_{e_t} | t \in [0,1)\right\}$, where $inf\{\delta_{e_t} | t \in [0,1)\} > 0$ by construction, then for any $z$ satisfying $0 < |z - z_0| < \delta$, then

$$M^*(1 + k_0)|z - z_0| < \frac{\varepsilon}{2},$$

$$\Rightarrow |h_t(z) - h_t(z_0)| \leq M^*(1+k_0)|z - z_0| + |e_t(z)| < \frac{\varepsilon}{2} + \frac{\varepsilon}{2} = \varepsilon.$$

Therefore, by a) and b), $\{h_t\}$ is equicontinuous in $z_0 \in \overline{\mathbb{C}} \setminus (J(f) \cup \mathcal{L})$.

■

## 4.2 Equicontinuity of $\{h_t\}$ in $J(f) \setminus \{\infty\}$

To prove equicontinuity for the remaining cases we need the next criteria of equicontinuity at one point:

LEMMA 4.2 *Let* $\mathcal{A} := \{h : \mathbb{D} \to \mathbb{C}\}$ *be a family of continuous injective maps such that* $\bigcup_{h\in\mathcal{A}} h(\mathbb{D})$ *avoids at least two points* $a, b \in \mathbb{C}$.

1. *Let* $(U_n)_{n\geq 0}$ *be a nested sequence of disk-like neighborhoods of the origin in the unit disk* $\mathbb{D}$ *such that* $A'_n = \mathbb{D} \setminus \overline{U_n}$ *is an annulus. If there exists an increasing subsequence* $\eta_n \to \infty$ *such that* $\forall h \in \mathcal{A}, \forall n \geq 0$, *then* $mod(h(A'_n)) \geq \eta_n$, *therefore* $\mathcal{A}$ *is equicontinuous at the origin.*

2. *Let* $A_n \subset \mathbb{D}$ *be a nested sequence of annuli such that, for all* $n$, $A_{n+1}$ *is contained in the component of* $\mathbb{C} \setminus \overline{A_n}$ *containing* $0$. *If there is* $M > 0$ *such that* $\forall h \in \mathcal{A}, \forall n \geq 0$, *then* $mod(h(A_n)) \geq M$, *therefore* $\mathcal{A}$ *is equicontinuous at the origin.*

Proof.

1. Let $\varepsilon > 0$. There exists $n \in \mathbb{N}$ such that $\frac{\pi^2}{2\varepsilon^2} < \eta_n \leq mod(h(A'_n))$ for $U_n$ and the remaining $h \in \mathcal{A}$. Let $\delta > 0$ such that the neighborhood $V_\delta(0) \subseteq U_n$, and let $z \in V_\delta(0) \setminus \{0\}$. Since $a, b \notin h(\mathbb{D})$, $h(A'_n)$ separates the points $h(0), h(z)$ from the points $a, b$. Therefore we use Lemma 6.1 from [*Lehto & Virtanen 1973*]: "If the ring domain $B$ separates the pair of points $a_1, b_1$, from the pair $a_2, b_2$, and if the spherical distance satisfies $d(a_i, b_i) \geq \alpha > 0$, $i = 1, 2$, then the module $M(B) \leq \frac{\pi^2}{2\alpha^2}$."

Let $\alpha := min(d_s(a,b), d_s(h(0), h(z)))$ then $mod(h(A'_n)) \leq \frac{\pi^2}{2\alpha^2}$, which implies $\frac{\pi^2}{2\varepsilon^2} < \frac{\pi^2}{2\alpha^2}$ and $\alpha < \varepsilon$. If $\varepsilon$ is small enough such that $0 < \varepsilon < d_s(a,b)$ thus

$$\alpha = d_s(h(0), h(z)) < \varepsilon,$$

for every $z \in V_\delta(0)$ and for the remaining $h \in \mathcal{A}$ taken above. Therefore, $\mathcal{A}$ is equicontinuous in 0.

2. By superadditivity of the module, as can be seen in Lemma 6.3 from [*Lehto & Virtanen 1973*]: "Let $A, A_1, A_2, \ldots$ be ring domains such that $A_1, A_2, \ldots$ are disjoint (among them) subdomains of $A$. If every $A_n$ separates $(-A)_1$ and $(-A)_2$ (the components of the complement of $A$), then $\sum_n M(A_n) \leq M(A)$." See Figure 4.1. In the notation of this lemma, let

$$A = B_n := \left(\bigcup_{k=1}^{n} A_k\right) \bigcup \left(\bigcup_{k=0}^{n} \overline{C_k}\right) \setminus \partial\mathbb{D},$$

where $C_k$ is the annular component between $A_k$ and $A_{k+1}$, and $C_0$ is the component between $\partial\mathbb{D}$ and $A_1$. So, $A_k \subseteq B_n$ for every $k \leq n$, and using Lemma 6.3 from [*Lehto & Virtanen 1973*] for every $A = B_n$:

$$nM \leq \sum_{k=1}^{n} mod\left(h\left(A_k\right)\right) \leq mod\left(h\left(B_n\right)\right),$$

which is the case 1, where $A'_n = B_n$.

■

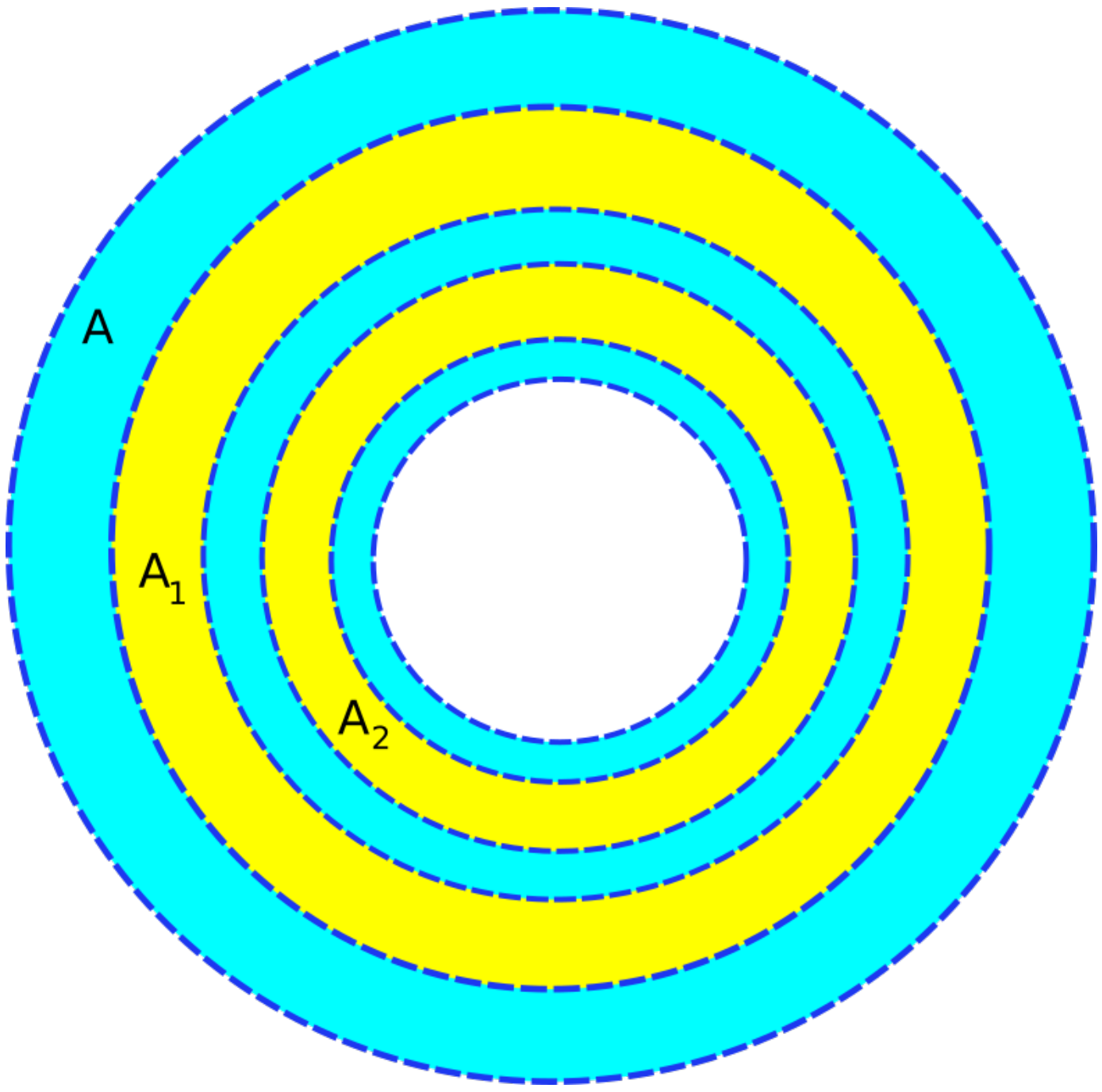

**Figure 4.1:** Lemma 6.3 from [*Lehto & Virtanen 1973*].

### 4.2.1 Control of the Moduli of Deformed Annuli

For the next lemma, check the above Definition 2.4 for the meaning of $\mathcal{B}$.

LEMMA 4.3 *Let $A \subset \mathbb{C}$ be a bounded annulus such that $\partial A \cap \mathcal{B} = \emptyset$. Then there is $m > 0$ (depending on $A$ but not on $t$) such that $mod\,(h_t(A)) \geq m$ for all $t$.*

To prove this lemma we need the next two lemmas.

KEY LEMMA: *There is a uniform constant $0 < c \leq 1$ with the following properties. Let $\eta : [0,1] \to \mathbb{C}$ (respectively $\eta : [0,1] \to \overline{\mathbb{C}}$) be a rectifiable curve with end points outside $\mathcal{B}$. Then*

$$l_{\varrho_e}(\eta) \geq c d_e\,(\eta(0), \eta(1)) \qquad \text{(resp. } l_{\varrho}(\eta) \geq c d\,(\eta(0), \eta(1))\text{)},$$

*where $d_e$ is the Euclidean metric (resp. $d$ the spherical metric), $\varrho_e$ (resp. $\varrho$) is the same metric but with zero density in $\mathcal{Y}$ (the support of $\sigma_t$), i.e.*

$$\varrho_e(z)\,|dz| = (1 - \chi_{\mathcal{Y}}(z))\,|dz|\,, \qquad \varrho(z)\,|dz| = \frac{1 - \chi_{\mathcal{Y}}(z)}{1 + |z|^2}\,|dz|\,,$$

*where* $\chi_{\mathcal{Y}}(z) = \begin{cases} 1 & \text{if } z \in \mathcal{Y}. \\ 0 & \text{if } z \notin \mathcal{Y}. \end{cases}$

Proof.

The proof is almost the same as [*Haissinsky & Tan, 2004*], with some modifications for this case. Let $I \subseteq [0,1]$ be a maximal open subinterval such that $\eta(I) \subset \mathcal{B}$, and we can analize any component of $\mathcal{B}$, i.e., for some blue strip $B(\lambda)$ where $\eta(I) \subset B(\lambda)$ and $\eta(\partial I) \subset \partial B(\lambda)$. Let us define $\eta'(I)$ as the geodesic joining the two ends of $\eta(I)$. So, we are going to prove that

$$l_{\varrho_e}(\eta\,(I)) \geq c l_e\,(\eta'\,(I)) \qquad \text{(resp. } l_{\varrho}(\eta) \geq c l\,(\eta'\,(I))\text{)},$$

for some $c > 0$. See Figure 4.2.

Case 1. If $\eta(I) \cap \mathcal{Y} = \emptyset$, then by definition of geodesic curve $l_{\varrho_e}(\eta\,(I)) = l_e(\eta\,(I)) \geq l_e(\eta'\,(I))$.

Case 2. If $\eta(I) \cap \mathcal{Y} \neq \emptyset$, let $I = (x_1, x_2)$, which give us two subintervals $I_1 = (x_1, x_1')$ and $I_2 = (x_2', x_2)$ satisfying $\eta\,(x_i) \in \partial B(\lambda)$ and $\eta\,(x_i') \in \partial\,(\mathcal{Y} \cap B(\lambda))$. We will show that there are constants $c_B, c > 0$, such that

$$
\begin{aligned}
l_{\varrho_e}(\eta(I)) = l_{\varrho_e}(\eta(I_1)) + l_{\varrho_e}(\eta(I_2)) = & \\
= l_e(\eta(I_1)) + l_e(\eta(I_2)) \geq & \quad \text{by definition of } \varrho_e(z)\,|dz|, \\
\geq |\eta(x_1) - \eta(x_1')| + |\eta(x_2) - \eta(x_2')| \geq & \quad \text{by definition of geodesic curve,} \\
\overset{*}{\geq} c_B\left[|\eta(x_1) - \partial\lambda_1| + |\eta(x_2) - \partial\lambda_1|\right] \geq & \quad \text{where } \partial\lambda_1 \text{is one point in } \partial\lambda, \\
\geq c_B\,|\eta(x_1) - \eta(x_2)| = & \quad \text{by triangle inequality} \\
= c_B l_e(\eta'(I)) \overset{**}{\geq} c l_e(\eta'(I)) & \quad \text{by definition of geodesic curve.}
\end{aligned}
$$

The inequalities marked by $*$ and $**$ are commented on below.

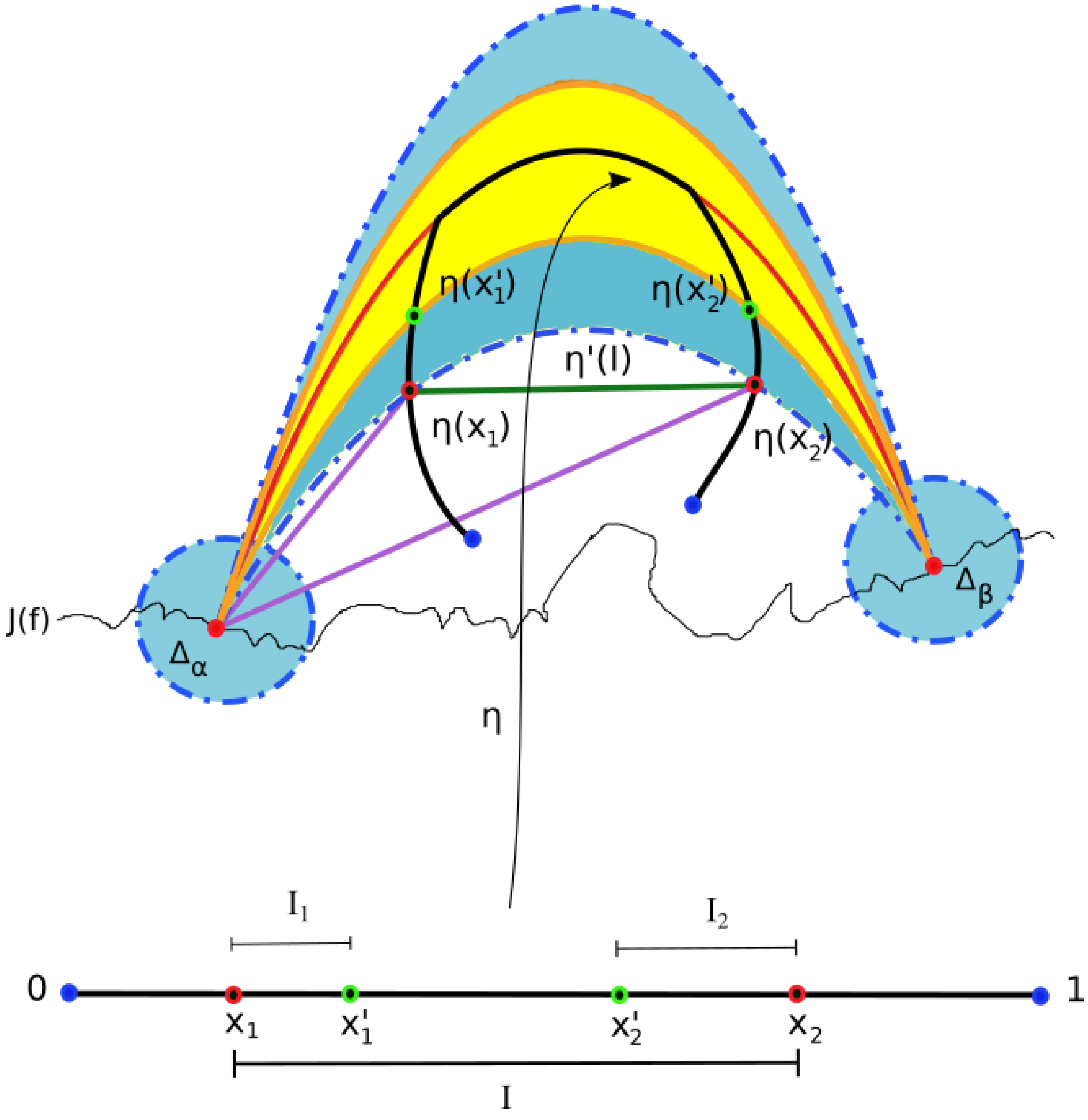

**Figure 4.2:** Proof of Lemma 4.3.

Proof of $*$. To prove it, we need the next Lemma which is the same version given by [*Haissinsky & Tan, 2004*].

LEMMA 4.4 *If* $g : \mathbb{D} \to \mathbb{C}$ *is a univalent function such that* $g(z) = \lambda z + O(z)$ *at the origin with* $|\lambda| > 1$, *and if* $\gamma_1$ *and* $\gamma_2$ *are two disjoint invariant arcs in* $\mathbb{D} \setminus \{0\}$ *landing at* $0$, *then* $\gamma_1 \cup \gamma_2$ *forms a "quasi-arc" in the following sense: there is a constant* $c_B > 0$, *depending only on the germ* $g$, *such that for any* $z \in \gamma_2$, $d_e(z, \gamma_1) \geq c_B d_e(z, 0)$.

To apply this lemma to our case, we need to observe that in the process of the pinching deformation, if $z_0$ is an endpoint of a leaf, it could happen that $f(z_0) \neq z_0$, but we apply the proof of the Lemma 4.4 given by [*Haissinsky & Tan, 2004*] to the function $g(z) = f(z) - f(z_0)$ to adjust it, i.e., there is no real change when $z_0$ is not being fixed. What we need is for $f$ to be semihyperbolic to apply theorem 2.3, that is, that $f$ is expansive in its Julia set.

To prove $**$, the arguments of [*Haissinsky & Tan, 2004*] using cross-ratio, are still valid in the same way as commented for Lemma 4.4.

□

Proof of Lemma 4.3. is a consequence of Key Lemma and there is no modification for this case, the proof given by [*Haissinsky & Tan, 2004*] using extremal length inequalities is the same.

■

### 4.2.2 One Good Annulus around each Julia Point

LEMMA 4.5 *Fix* $r > 0$ *(this will be the constant for the definition of semihyperbolicity).*

*1. For any* $z \in J(f) \setminus \overline{\mathcal{L}}$, *there are two open neighborhoods* $N'(z)$ *and* $N(z)$ *of* $z$ *in* $D_{r/4}(z)$ *and* $m > 0$ *such that mod* $\left(h_t\left(N(z) \setminus \overline{N'(z)}\right)\right) \geq m$ *for all* $t \in [0, 1)$.

*2. For any* $z \in \partial\lambda$, *for some* $\lambda \in \mathcal{L}$, *and* $B(\lambda)$ *the corresponding* $\mathcal{B}$ *− component, there are two open neighborhoods* $N'(\lambda)$ *and* $N(\lambda)$ *of* $\lambda$ *in* $\left(D_{r/4}(z) \cup F(f)\right) \cap B(\lambda)$ *and* $m > 0$ *such that mod* $\left(h_t\left(N(\lambda) \setminus \overline{N'(\lambda)}\right)\right) \geq m$ *for all* $t$.

Proof of point *1.* of Lemma 4.5 is the same by [*Haissinsky & Tan, 2004*], but we will detail one inequality through Lemma I.6.3 by [*Lehto & Virtanen 1973*] in * below, which states that for two annulus $A, B$ such that $A \subseteq B$ then $mod\, A \leq mod\, B$:

Let $z \in J(f) \setminus \mathcal{L}$. We choose simply connected neighborhoods with Jordan curves as boundaries $\overline{N'}, N$ such that $\overline{N'} \subset N \subset D_{r/4}(z)$, with no $\mathcal{B}^* - component$ would have a closure that intersects both boundaries of $\overline{N'}$ and $N$, and with the sets in the contention below being simply connected. Then

$$\overline{N'} \cup \left( \bigcup_{\substack{B(\lambda)\cap\partial N' \neq \emptyset \\ B(\lambda)\in\mathcal{B}}} \overline{B(\lambda)} \cup \partial N' \right) \subset N \setminus \left( \bigcup_{\substack{B(\lambda)\cap\partial N \neq \emptyset \\ B(\lambda)\in\mathcal{B}}} \overline{B(\lambda)} \cup \partial N \right)$$

with $\lambda \in \mathcal{L}$. See Figure 4.3.

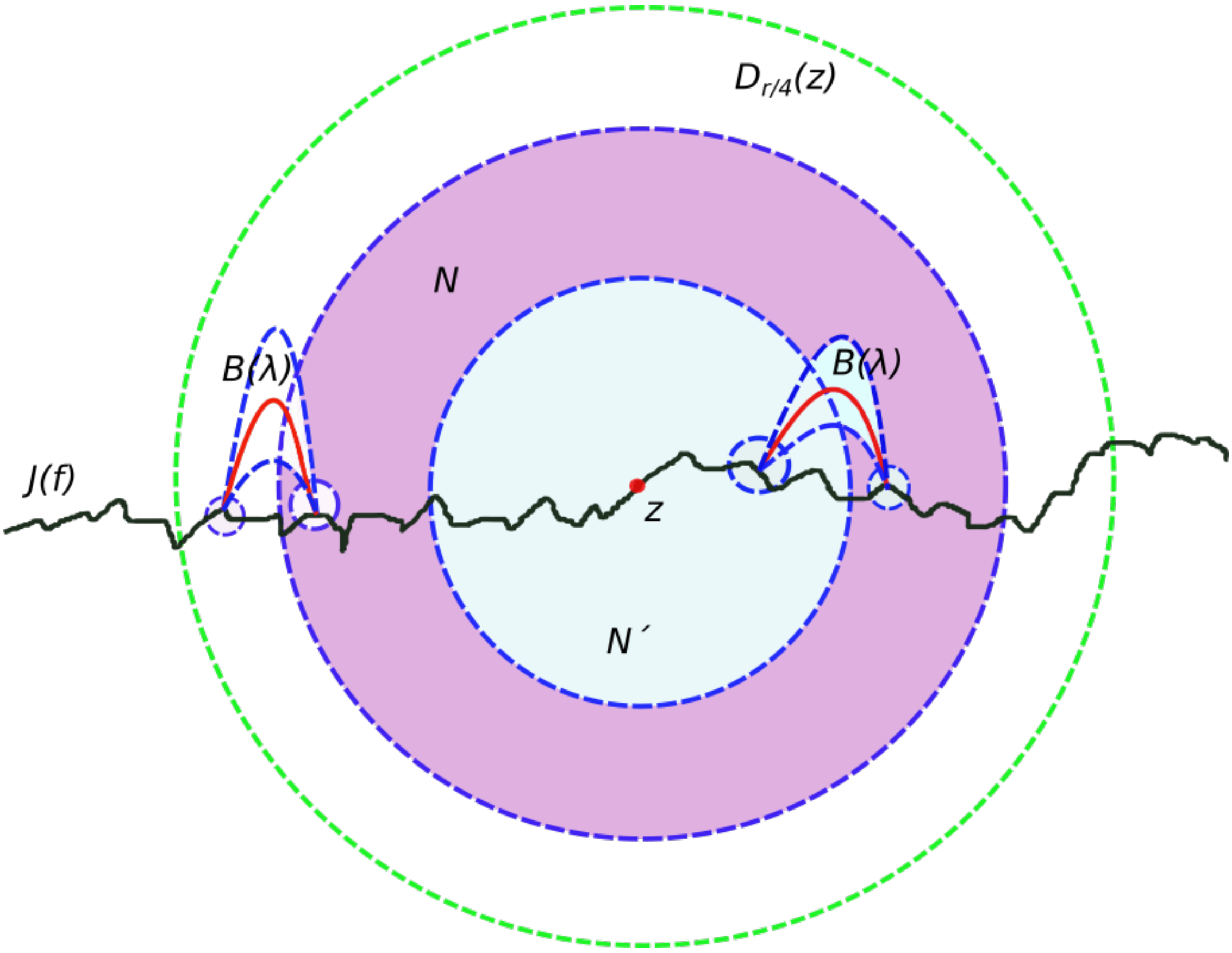

**Figure 4.3:** Proof of Lemma 4.4.

The left-hand set is compact and the right-hand set is open, so their difference is an annulus and it contains an annular component $A$ satisfying $\partial A \cap \mathcal{B} = \emptyset$.

Thus for Lemma 4.3, there is $m > 0$ such that

$$mod\, h_t\left(N \setminus \overline{N'}\right) \overset{*}{\geq} mod\, h_t(A) \overset{\text{Lemma 4.3}}{\geq} m.$$

For the second part, let $z \in \partial\lambda$ for some $\lambda \in \mathcal{L}$. We choose $N$ and $N'$ similarly, but as neighborhoods of $B(\lambda)$ and as subsets of $(\mathcal{B} \cap F(f)) \cup D_{r/4}(\partial\lambda)$. And we apply the same inequalities.

■

PROPOSITION 4.2 *Equicontinuity of* $\{h_t\}$ *at* $z_0 \in J(f) \setminus \{\infty\}$.

Proof.

Following [*Haissinsky & Tan, 2004*]. Let us normalize $f(z)$ such that $\infty$ is interchanged by some point of $F(f)$, and we denote this point as $\infty_{\mathbb{C}} \in \mathbb{C}$, thus $J(f)$ is a compact set.

The sets $N'(z)$ from above, such that $z \in J(f) \setminus \{\infty_{\mathbb{C}}\}$, define an open cover on $J(f)$. We extract a finite subcover $N'(z_i)$, with $i = 1, ..., l$. By Lemma 4.5 there exists $m > 0$ such that, for any $t \in [0, 1)$ and any $i \in \{1, ..., l\}$,

$$mod\, h_t\left(N(z_i) \setminus \overline{N'(z_i)}\right) \geq m.$$

Let $z_0 \in J(f) \setminus \{\infty_{\mathbb{C}}\}$. By semihyperbolicity, there are infinitely many $n$ such that $f^n$ maps a neighborhood $U$ of $z_0$ to $D_r(z_0)$ with a degree at most $\delta$. And, there is $i(n)$ such that $f^n(z_0) \in N'\left(z_{i(n)}\right)$. Taking a subsequence $f^{n_k}(z_0)$ if necessary by compacity, we may assume $i(n_k) \equiv i$ for some $i \in \{1, ..., l\}$.

Then we have two cases. First case, $z_i \in \partial\lambda$ for some $\lambda \in \mathcal{L}$; then we choose $N(z_i)$ contained in a $\mathcal{B}-$ *component*, denoted by $B(\lambda)$, and by semihyperbolicity, $deg\left(f^{n_k}|_{B(\lambda)}\right) \leq \delta$.

Second case, $z_i \notin \partial\lambda$, thus

$$f^{n_k}(z_0) \in N'(z_i) \subset N(z_i) \subset D_{r/4}(z_i) \subset D_{r/2}(f^{n_k}(z_0)) \subset D_r(f^{n_k}(z_0)),$$

for $n_k \in \mathbb{N}$. See Figure 4.4.

Let $E_{n_k}$ and $U_{n_k}$, be the respective components of $f^{-n_k}\left(\overline{N_i'}\right)$ and $f^{-n_k}(N_i)$ contain-ing $z_0$, where $E_{n_k}\subset U_{n_k}$. Let $A_{n_k}:=U_{n_k}\setminus E_{n_k}$. Then

$$mod\, h_t\left(A_{n_k}\right)\overset{*}{\geq}\frac{1}{\delta}mod\, h_t\left(N\left(z_i\right)\setminus\overline{N'\left(z_i\right)}\right)\overset{\text{Lemma 4.5}}{\geq}\frac{m}{\delta}.$$

The proof of the *-inequality is in the proof of Lemma 2.1 by [*Shishikura & Tan 2000*].

Taking again a subsequence if necessary, we may assume that the annuli $A_{n_k}$ are disjoint, and they are nesting down to $z_0$, this by Theorem 2.3 and by the annulus $D_r\left(f^{n_k}\left(z_0\right)\right)\setminus D_{r/2}\left(f^{n_k}\left(z_0\right)\right)$. By Lemma 4.2, $\{h_t\}$ is equicontinuous at $z_0$.

■

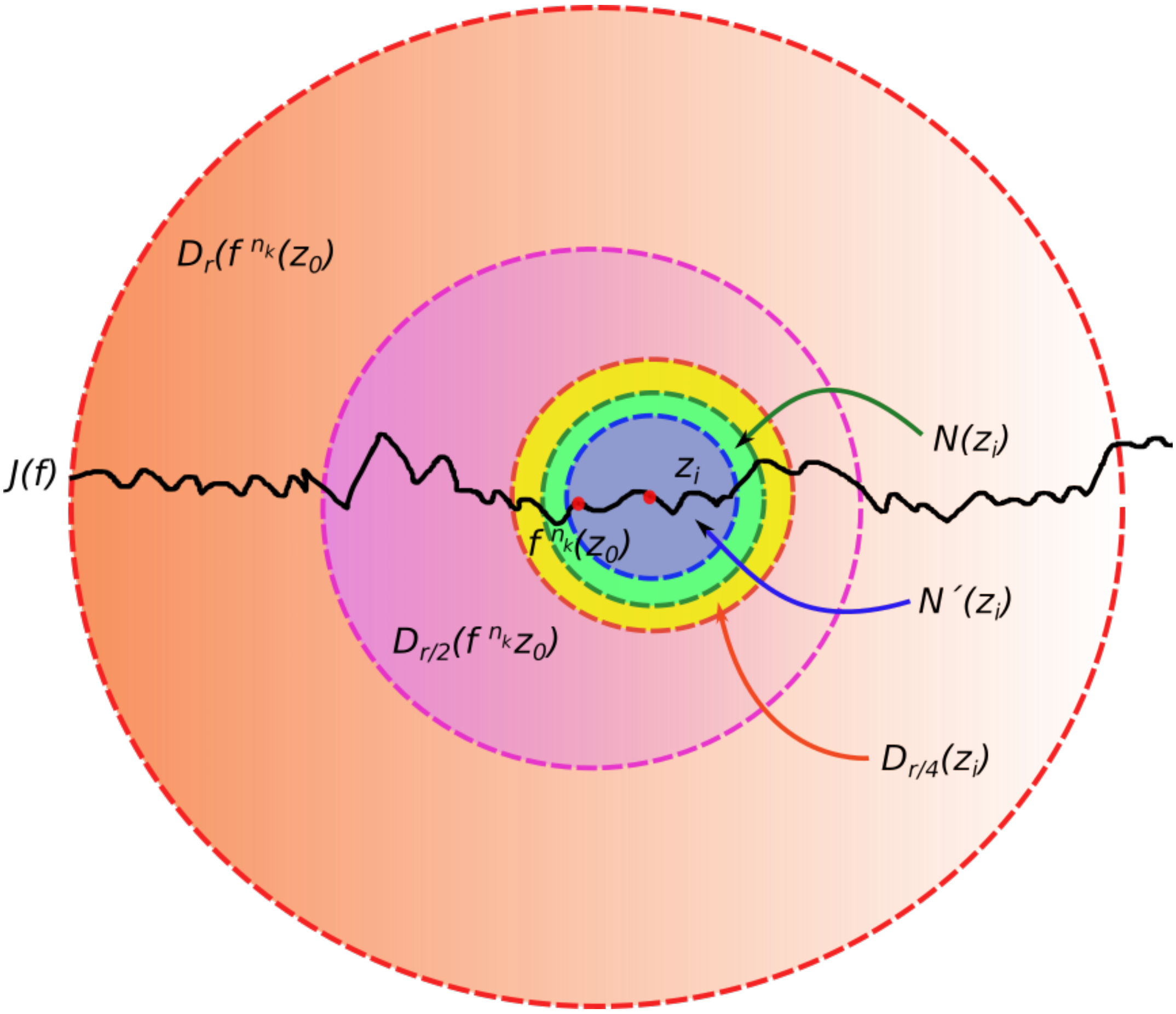

**Figure 4.4:** Proof of Proposition 4.2.

## 4.3 Equicontinuity of $\{h_t\}$ at $\infty$

PROPOSITION 4.3 *The sequence $\{h_t\}$ is equicontinuous at $\infty$.*

Proof.

As we noted in section 2.2, $f$ is not semihyperbolic at $\infty$. In like manner of the proof of Proposition 4.2, let us normalize $f(z)$ such that $\infty$ is interchanged by some point of $F(f)$, and we denote this point as $\infty_{\mathbb{C}} \in \mathbb{C}$. Since there is no leaf with $\infty$ as an endpoint in the lamination, i.e., $\lambda_{a,\infty} \notin \mathcal{L}$, thus we can choose by Corollary 2.1 simply connected neighborhoods about $\infty_{\mathbb{C}}$, with Jordan curves as boundaries of $\overline{N'}, N$ such that $\overline{N'} \subset N$ and with no $\mathcal{B}^* - component$ would have a closure that intersect the boundaries of $\overline{N'}$ and $N$, and with the sets in the contention below being simply connected. Then, like in the proof of Lemma 4.5,

$$\overline{N'} \cup \left( \bigcup_{\substack{B(\lambda)\cap\partial N' \neq \emptyset \\ B(\lambda)\in\mathcal{B}}} \overline{B(\lambda)} \cup \partial N' \right) \subset N \setminus \left( \bigcup_{\substack{B(\lambda)\cap\partial N \neq \emptyset \\ B(\lambda)\in\mathcal{B}}} \overline{B(\lambda)} \cup \partial N \right),$$

with $\lambda \in \mathcal{L}$. The left-hand set is compact and the right-hand set is open, so their difference is an annulus and it contains an annular component $A$ satisfying $\partial A \cap \mathcal{B} = \emptyset$.

Thus, there is $m > 0$ such that

$$mod\, h_t\left(N \setminus \overline{N'}\right) \overset{*}{\geq} mod\, h_t(A) \overset{\text{Lemma 4.3}}{\geq} m.$$

In this case, we do not need the dynamics to build the nesting annuli around $\infty_{\mathbb{C}}$. As long as $\lambda_{a,\infty} \notin \mathcal{L}$ we can build, in the same way, the annuli $A_n$ nesting down at $\infty_{\mathbb{C}}$ such that $mod\, h_t(A_n) \overset{\text{Lemma 4.3}}{\geq} m$, for all $t \in [0, 1)$ and $n \in \mathbb{N}$. This is a direct consequence of the fact that the quasiconformal mappings $h_t$ are only deforming the grand orbit of the lamination $\mathcal{L}$ and the grand orbit $\mathcal{Y}$. The grand leaves of $\mathcal{L}$ may be accumulating at $\infty_{\mathbb{C}}$ since $\lambda_{a,\infty} \notin \mathcal{L}$ we can find the corresponding annuli $A_n$, where $h_t(A_n)$ is isotopic to the identity mapping $Id$.

By Lemma 4.2, $\{h_t\}$ is equicontinuous at $\infty_{\mathbb{C}}$.

■

## 4.4 Equicontinuity of $\{h_t\}$ in $\mathcal{L}$

PROPOSITION 4.4 The sequence $\{h_t\}$*is equicontinuous in* $\mathcal{L}$*, where* $\lambda_{a,\infty} \notin \mathcal{L}$*.*

Proof

Let $z_0 \in \mathcal{L}$.

We prove the proposition in two parts:

a) Let $z \in \lambda \subseteq \mathcal{L}\backslash\{\infty\}$ and let $\varepsilon > 0$, we are going to find the $\delta$ of the equicontinuity based on how $z$ can approach $z_0$.

We are taking some parts of the proof of Lemma 2.3. Let $h_t$ be the integrating maps of Definition 2.5, by Theorem 2.2 (Ahlfors-Bers-Morrey) the pullback of the measurable field of circles under the mappings $\widetilde{P}_t \circ \psi_\pm$ is the same as the pullback of the measurable field of circles under $h_t$.

By definition, $\phi_\pm := \Psi \circ M \circ exp \circ S_\pm$, where $\Psi, M, exp,$ and $S_\pm$ are conformal mappings. Hence, the pullback of the measurable field of circles under their inverse branches $\psi_\pm$ is a measurable field of circles. And the actual deformation is provided by the pullback under $\widetilde{P}_t$.

The map $\tilde{P}_t : \mathbb{R} \times [L_b, L_r] \to \mathbb{R} \times [L_b, \tau(t)]$ with $t \in [0,1]$ is defined by

$$\tilde{P}_t(x+iy) = x + iv_t(y),$$

with $\tau(t)$ an increasing function and $v_t(y)$ defined in the same section. This map $\tilde{P}_t$ is a non-decreasing mapping for constant $y$ and $t \in (0,1)$ by construction. Let $Q$ be a cuadrilateral intersecting a grand leaf $\lambda \in \mathcal{L}$, such that $z_0$ and $z$ are vertices in $Q$. Take $z_1 \in Q$ as another vertex such that $\psi_\pm(z_1) = x_1 + iy_1$.

By definition of quasiconformal dilatation in Section 2.2, we have that $K\left(\widetilde{P}_t\right) = \frac{1}{\partial_y v_t(y)}$. Thus $K(h_t) = K\left(\widetilde{P}_t \circ \psi_\pm\right) = K\left(\widetilde{P}_t\right) = \frac{1}{\partial_y v_t(y)}$, which has support in the grand orbit of $\mathcal{Y} = \bigcup_{\lambda\in\Lambda} Y(\lambda)$. Since $\partial_y v_t(y_1) \to 0$ as $y_1 \to L_r$ and $t \to 1$, then $K(h_t) \to \infty$.

Therefore, for $z_0, z \in \lambda$, being $Q$ a quadrilateral with two opposite sides intersecting the segment $\lambda$ in $z_0$ and $z$, we have that $d\left(h_t\left(z_0\right), h_t\left(z\right)\right) \to 0$ as $t \to 1$. Note that we have a euclidean diameter, since $\infty \notin \lambda$. See the left-hand side of Figure 2.13.

So the integrating map $h_t$ of $\tilde{P}_t \circ \psi_\pm$ resulting of the pullback of the field of circles is a decreasing contraction on $\lambda \in \mathcal{L} \setminus \{\infty\}$ as $t$ grows and $\lim_{t\to 1} diam_s h_t\left(\bar{\gamma}\right) = 0$.

Then for any $\varepsilon > 0$ there exists $t' \in [0, 1)$ such that if $t > t$' thus $|h_t(z) - h_t(z_0)| < \varepsilon$.

Now, by Lemma 2.4 by [*Haissinsky & Tan, 2004*]:

"Lemma 2.4. The family $\{h_t\}$ is uniformly equicontinuous if and only if it is pointwise equicontinuous, i.e., for any $z_0 \in \overline{\mathbb{C}}$ and any $\varepsilon > 0$, there exists $\eta > 0$ and $t_0 < 1$ such that for any $t \in [t_0, 1)$ and any $z$ with $d_s(z, z_0) \leq \eta$, we have that $d_s(h_t(z), h_t(z_0)) \leq \varepsilon$."

Choosing $\eta = diam_s(\lambda)$ and $t_0 = t'$, by this lemma, $\{h_t\}$ is equicontinuous in $z_0 \in \mathcal{L}$ when $z \in \lambda \in \mathcal{L}$.

b) Let $z \notin \mathcal{L} \setminus \{\infty\}$, such that $\psi_\pm(z) = x + iy$ with $y \in [L_y, L_r)$.

First, let $\gamma_{L_y} := \{z \in \mathbb{C} | \psi_\pm(z) = x + iL_y\}$. Note that $\gamma_{L_y} \subseteq \partial\mathcal{Y}$. By part *a*) of this proof, $h_t$ is a decreasing contraction on $\lambda \in \mathcal{L} \setminus \{\infty\}$ as $t$ grows, and $\underset{t\to 1}{lim} diam_s h_t(\bar{\gamma}) = 0$. By definition, for any $z \in \gamma_{L_y}$,

$$|z - z_0| \leq diam(Y(\lambda)) \leq A'$$

for some $A' \in \mathbb{R}$ since $Y(\lambda) \cap \{\infty\} = \infty$.

Similarly, let $\gamma_{L'} := \{z \in \mathbb{C} | \psi_\pm(z) = x + iL' \text{ and } L_y < L' < L_r\}$. Let

$$Y_{L'}(\lambda) = \phi_+(\mathbb{R} \times [L', L_r]) \cup \phi_-(\mathbb{R} \times [L', L_r]),$$

where $Y_{L'}(\lambda) \subset Y(\lambda)$ and for any $z \in \gamma_{L'}$,

$$|z - z_0| \leq diam(Y_{L'}(\lambda)) < diam(Y(\lambda)).$$

Since $h_t$ is a homeomorphism for every $t \in [0, 1)$, $h_t(Y_{L'}(\lambda)) \subset h_t(Y(\lambda))$ properly and

$$|h_t(z) - h_t(z_0)| < diam(h_t(Y_{L'}(\lambda))) < diam(h_t(Y(\lambda))) := A_t.$$

Since $\lambda$ is compact in $\mathbb{C}$, there exists $A := sup\{A_t | t \in [0, 1]\}$.

Second, the equicontinuity. Let $\varepsilon > 0$.

Since $h_t$ is continuous for every $t \in [0, 1)$, there exists $\delta_t(\varepsilon) > 0$ such that if $z \in B_\delta(z_0)$ it holds that $|h_t(z) - h_t(z_0)| < \varepsilon$.

For the equicontinuity of $\{h_t\}$, we need to find $\delta > 0$ independent of $t$. We can find a $\delta(A)$ depending on $A$ only. Because $A$ is bound of the pinching deformation on $Y(\lambda)$ where the integrating map $h_t$ is the result of integrating the Beltrami coefficient of $\widetilde{P}_t \circ \psi_\pm$ which was built differentiably as a differentiable isotopy as $t$ moves increasingly in $[0, 1)$, and since $\{h_t\}$ is a holomorphic motion and which is isotopic to the identity when $t \to 1$ in the complement $\overline{\mathbb{C}} \setminus \mathcal{L}$. See the Extended $\lambda$-Lemma by [*Astala & Martin 2001*].

■

*Remark:* By Theorem 2.2 (Ahlfors-Bers-Morrey), the integrating maps $h_t$ fix three points in $\overline{\mathbb{C}}$. Since $\infty$ is an essential singularity of entire transcendental functions, $h_t(\infty) = \infty$ for any $t \in [0, 1)$. We can normalize another fix point, let $h_t(z_t) = z_t$ for every $t \in [0, 1)$, let $g_t(z) := h_t(z + z_t) - z_t$ where $g_t(0) = 0$ and $g_t(\infty) = \infty$.

By Theorem 4.2 from [*Lehto & Virtanen 1973*]: "A family $\mathcal{H}$ of $K-$quasiconformal mappings of the domain $D$ is equicontinuous in $D$ if there exists $d > 0$ such that for three fixed points $a_1, a_2, a_3 \in D$ the (spherical) distances $d_s(h(a_i), h(a_j))$ for $i, j = 1, 2, 3$, $i \neq j$, are greater than $d$ for every mapping $h \in \mathcal{H}$".

Thus, it is enough to find another fix point $z_0$ under $g_t$ independent from $t \in [0, 1)$, to conclude that $\{g_t\}$ is equicontinuous, implying that $\{h_t\}$ is equicontinuous. We conjecture that this fix point $z_0$ can be found in some Fatou component in the case that $z \in \lambda \subseteq \mathcal{L} \setminus \{\infty\}$, due to the compactness of $\lambda$ in $\mathbb{C}$ and that the pinching deformation outside of $\mathcal{V}$ is a holomorphic motion and which is isotopic to the identity $z \mapsto z$ when $t \to 1$, and under some extra conditions, possibly. This would imply that the pinching deformation converges uniformly in a more general setting, and the pinching deformation would be a more effective tool to deform the dynamics of entire transcendental functions.

## 4.5 The Nontrivial Fibers of $H$ are Leaves in $\mathcal{L}$

In this section we will see that if $H$ is any convergent map from the equicontinuity of $\{h_t\}$, its nontrivial fibers are $\mathcal{L}-$*components*, i.e., the grand leaves on the lamination.

Let the sets $A, B, C \subset \overline{\mathbb{C}}$, we say that $A$ *separates* $B$ *from* $C$ if $B$ and $C$ are contained in different components of $\overline{\mathbb{C}} \setminus A$.

Following [*Gardiner 2000*], let $\Gamma$ be a set of curves on a Riemann surface $S$. Let $\gamma \in \Gamma$ be a countable union of open arcs or closed curves, the extremal length, $\Lambda(\Gamma)$, of $\Gamma$ is a sort of average minimum length of the curves in $\Gamma$. It is an important quantity because it is invariant under conformal maps and quasi-invariant under quasiconformal maps. To see this let us introduce the next concept.

A metric $\rho(z)\,|dz|$ is *allowable* if

- it is invariantly defined for different local parameters $z$, i.e., $\rho_1(z_1)\,|dz_1| = \rho_2(z_2)\,|dz_2|$, where $\rho_1$ and $\rho_2$ are representatives of $\rho$ in terms of the parameters $z_1$ and $z_2$;
- $\rho$ is locally $L_2$ and $\geq 0$ everywhere; and
- $A(\rho) = \iint_S \rho^2 dxdy \neq 0$ or $\infty$, where the integral is taken over the whole Riemann surface.

Due to the first condition, $A(\rho)$ is well defined. For such allowable metric $\rho$, we define

$$L_\gamma(\rho) := \int_\gamma \rho\,|dz|$$

if $\rho$ is measurable along $\gamma$; otherwise we define $L_\gamma(\rho) := \infty$. Let $L(\rho) = L(\rho, \Gamma) := \inf_{\gamma\in\Gamma} \{L_\gamma(\rho)\}$.

DEFINITION 4.1 *Let $\Gamma$ be a set of curves on a Riemann surface with $\gamma \in \Gamma$ be a countable union of open arcs or closed curves, the* ***extremal length of*** $\Gamma$ *is defined as*

$$\Lambda(\Gamma) := \sup_\rho \left\{ \frac{L(\rho)^2}{A(\rho)} \right\}$$

*where the supremum is taken over all allowable metrics.*

Note that the extremal length is invariant if $\rho$ is multiplied by a positive scalar, so we can scale $\rho$ in any way we want. The proof of the next lemma is in [*Gardiner 2000*].

LEMMA 4.6 *Suppose $f$ is a quasiconformal map with dilatation $K$ taking a Riemann surface $R$ onto a Riemann surface $R'$ and taking a set of curves $\Gamma$ onto a set of curves*

$\Gamma'$. *Then*

$$\frac{\Lambda(\Gamma)}{K} \leq \Lambda\left(\Gamma'\right) \leq K\Lambda(\Gamma).$$

In addition, we need the next lemma for uniform control of lengths, and it will also be used to prove the semihyperbolicity of limit functions $\{f_t\}$. Let $z \in \overline{\mathbb{C}}$, we assign to $z$ a compact subset $K(z)$ as follows. If $z \in \mathcal{B}^*$, there is a unique $\mathcal{B}^* - component$, $B^*$, and a unique $\mathcal{Y}^* - component$, $Y^*$, such that $z \in B^*$ and $Y^* \subset B^*$. We set $K(z) := Y^*$. If $z \notin \mathcal{B}^*$, we set $K(z) := \{z\}$.

Let $\mathcal{Q}$ be the set of distinct quadruples $q := (z_1, z_2, z_3, z_4)$, with $z_i \neq z_j$ if $i \neq j$. If $q \in \mathcal{Q}$, we define $\Gamma_q$ to be the set of rectifiable curves which separates $(z_1, z_2)$ from $(z_3, z_4)$.

For $r > 0$, we define $\mathcal{Q}_r \subset \mathcal{Q}$ as the set of quadruples such that $d\left(K\left(z_1\right), K\left(z_2\right)\right) \geq r$, $d\left(K\left(z_3\right), K\left(z_4\right)\right) \geq r$, $d\left(z_1, z_2\right) \geq r$ and $d\left(z_3, z_4\right) \geq r$.

LEMMA 4.7 *For all $r > 0$, there is a constant $l = l(r) > 0$ such that, for any $q \in \mathcal{Q}_r$ and any $\gamma \in \Gamma_q$, then $l_\varrho(\gamma) \geq l$, where $l_\varrho$ is given in the Key Lemma.*

The proof is exactly as Lemma 2.10 by [*Haissinsky & Tan, 2004*].

LEMMA 4.8 *Let $(z_1, z_2, z_3, z_4)$ be four different points such that no two belong to the same $\mathcal{L} - component$. Then, for $\Gamma$ the set of Jordan curves which separate $(z_1, z_2)$ from $(z_3, z_4)$ there is $m > 0$ such that $\Lambda\left(h_t\left(\Gamma\right)\right) \geq m$ for all $t \in [0, 1)$.*

The proof of Lemma 4.8 is the same as in Lemma 2.8 by [*Haissinsky & Tan 2004*] for this case. Henceforth the set of Jordan curves $\Gamma$ which separate $(z_1, z_2)$ from $(z_3, z_4)$ will be denoted as $\Gamma_{(z_1,z_2),(z_3,z_4)}$.

The appendix of [*Haissinsky & Tan, 2004*] states two results that we will use in the proof of Proposition 4.5. Lemma B.1 in [*Haissinsky & Tan, 2004*] shows that if $z, w \in \mathbb{C} \setminus \{0\}$ in such a way that $|z - w| < |w|$ and if $\Gamma$ is the set of rectifiable curves which separates $\{0, \infty\}$ from $\{z, w\}$, i.e., $\Gamma_{(0,\infty),(z,w)}$, then

$$|z - w| > |w|\, exp\left(\frac{-2\pi}{\Lambda(\Gamma_{(0,\infty),(z,w)})}\right).$$

As Corollary B.2 from this lemma we have that if $z_t, w_t, a, b \in \mathbb{C}$, are four distinct points with $z_t, w_t$ depending on a parameter $t$; if we assume that $d\left(w_t, \{a, b\}\right) \geq$ $C > 0$ for all $t$. And if $d(z_t, w_t) \to 0$ then $\Lambda\left(\Gamma_{(z_t,w_t),(a,b)}\right) \to 0$. As contrapositive we have that if $\Lambda\left(\Gamma_{(z_t,w_t),(a,b)}\right) \geq C' > 0$, then $d(z_t, w_t) \geq C'' > 0$.

PROPOSITION 4.5 *For any limit map $H$ of $\{h_t\}$, the nontrivial fibers are $\mathcal{L}-$ components.*

Proof.

Following [*Haissinsky & Tan, 2004*].

As $\{h_t\}$ is equicontinuous there exists a subsequence $\{h_{t_n}\} \rightrightarrows H$. By Proposition 4.4, $H$ maps each $\mathcal{L} - component$ to a point. We will show that they are the only nontrivial fibers of $H$.

Let $z \neq w$. such that $z, w$ are not in the same $\mathcal{L} - component$. We may assume that $H(z)$ and $H(w)$ are not fixed points under $h_t$ in the complement of $\mathcal{B}^*$, and that $\{z, w, a, b\}$ are in different $\mathcal{L} - components$. By Lemma 4.8 $\Lambda\left(h_t\left(\Gamma_{(z,w),(a,b)}\right)\right) \geq m > 0$ for all $t$. Thus by the contrapositive of Corollary B.2 mentioned above $d(h_t(z), h_t(w)) \geq m' > 0$ for all $t$ when $t \to 1$. Therefore $H(z) \neq H(w)$.

■

Therefore, sections 4.1, 4.2, 4.3, 4.4, and 4.5 prove Theorem B.

■

## 4.6 Equicontinuity of $\{f_t\}$

In the appendix of [*Haissinsky & Tan, 2004*] is proved the next lemma related to uniform convergence and was cited in Section 2.3:

LEMMA 2.2 *Let $f : \mathbb{C} \to \mathbb{C}$ be a continuous surjective map. For $t \in [0, 1)$, let $G_t, R_t : \overline{\mathbb{C}} \to \overline{\mathbb{C}}$ be two families of homeomorphisms of $\overline{\mathbb{C}}$. Assume that, as $t \to 1$, $G_t$ and $H_t$ converge uniformly to continuous maps $G, R$, respectively, and $f$ maps each fiber of $G$ into a fiber of $R$. Then $f_t := R_t \circ f \circ \left(G_t^{-1}\right) : \mathbb{C} \to \mathbb{C}$ converges uniformly to a continuous map $F$ and $F \circ G = R \circ f$.*

PROPOSITION 4.6 *The family $\{f_t\}$ is equicontinuous in $\overline{\mathbb{C}}$, where $f_t := h_t \circ f \circ \left(h_t^{-1}\right)$ is an entire transcendental function.*

Proof.

The proof is the same as the one given by [*Haissinsky & Tan, 2004*]. We know that $\{h_t\}$ is equicontinuous, then there exists a subsequence $\{h_{t_n}\} \rightrightarrows H$. Thus from proposition 4.5 each $\mathcal{L}-$ *component* is a fiber of $H$ and all the other fibers are single points. So, $f$ maps any fiber of $H$ into a fiber of $H$. Replacing both $R_t$ and $G_t$ by $h_{t_n}$ and using the same $f$ in the notation in Lemma 2.2, we conclude that $f_{t_n} = h_{t_n} \circ f \circ \left(h_{t_n}^{-1}\right)$ converges uniformly to a limit function $F$.

Finally, as $f$ is an entire transcendental function, by Lemma 2.1, so it is $f_t$.

■

Note that any limit $F$ of $\{f_t\}$ is also an entire function by the Weierstrass theorem on uniform convergence. This will be used in Section 4.8.

## 4.7 Semihyperbolicity of $F$

Now we will prove that any limit function $F$ of the sequence $\{f_t\}$ is semihyperbolic. Since the mappings $h_t$ are homeomorphisms of the sphere $\overline{\mathbb{C}}$, if $f$ is semihyperbolic, it follows that $f_t = h_t \circ f \circ \left(h_t^{-1}\right)$ is semihyperbolic. Let us normalize $f$ in such a way that 0 is a fixed point of $h_t$ in $F(f)$. Let $C := \inf_t dist_e\left(0, J\left(f_t\right)\right)$. Since $0 \in F(f)$ and is disjoint from $J(f) \cup \mathcal{L}$, and $\{h_t\}$ is equicontinuous, then $C > 0$. And we have the following Lemma:

LEMMA 4.10 *If $r < C/2$ then there is some constant $r' > 0$ such that, for all $z \in J(f) \setminus \overline{\mathcal{L}}$ and any $w \notin \mathcal{B}^* \cup D_r(z)$, we have that*

$$|H(z) - H(w)| \geq r'.$$

The proof is the same as in Lemma 2.11 by [*Haissinsky & Tan, 2004*].

PROPOSITION 4.7 *Let $F$ be a limit function of $\{f_t\}$ where $f_t = h_t \circ f \circ \left(h_t^{-1}\right)$ and $f$ is a semihyperbolic function such that $sing\left(f^{-1}\right) \cap \mathcal{L} = \emptyset$. Then $F$ is semihyperbolic.*

Proof.

There are some modifications to the proof of [*Haissinsky & Tan, 2004*] for this case.

As $f$ is semihyperbolic, we denote its constants as $(r_0, N)$.

Let $\beta \in \partial\lambda$ for some grand leaf $\lambda \in \mathcal{L}$, let $n_0 \in \mathbb{N}$ and let $\mathcal{B}^*\left(n_0\right)$ be the union of all the $\mathcal{B}^* -$ *components* such that any $\mathcal{B}^* \setminus \mathcal{B}^*\left(n_0\right) -$ *component* does not exceed $min\left\{r_0, r'_\beta\right\}/3$ where $r'_\beta$ is given by Definition 2.4.

Define

$$r_\beta := min\left\{\frac{r_0}{3}, d\left(\beta, \mathcal{B}^*\left(n_0\right) \setminus B^*\left(\beta\right)\right), \frac{r'_\beta}{3}\right\}.$$

Let us rename $B(\lambda) := B'(\lambda) \cup \Delta_\alpha \cup \Delta_\beta$ from Definition 2.4 as $B_\beta$. Set the annulus $A\left(\beta\right) := \Delta_\beta \setminus \Delta_{\beta_k}$, where $\Delta_{\beta_k}$ is the component of $f^{-k}\left(\Delta_\beta\right)$ contained in $\Delta_\beta$ by density in the inverse images of the non-exceptional points in the Julia set and shrinking due to Lemma 2.4 for some $k \in \mathbb{N}$. Then there exists $r_2 > 0$ such that for any $z \in J(f) \cap A\left(\beta\right)$ the disk $D_{r_2}(z)$ is disjoint from $\mathcal{B}^*\left(n_0\right)$. Note that $r_2 \leq |z - \beta| \leq r_\beta/2$.

Define

$$r < \left\{\frac{C}{2}, \frac{r_0}{3}, r_2\right\},$$

and consider a point $x \in J(f) \setminus \mathcal{L}$. By semihyperbolicity, for all $n \in \mathbb{N}$ and for all components $V$ of $f^{-n}\left(D_r(f^{-n}(x))\right)$, $deg\left(f^n : V(x) \to D_{r_0}\left(f^n(x)\right)\right) \leq N$, where $f^n(x) = a$ in the definition of semihyperbolicity. Now, we will construct inductively a sequence $\left\{n'_p\right\}$ such that $deg\left(F^{n'_p} : V'(H(x)) \to D_{r'}\left(F^{n_p}(H(x)\right)\right) \leq N$, where $r'$ is associated with $r$ by Lemma 4.10.

We assume that we have already constructed $n'_1, \ldots, n'_{p-1}$. Let $n$ be the smallest natural so that $n > n'_{p-1}$. We distinguish two cases. The first follows Lemma 4.10 when the point $x$ is far from $\mathcal{B}^* - components$ with large diameter, the second from Lemma 2.4.

Case 1, if $D_r\left(f^n(x)\right) \cap \mathcal{B}^*\left(n_0\right) = \emptyset$, then we define $D'_n$ to be the union of $D_r\left(f^n(x)\right)$ with all $\mathcal{B}^* - components$ $B^*$ such that $B^* \cap D_r\left(f^n(x)\right) \neq \emptyset$. Since $r < r_0/3$, it follows that $D'_n \subset D_{r_0}\left(f^n(x)\right)$ by the definition of $\mathcal{B}^*\left(n_0\right)$ above, so it is also the case for the fill-in of $D'_n$, that we denote as $D_n$. It follows from Lemma 4.10 that $H\left(D_n\right) \supset D_{r'}\left(H\left(f^{n_k}(x)\right)\right)$. Let $V(H(x))$ be the connected component of $F^{-n}\left(D_{r'}\left(F^n(H(x))\right)\right)$ which contains $H(x)$, therefore the degree is at most $N$ since $D_n \subset D_{r_0}\left(f^n(x)\right)$. Thus, we set $n'_p = n$.

Case 2, if $D_r\left(f^n(x)\right) \cap \mathcal{B}^*\left(n_0\right) \neq \emptyset$, then there is a $\mathcal{B}^*\left(n_0\right) - component$ $G^*$ such that $f^n(x) \in G^*$. Then, there is a minimal iterate $j$ such that $f^{n+j}(x) \in A(\beta)$. We set $n'_p = n + j$. Since $r < r_2$, it follows that $D_r\left(f^{n+j}(x)\right) \cap \mathcal{B}^*\left(n_0\right) = \emptyset$. Define $D'$ as the union of $D_r\left(f^{n+j}(x)\right)$ with all $\mathcal{B}^* - components$ $B^*$ such that $B^* \cap D_r\left(f^{n+j}(x)\right) \neq \emptyset$. Let $w \in D_r\left(f^{n+j}(x)\right)$. Then

$$\begin{aligned}
|w-\beta| \leq \quad & |w-f^{n+j}(x)| + |f^{n+j}(x)-\beta| \leq & \text{By Triangle inequality,}\\
& \leq \left(r + \frac{r'_\beta}{3}\right) + \frac{r_\beta}{2} \leq & \text{By the definition of } r_2,\\
& \leq \left(r_2 + \frac{r'_\beta}{3}\right) + \frac{r_\beta}{2} \leq & \text{By definition of } r,\\
& \leq \frac{r'_\beta}{6} + \frac{r'_\beta}{3} + \frac{r'_\beta}{6} & \text{Since } r_2 \leq \frac{r_\beta}{2} \leq \frac{r'_\beta}{6} \text{by Definition 2.4.}
\end{aligned}$$

It follows that $D' \subset \Delta_\beta$ and the fill-in of $D'$, denoted as $D$, is also contained in $\Delta_\beta$. Therefore $H\left(D_r\left(f^{n+j}(x)\right)\right) \supset D_{r'}\left(F^{n+j}(H(x)\right)$, and the degree of the restriction of $F^{n+j}$ to any component of $F^{-(n+j)}\left(D_{r'}\left(F^{n+j}(H(x)\right)\right)$ is at most $N$ since $f$ is semihyperbolic.

Now, since $f$ is semihyperbolic, it is a proper function in $V$, and $H$ is a proper map as a consequence of Proposition 4.5, $F = H \circ f \circ H^{-1}$, thus $F$ is a proper function in $H(V)$.

Therefore, $F$ is semihyperbolic.

■

## 4.8 Proof of Theorem C

To finish the proof of Theorem C, we will use a theorem on rigidity by [*Skorulski & Urbanski 2012*]. Let $X$ and $Y$ be two metric spaces, we say that a homeomorphism $H : X \to Y$ is *locally bi-Lipschitz* if each point $x \in X$ has some open neighborhood $U$ such that both the restriction $H\mid_U : U \to H(U)$ and its inverse $\left(H\mid_U\right)^{-1} : H(U) \to U$ are Lipschitz continuous. In particular, this theorem needs the concept of tame function:

DEFINITION 4.2 *An entire transcendental function $f : \mathbb{C} \to \mathbb{C}$ is called **tame** if its postsingular set $P(f) = \overline{\bigcup_{n=0}^{\infty} f^n\left(Sing\left(f^{-1}\right)\right)}$ does not contain the Julia set of $f$.*

However, we will not need to add this concept to the hypothesis of Theorem C due to the next lemma:

LEMMA 4.11 *If an entire transcendental function $f$ is semihyperbolic, then $f$ is tame.*

Proof.

We prove the contrapositive: if $f$ is not tame, then $f$ is not semihyperbolic.

As $f$ is not tame, $J(f) \subseteq P(f)$. Let $z \in J(f)$, thus $z \in P(f)$ and there exists $n = n(z) \in \mathbb{N}$ such that $z \in \overline{f^n\left(Sing\left(f^{-1}\right)\right)}$, so $\{f^{-n}(z)\} \subseteq \overline{Sing\left(f^{-1}\right)}$, where $\{f^{-n}(z)\} \subseteq J(f)$. Hence $J(f) \cap \overline{Sing\left(f^{-1}\right)} \neq \emptyset$, but $f$ is not semihyperbolic in $\overline{Sing\left(f^{-1}\right)}$, therefore $f$ is not semihyperbolic.

■

THEOREM 4.1 [*Skorulski & Urbanski 2012*] *If the restrictions to the Julia sets of two tame transcendental meromorphic functions $f : \mathbb{C} \to \overline{\mathbb{C}}$ and $g : \mathbb{C} \to \overline{\mathbb{C}}$ are topologically conjugate by a locally bi-Lipschitz homeomorphism $H : J(f) \to J(g)$, then this conjugacy extends to an affine linear map $(z \mapsto az + b)$conjugacy from $\mathbb{C}$ to $\mathbb{C}$ between the meromorphic functions $f, g$.*

On the other hand, we will need that the Lebesgue measure $m\left(J(f)\right)$ to be 0, to guarantee this we use a theorem by [*Stallard 1990*]. A set $E \subseteq \mathbb{C}$ is called *thin at* $\infty$ if there exist $R > 0$ and $\varepsilon > 0$ such that for any $z \in \mathbb{C}$ and any $D_r(z)$ with $r > R$,$density\left(E, D_r(z)\right) < 1 - \varepsilon$ in all sufficiently large discs, where the *density of* $E$ is defined by

$$density\left(E, D_r(z)\right) := \frac{m\left(E \cap D_r(z)\right)}{m\left(D_r(z)\right)}.$$

THEOREM 4.2 [*Stallard 1990*] *If $f$ is an entire transcendental function such that $d\left(P(f), J(f)\right) > 0$ and $E$ is a measurable completely invariant subset of $J(f)$ such that $E$ is thin at $\infty$, then $m(E) = 0$.*

Now we end the proof of one of the main results of this thesis.

THEOREM C *Let $f : \mathbb{C} \to \mathbb{C}$ be a semihyperbolic entire transcendental function with a Baker domain $U$ that satisfies property $P$ and with $J(f)$ thin at $\infty$. Let $\mathcal{L}$ be a grand orbit of a Baker lamination $\Lambda$ in $U$ that does not contain a leaf of the type $\lambda_{a,\infty}$ and $sing\left(f^{-1}\right) \cap \mathcal{L} = \emptyset$, then there exists a uniformly convergent continuous pinching deformation $f_t = h_t \circ f \circ h_t^{-1}$ to an entire function $F$. The mappings $h_t$ are quasiconformal mappings that converge uniformly to a map $H$, whose non-trivial fibers are the $\mathcal{L} -$ components.*

Proof.

Following [*Haissinsky & Tan, 2004*] with some significant modifications.

We know from Theorem B that $\{h_t\}$ is equicontinuous, and from Proposition 4.6 that so it is $\{f_t\}$. To prove the uniform convergence of the deformation, we will prove the uniqueness of the limits of the uniform convergent subsequences of $\{h_t\}$ as $t \to 1$.

Let us suppose that there are subsequences $\{t_n\}$ and $\{s_n\}$ from $t$ tending to 1 such that $\{h_{t_n}\} \rightrightarrows H_1$, $\{h_{s_n}\} \rightrightarrows H_2$, $\{f_{t_n}\} \rightrightarrows F_1$ and $\{f_{s_n}\} \rightrightarrows F_2$, where $\rightrightarrows$ denotes uniform convergence. By Proposition 4.5, $H_1$ and $H_2$ have the same fiber systems, then there exists a homeomorphism $\phi$ making the following diagram commute:

$$
\begin{array}{ccccc}
 & & Id & & \\
 & \overline{\mathbb{C}} & \to & \overline{\mathbb{C}} & \\
H_1 & \downarrow & & \downarrow & H_2 \\
 & \overline{\mathbb{C}} & \to & \overline{\mathbb{C}} & \\
 & & \phi & &
\end{array}
$$

Now, let us see that $H_i\left(J(f) \cup \mathcal{L}\right)$ is completely invariant by $F_i$. We see that

$$
\begin{array}{rcccc}
F_i\left(H_i\left(J(f) \cup \mathcal{L}\right)\right) & = & F_i\left(H_i\left(J(f)\right) \cup H_i\left(\mathcal{L}\right)\right) & = & F_i\left(H_i\left(J(f)\right)\right) \\
 & = & H_i\left(f(J(f))\right) & = & H_i(J(f)) \\
 & = & H_i\left(J(f) \cup \mathcal{L}\right), & &
\end{array}
$$

since $F_i \circ H_i = H_i \circ f$. Similarly $F_i^{-1}\left(H_i\left(J(f) \cup \mathcal{L}\right)\right) = H_i\left(J(f) \cup \mathcal{L}\right)$. Furthermore, as a consequence of the properties of the Julia set and that $F_i \circ H_i = H_i \circ f$, $H_i\left(J(f) \cup \mathcal{L}\right)$ has no isolated points, and periodic points under $F_i$ are dense in $H_i\left(J(f) \cup \mathcal{L}\right)$. Therefore $H_i\left(J(f) \cup \mathcal{L}\right) = J\left(F_i\right)$.

Since

$$
\begin{array}{ccccccc}
 & & & & F_1 & & \\
 & \overline{\mathbb{C}} & & \overline{\mathbb{C}} & \to & \overline{\mathbb{C}} & \\
 & \downarrow & H_1 & \uparrow & & \uparrow & H_1 \\
\phi & \downarrow & & \overline{\mathbb{C}} & \xrightarrow{f} & \overline{\mathbb{C}} & \\
 & \downarrow & H_2 & \downarrow & & \downarrow & H_2 \\
 & \overline{\mathbb{C}} & & \overline{\mathbb{C}} & \to & \overline{\mathbb{C}} & \\
 & & & & F_2 & &
\end{array},
$$

this implies that $\phi$ is a topological conjugacy from $F_1$ to $F_2$, and when we restrict $\phi$ to the Fatou set $\mathcal{F}(F_1)$, we have that $\phi = H_2 \circ Id \circ H_1^{-1} = H_2 \circ H_1^{-1}$.

Due to the Technical Assumption of Section 2.3 on the complex structure $\sigma_t$, on any compact set of $H_1^{-1}(\mathcal{F}(F_1))$ the two maps $H_1$ and $H_2$ integrate the same complex structure. Thus by Theorem 2.2 and the uniform convergence on $\{h_{t_n}\}$ and $\{h_{s_n}\}$, $H_1$ and $H_2$ are the same map up to a composition with an automorphism of $\overline{\mathbb{C}}$ (an affine map $z \mapsto a'z + b'$), i.e., the map $\phi$ is conformal in $\mathcal{F}(F_1)$ by Lemma 2.1 and $H_2 = (a'z + b') \circ H_1$, thus $\phi\mid_{\mathcal{F}(F_1)} = a'z + b'$.

By construction of $P_t$ and $\psi$ in Section 2.3, these maps are locally bi-Lipshitz in their domains and so its integrating map $h_t$ is in $J(f_t)$, then $H_1$ and $H_2$ are locally bi-Lipschitz in $J(F_1)$ and $J(F_2)$, respectively, so it is $\phi$ in $J(F_1)$.

Applying Theorem 4.1, $\phi(z) = c'z + d'$ for all $z \in \mathbb{C}/J(F_1)$, and by holomorphic rigidity, $a' = c', b' = d'$ in $\mathcal{F}(F_1)$. But $h_t$ fixes $a, b, c$, due to normalization and equicontinuity of $\{h_t\}$, thus $H_1, H_2$, and $\phi$ fix the same points. Since $\phi$ is a Moebius transformation, $\phi = Id$ on $\mathbb{C} \setminus F_1$. Now, we need to extend $\phi$ to $\mathbb{C}$.

Since $f$ is semihyperbolic, by Proposition 4.7 $F_1$ is semihyperbolic, which implies $d(P(F_1), J(F_1)) > 0$. Now, we prove that $F_1$ is thin at $\infty$.

Due to $H_1$ being conformal and injective in $\mathbb{C} \setminus \mathcal{V}$, the derivative of $H_1$ is non-zero. Thus the one-quarter Koebe theorem prevents the components of $\mathcal{U} \setminus \mathcal{V}$ to shrink excessively. Since $m(\mathcal{L}) = 0$, we have that $J(f)$ and $J(F_1)$ have the same density for sufficiently large discs $D_r(z)$. Then $J(F_1)$ is thin at $\infty$.

Therefore by Theorem 4.2, $m(J(F_1)) = 0$ and $\phi = Id$ on all $\mathbb{C}$ except in a set of Lebesgue measure zero, implying that $\phi$ is holomorphic in all $\mathbb{C}$, particularly in $J(F_1)$, and $\phi(z) = z$ for any $z \in \mathbb{C}$.

Therefore $H_1 = H_2$, and $F_1 = F_2$. And this proves Theorem C.

■

# 5 Wandering Domain to a positive distance from $P(f)$

One of the main reasons to pursue Theorem C is to solve an open problem in the area. As commented in [*Bergweiler 1995*] and [*Lauber 2004*], he asks if there exists an entire transcendental function with a wandering domain at a positive distance from its postsingular set $P(f)$.

Until the present time, we are not aware if this has been solved by someone else, but we answer in positive the problem using the Corollary D of Theorem C. We will show the existence of such a function in the next situation.

In [*Bergweiler 1995*] the function $f(z) = 2 - log(2) + 2z - exp(z)$ is introduced as an example of a function with Baker domain to a positive distance from $Sing\,(f^{-1})$. This function is the lifting of the function $g(w) = \frac{1}{2}w^2 e^{2-w}$ with $w = exp(z)$. The basin of attraction $\mathcal{A}'$ of the super attracting fixed point 0 of $g$ is lifted to a Baker domain of hyperbolic type-I, $U \supseteq \{z : Re\, z < -2\}$, lifting 0 to $\infty$. Also, $f(z)$ has a superattracting fixed point at $log(2)$ and a repulsing fixed point at $x_0 \lessapprox -0.900477$ in the boundary of the basin of attraction $\mathcal{A}$ of $ln(2)$. See Figure 5.1.

As commented in [*Bergweiler 1995*], there are no finite asymptotic values in $f(z)$. Hence, $Sing\,(f^{-1}) = \{log(2) + 2\pi ki \mid k \in \mathbb{Z}\}$ and since

$$Sing\left(f^{-1}\right) = P(f) = \{log(2) + 2\pi ki \mid k \in \mathbb{Z}\}.$$

In addition $P(f)$ is contained in a wandering domain $W$, because, if $\mathcal{A}$ is the basin of attraction of $log(2)$, then

$$W = \{z + 2\pi ki \mid z \in \mathcal{A}, k \in \mathbb{Z} \setminus \{0\}\}.$$

This implies that $d\,(P(f), U) > 0$ and that $f$ is semihyperbolic. Also, $J(f)$ is thin at $\infty$ as a direct consequence of the existence of its Baker domain.

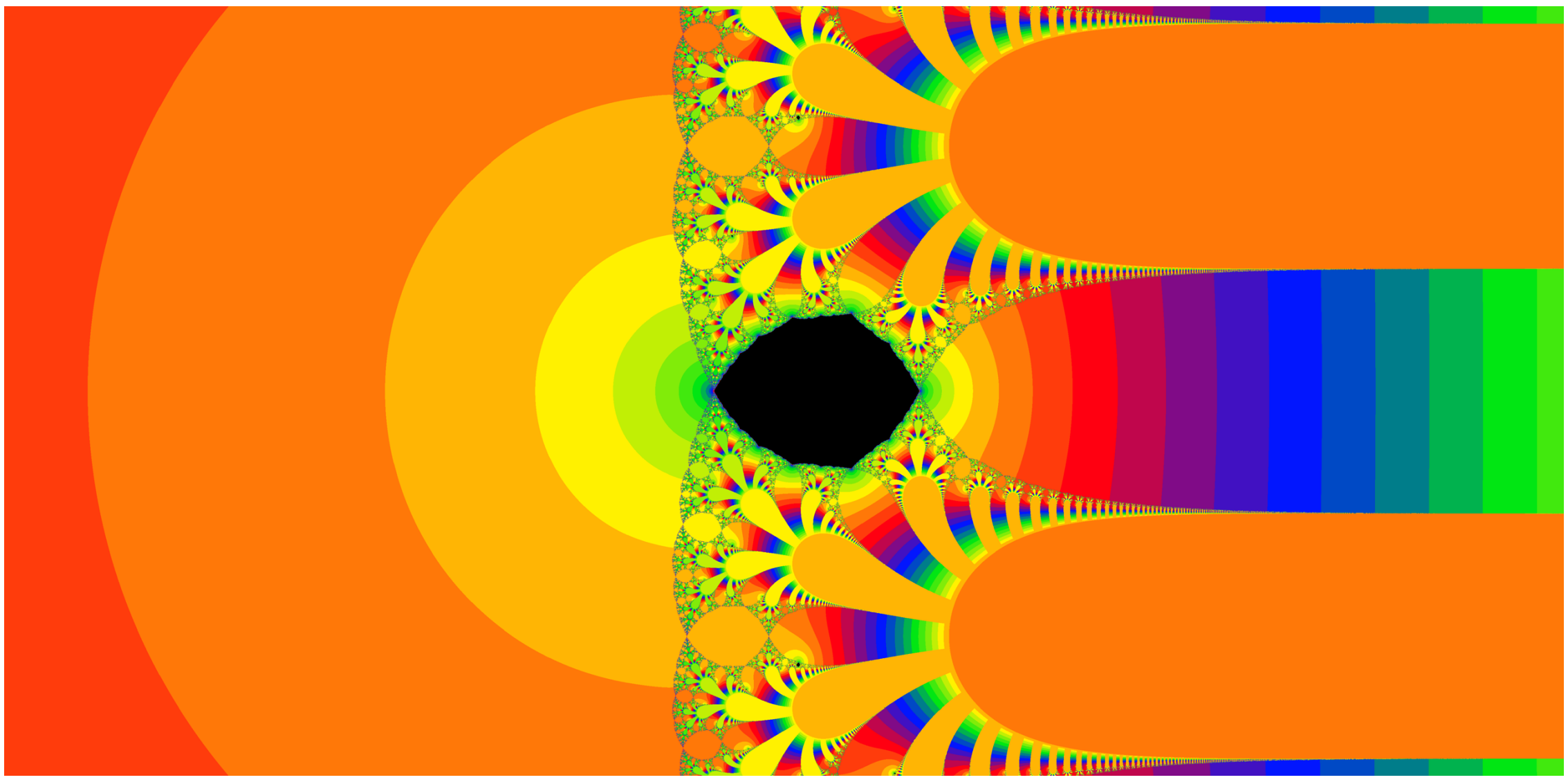

**Figure 5.1:** Dynamics of the function $f(z) = 2 - log(2) + 2z - exp(z)$.

These are some of the properties of $f(z)$, and we will use this function to prove the next theorem.

THEOREM D *There exists an entire function $F$ with a wandering domain $W$ such that $d\,(P(F), W) > 0$.*

Proof.

We will do a pinching deformation to $f(z) = 2 - log(2) + 2z - exp(z)$ to prove the existence of the function $F$. First, let us apply Theorem C to this function, by the remarks above, it only rests to see that $f(z)$ is semihyperbolic and to give an adequate grand orbit $\mathcal{L}$.

Let $x_0 \lessapprox -0.900477$ be the repulsing fixed point of $f$ mentioned above, we see that

$$f\,(x_0 + 2\pi ki) = x_0 + 2^2\pi ki,$$

with $k \in \mathbb{Z}$. Since $f(z) \approx 2 - log(2) + 2z$ when $Re(z) \to -\infty$, by the uniformization of univalent Baker domains in Theorem 1.2, there exists a topological semi-annulus $R$ as a fundamental domain of $f$ in the Baker domain $U$ with vertices in $\{x_0 \pm 2\pi i, x_0 \pm 4\pi i\}$. We choose two geodesics $\lambda_+$ and $\lambda_-$ in $R$ such as their end-points are in the sides of $\partial R$, which goes from $x_0 + 2\pi i$ to $x_0 + 4\pi i$, and from $x_0 - 2\pi i$ to $x_0 - 4\pi i$, respectively. Now we take the grand orbit of $\lambda_+$ and $\lambda_-$ under $f$ to build a grand orbit $\mathcal{L}$ under iteration as described in section 2.2. See Figure 5.2.

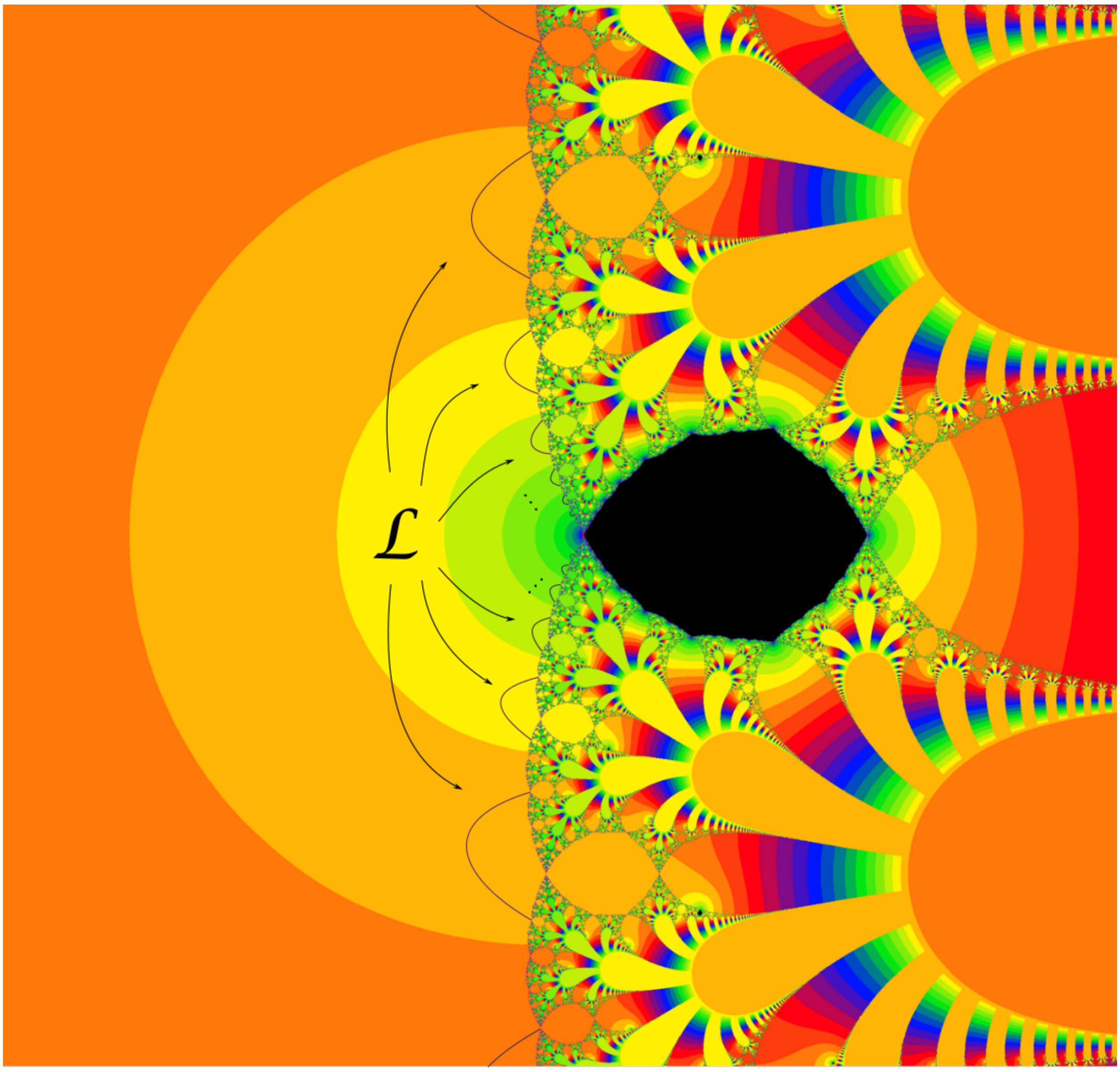

**Figure 5.2:** The grand orbit $\mathcal{L}$ of $f(z) = 2 - log(2) + 2z - exp(z)$.

By Theorem C, we can pinch $\mathcal{L}$ via quasiconformal conjugations $h_t$ of $f$ where $f_t = h_t \circ f \circ h_t^{-1}$ uniformly converges to a function $F$ when $t \to 1$. For every $t \in [0, 1)$ we have a Baker domain $U_t = h_t(U)$ for $f_t$. When $t = 1$, the grand leaves of $\mathcal{L}$ collapse in one point creating a new set $h_1(U) = U_1 \cup W_1$ where $U_1 \cap W_1 = \emptyset$. Here, $U_1$ is a Baker domain of hyperbolic type-I of $F = f_1$. Let us see what happens with $W_1$.

By Theorem 2 of [*Bergweiler 1995*], the boundary of the Baker domain $U$ of $f$ is a Jordan curve in $\overline{\mathbb{C}}$. Let us denote by $J_{\lambda_+}(f)$ the subset of $J(f)$ going from $x_0 + 2\pi i$ to $x_0 + 4\pi i$, and by $J_{\lambda_-}(f)$ the subset of $J(f)$ going from $x_0 - 2\pi i$ to $x_0 - 4\pi i$, both in the sides of the semiannulus $R$, thus $\lambda_+ \cup J_{\lambda_+}(f)$ and $\lambda_- \cup J_{\lambda_-}(f)$ are Jordan curves in $\overline{R}$. Since $h_t$ is an homeomorphism, $h_t\left(\lambda_+ \cup J_{\lambda_+}(f)\right)$ and $h_t\left(\lambda_- \cup J_{\lambda_-}(f)\right)$ are Jordan curves for every $t \in [0,1]$.

Especially $h_1\left(\lambda_+ \cup J_{\lambda_+}(f)\right) = h_1\left(J_{\lambda_+}(f)\right)$ and $h_1\left(\lambda_- \cup J_{\lambda_-}(f)\right) = h_1\left(J_{\lambda_-}(f)\right)$ are Jordan curves whose bounded components in their complements are topological discs and the equality is because $\lambda_+$ and $\lambda_-$ are in the fibers of the pinching deformation by Theorem C. Hence, by construction, $W_1$ is made by disjoint discs who share one point in their boundary with $\partial U_1$.

But $W_1$ has inherited the dynamical behavior of $f$ through the fundamental domain, so every topological disc in $W_1$ moves to the "next" topological disc under iteration of $F$, as the annuli did in the dynamics of $f(z)$, thus $W_1$ is a wandering domain. Furthermore, $P(F)$ remains in $h_1(W)$, the wandering domain of $F$ which is topologically the same wandering domain $W$ of $f(z)$, since the pinching deformation only affects the Baker domain $U$ of $f$ and there are no fibers of $h_1$ of the lamination contained in it.

Therefore, $F$ is an entire transcendental function with one univalent Baker domain $U_1$, a superattracting basin $h_1(\mathcal{A})$, and two wandering domains $W_1$ and $h_1(W)$ such that $d(W_1, P(F)) > 0$. See Figure 5.3.

■

It is worth mentioning that Theorem D is in accordance with corollary C by [*Baransky et al., 2020*] which states:

COROLLARY C [*Baransky et al., 2020*]

*Let $f$ be a topologically hyperbolic meromorphic map* (i.e., *dist*$(P(f), J(f) \cap \mathbb{C}) > 0$) *and $U$ be a Fatou component of $f$. Denote by $U_n$ the Fatou component such that $f^n(U) \subset U_n$ and suppose that $U_n \cap P(f) = \emptyset$ for $n > 0$. Then, for every compact set $K \subset U$, every $z \in K$ and every $r > 0$ there exists $n_0$ such that for all $n \geq n_0$, $D_r(f^n(z)) \subset U_n$. In particular, diam$(U_n) \to \infty$ and dist$(f^n(z), \partial U_n) \to \infty$, for every $z \in U$, as $n \to \infty$.*

In particular, the wandering domain $W_1$ of $F$ satisfies this corollary, because the pinching deformation only deforms the grand orbit $\mathcal{L}$ to point components, the neighborhood $\mathcal{Y} \setminus \mathcal{L}$ is bounded deformed, and $\mathbb{C} \setminus \mathcal{Y}^*$ is not deformed. This implies that the Baker domain $U_1$ of $F$ is practically the same as the Baker domain $U$ for $f$, which both are hyperbolic type I. Thus their fundamental domains grow conjugated to $z \mapsto az$, where in the case of $f$, $a = 2$. As the leaves of the lamination $\mathcal{L}$ are in the fundamental domains of $U$, when the pinching deformation is realized the wandering domain $W_1$ inherits this growth, satisfying corollary C from [*Baransky et al., 2020*].

To the knowledge of the author, the known examples of wandering domains do not grow without limit. It is relevant to ask if the return of this Corollary C by [*Baransky et al., 2020*] is valid. In the sense that the set of singularities of the inverse function, $sing\,(f^{-1})$, and the postsingular set $P(f)$ are related with many considerable results in holomorphic dynamics, for instance, in any attracting or parabolic basin, there exists a critical point.

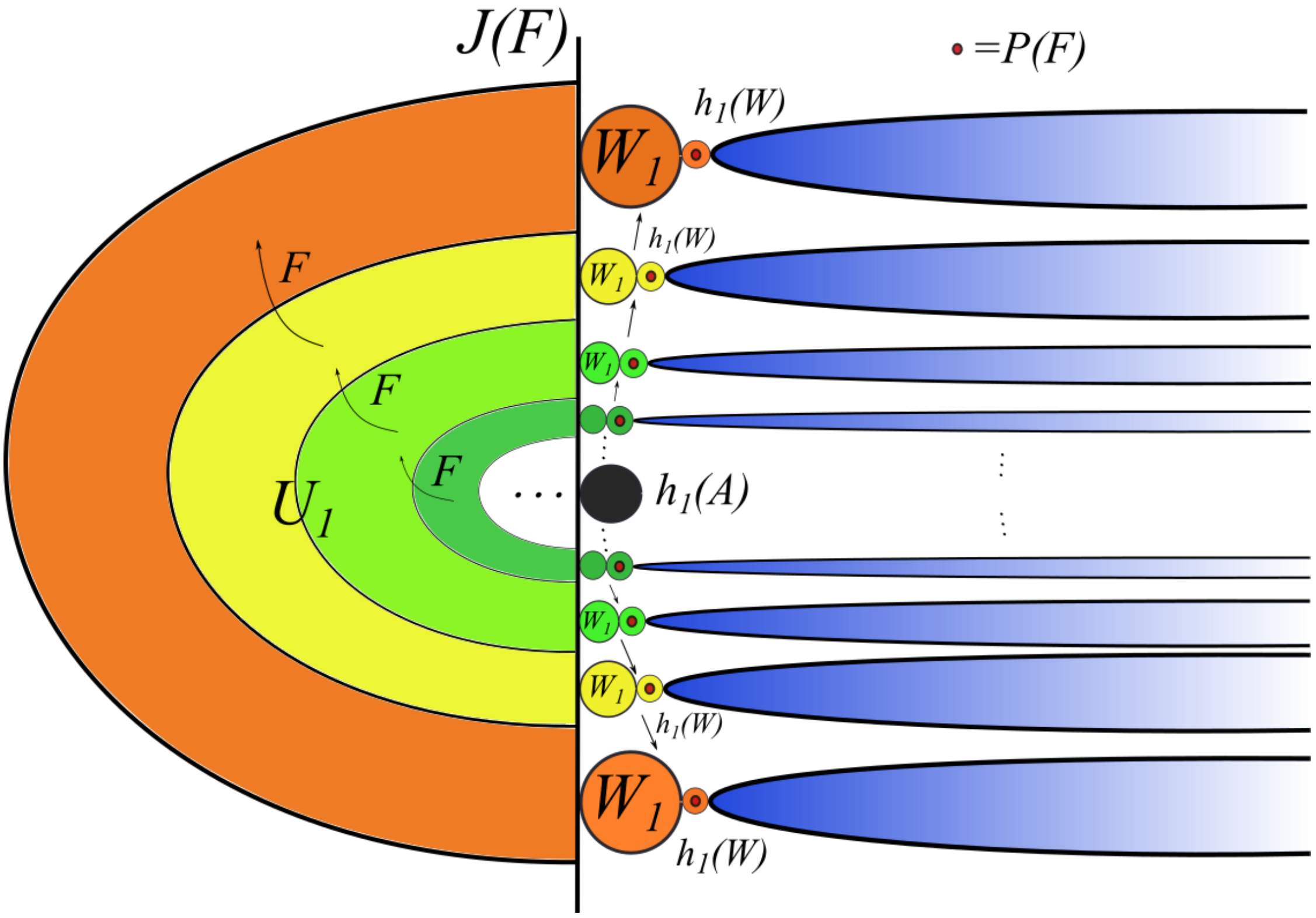

**Figure 5.3:** The wandering domain $W_1$ to a positive distance of $P(F)$. With Baker domain $U_1$ and its preimages in blue.

# References

[*Ahlfors 1966*] L. Ahlfors, "Lectures on Quasiconformal Mappings", Wadsworth & Brooks / Cole mathematics series, 1966.

[*Alexander 1994*] D. S. Alexander, "A History of Complex Dynamics: From Schroder to Fatou and Julia", Vieweg, 1994.

[*Astala & Martin 2001*] K. Astala and G. J. Martin, *Holomorphic Motions,* Papers on Analysis, Report Univ. Jyväskylä **83**, 27-40, 2001.

[*Baker 1975*] I. N. Baker, *The Domains of Normality of an Entire Function*, Ann. Acad. Sci. Finn. Ser. A I Math. **12**, 277-283, 1975.

[*Baker 1988*] I. N. Baker, *Infinite limits in the iterations of entire functions*, Ergodic Theory Dyn. Syst. **8**, 503-507, 1988.

[*Baker & Domínguez 1999*] I. N. Baker and P. Domínguez, *Boundaries of unbounded Fatou components of entire functions*, Ann. Acad. Sci. Fenn. **24**(2), 437-464, 1999.

[*Baker & Weinrich 1991*] I.N. Baker and J. Weinrich, *Boundaries which arise in the iteration of entire functions*, I.N. Baker and J. Weinrich, Rev. Roumaine Math. Pures Appl. 36, pp413-420, 1991.

[*Baransky & Fagella 2001*] K. Baransky and N. Fagella, *Univalent Baker Domains*, Nonlinearity **14**, pp411-429, 2001.

[*Baransky et al., 2020*] K. Baransky, N. Fagella, X. Jarque and B. Karpinska, *Fatou components and singularities of meromorphics functions.* Proceedings of the Royal Society of Edinburgh Section A, **150**(2), 633-654, 2020.

[*Bergweiler 1993*] W. Bergweiler, *Iteration of Meromorphic Functions*, Bulletin (New Series) of the American Mathematical Society, Vol. **29**, number 2, 151-188, 1993.

[*Bergweiler 1995*] W. Bergweiler, *Invariant Domains and Singularities*, Mathematical Proceedings of the Cambridge Philosophical Society, Vol **117**, Issue 03, pp525-532, 05/1995.

[*Bergweiler & Eremenko 2007*] W. Bergweiler and A. E. Eremenko, *Direct Singularities and Completely Invariant Domains of Entire Functions*, Illinois Jour. of Maths. **52**, 243-259, 2007.

[*Bergweiler & Morosawa 2002*] W. Bergweiler and S. Morosawa, *Semihyperbolic Entire Functions*, IOP Publishing Ltd & London Mathematical Society, Nonlinearity, Volume **15**, Number 5, 1673-1684, 2002.

[*Bergweiler & Zheng 2012*] W. Bergweiler and J. Zheng, *Some examples of Baker Domains*, IOP Publishing Ltd & London Mathematical Society, Nonlinearity, Volume **25**, Number 4, 1033-1045, 2012.

[*Bernstein & Chivian 2008*] A. Bernstein and E. Chivian "Sustaining Life, How Human Health Depends on Biodiversity", Oxford University Press, 2008.

[*Bers 1962*] L. Bers, *The equivalence of two definitions of quasiconformal mappings*, Comment. Math. Helv. **37**, 148-154, 1962/63.

[*Bers 1974*] L. Bers, *On spaces of Riemann surfaces with nodes,* Bulletin of American Society, Vol **80**, No. 6, 1219-1222, 1974.

[*Branner & Fagella 2014*] B. Branner and N. Fagella, "Quasiconformal Surgery in Holomorphic Dynamics", Cambridge University Press, Cambridge studies in advanced mathematics **141**, 2014.

[*Cayley 1879*] A. Cayley, *Applications of the Newton-Fourier Method to an Imaginary Root of an Equation*, Quarterly Journal of Pure and Applied Mathematics, **16**, 179-185, 1789.

[*Cremer 1932*] H. Cremer, *Uber die Schrodersche Funktionalgleichung und das Schwarsche Eckenabbildungsproblem*, Ber. Verh. Sachs. Akad. Wiss. Leipzig, Math-Phys. Kl. **84**, 291-324, 1932.

[*Devaney 1989*] R. L. Devaney, "An Introduction to Chaotic Dynamical Systems" 2nd edition, Addison Wesley Publishing Company, Inc., 1989.

[*Devaney 2010*] R. L. Devaney, *Complex Exponential Dynamics* in "Handbook of Dynamical Systems", Broer et al. eds. Elsevier, 125-223, 2010.

[*Dominguez & Sienra 2015*] P. Dominguez and G. Sienra, *Some pinching deformations of the Fatou function*, Fundamenta Mathematicae **228**(1), 1-15, 2015.

[*Douady & Hubbard 1984-85*] A. Douady and J. H. Hubbard, "Etude Dynamique des Polynomes Complexes I and II," Publ. Math Orsay, 1984-1985.

[*Eremenko & Lyubich 1992*] A. E. Eremenko and M. Yu. Lyubich, *Dynamical Properties of Some Classes of Entire Functions*, Annales de l´Institut Fourier, tome **42**, no. 4, 989-1021, 1992.

[*Fagella & Henriksen 2006*] N. Fagella and C. Henriksen, *The Teichmuller Space of an Entire Function*, Complex Dynamics: Families and Friends, AK Peters, pp450-456, 2008.

[*Fatou 1919*] P. Fatou, *Sur les equationes fonctionnelles*, Bulletin Soc. Math. France tome **47**, 161-271, 1919, and tome 48 , 33-94, 208-314,1920.

[*Gardiner 2000*] F. Gardiner, "Quasiconformal Teichmuller Theory", American Mathematical Society, Mathematical Surveys and Monographs **76**, 2000.

[*Haissinsky 2002*] P. Haissinsky, *Pincement de Polynomes*, Comment. Math. Helv. **77**, pp1-23, 2002.

[*Haissinsky & Tan, 2004*] P. Haissinsky and L. Tan, *Convergence of Pinching Deformations and Matings of Geometrically Finite Polynomials*, Fundamenta Mathematicae **181**, no.2, pp143-188, 2004.

[*Julia 1918*] G. Julia, *Memoire sur l´Itération des Fonctions Rationnelle*s, J. Math. Pures Appl. **8** no.4, 47-245, 1918.

[*Lauber 2004*] A. Lauber, *On the Stability of Julia Sets of Functions having Baker Domains*, Ph.D. Thesis, University at Gottingen, 2004.

[*Lehto 1987*] O. Lehto, "Univalent Functions and Teichmuller Spaces", Springer Verlag, 1987.

[*Lehto & Virtanen 1973*] O. Lehto and K. I. Virtanen "Quasiconformal Mappings in the plane", Springer Verlag, 1973.

[*Lyubich & Minsky 1997*] M. Lyubich and Y. Minsky, *Laminations in holomorphic dynamics*, J. Differential Geom. 47, 17-94, 1997.

[*Makienko 2000*] P. Makienko, *Unbounded Components in Parameter Space of Rational Maps*, Conformal Geometry and Dynamics Vol 4, 1-21, 2000.

[*Mandelbrot 1982*] B. Mandelbrot, "The Fractal Geometry of Nature", W.H. Freeman and Co., 1982.

[*Mañe 1993*] R. Mañe, *On a Theorem of Fatou*, Bol. Soc. Bras. Mat., Vol. 24, **1**, 1-11, 1993.

[*Mc Mullen 1994*] C. Mc Mullen, "Complex dynamics and renormalization", Princeton University Press, 1994.

[*Montel 1927*] P. Montel, "Familles normale," reprint, Chelsea Publishing Company, New York, 1974.

[*Morosawa et al. 1998*] S. Morosawa, Y. Nishimura, M. Taniguchi, and T. Ueda, "Holomorphic Dynamics", Cambridge University Press, 1998.

[*Rees 1986*] M. Rees, *Realization of matings of polynomials as rational maps of degree two*, unpublished manuscript, 1986.

[*Rippon & Stallard 1999*] P. J. Rippon and G. M. Stallard, *Families of Baker domains I*, Nonlinearity, Volume 12, Number **4**, pp1005, 1999.

[*Robles & Sienra 2022*] R. Robles and G. *Sienra, Baker Domains and Non Convergent Deformations,* Fractal Geom. 9, no. 1/2, pp. 1–23, 2022.

[*Schroeder 1870*] E. Shroeder, *Ueber unendlich viele Algorithmen zur Auflosung de Gleichungen*, Mathematsiche Annalen, **2**, 317-365, 1870.

[*Schroeder 187*1] E. Shroeder, *Ueber iterite Functionen*, Mathematische Annalen, **3**, 296-322, 1871.

[*Shishikura 2000*] M. Shishikura, *On a theorem of M. Rees for matings of polynomials.* The Mandelbrot set, Theme and Variations, 289-305, ed. Tan Lei, LMS Lecture Note Series 274, Cambridge University Press, 2000.

[*Shishikura & Tan 2000*] M. Shishikura and Tan L., *An Alternative Proof of Mañe´s Theorem on Non-expanding Julia sets,* The Mandelbrot set, Theme and Variations, 265-269, ed. Tan Lei, LMS Lecture Note Series 274, Cambridge University Press, 2000.

[*Siegel 1942*] C. L. Siegel, *Iteration of Analytic Functions*, Ann. of Math. **43**, 607-612, 1942.

[*Sienra 2006*] G. Sienra, Surgery and Hyperbolic Univalent Baker Domains. Nonlinearity 19, pp.265-279, 2006.

[*Skorulski & Urbanski 2012*] B. Skorulski and M. Urbanski, *Dynamical Rigidity of*

*Transcendental Meromorphic Functions*, Nonlinearity Volume **25** Number 8, 2337-2348, 2012.

[*Stallard 1990*] G. Stallard, *Entire Functions with Julia sets of zero measure*, Math. Proc. Cambridge Philos. Soc., Vol **108**, issue 3, 551-557, 1990.

[*Sullivan 1985*] D. Sullivan, *Quasiconformal homeomorphisms and dynamics I: Solution of the Fatou-Julia problem on wandering domains*, Ann. Math. **122**, 401-418.

[*Tan 1990*] L. Tan, *Matings of quadratic polynomials*, Ergodic Th. Dyn. Sys., **12**, 589-620, 1990.

[*Tan 2002*] L. Tan, *On Pinching Deformations of Rational Maps*, Annales Scientifiques de l´Ecole Normale Supérieure, 4e serie, **35**, 353-370, 2002.

[*Wilson 1999*] E. O. Wilson, "The Diversity of Life," WW Norton & Company, 1997.

[*Zakeri & Zeinalian 1996*] S. Zakeri and M. Zeinalian, *When Ellipses look like Circles: The Measurable Riemann Mapping Theorem*, Nashr-i Riazi **8**, 5-14, 1996.